\documentclass{article}
\usepackage{epsfig,epic,eepic,latexsym, amssymb}

\input xy
\xyoption{all}
\xyoption{2cell}
\UseAllTwocells
\xyoption{all}
\makeindex

\def \pf{PROOF:}
\def    \QED    {\hfill\hbox{\hskip 4pt
                \vrule width 5pt height 6pt depth 1.5pt}}
\def\epf{\QED\\}

\newcommand{\Cat}{ {\bf Cat}}

\newcommand{\SMC}{SMC}

\newcommand{\A}{ {\cal A} }
\newcommand{\B}{ {\cal B} } 
\newcommand{\C}{ {\cal C} }
\newcommand{\D}{ {\cal D} }
\newcommand{\X}{ {\cal X}}

\newcommand{\unc} { {\cal I} }

\newcommand{\Gra}{ {\cal H} }
\newcommand{\FGd}{ {\cal FG}} 

\newcommand{\Obj}{ Obj }

\newcommand{\gen}{\star}

\newcommand{\ac}{ ass } 
\newcommand{\lc}{ l }
\newcommand{\rc}{ r }

\newcommand{\syc}{ s } 
\newcommand{\un}{I}

\newcommand{\Res}{Rn}
\newcommand{\Ext}{En}

\newcommand{\Post}{Post}

\newcommand{\Dual}{D}
\newcommand{\ev}{ev}

\newcommand{\vv}{v}
\newcommand{\vs}{v^*}

\newcommand{\deuxnateva}{11.2}
\newcommand{\waxiomdouze}{19.10}
\newcommand{\tbctrois}{10.9}
\newcommand{\evdual}{11.9}
\newcommand{\homevtrois}{11.3}
\newcommand{\homevtroisdual}{11.4}
\newcommand{\calcdeux}{19.6}
\newcommand{\lemsymofundeux}{7.1}
\newcommand{\impronatPosttrois}{9.11}
\newcommand{\impronatPost}{9.10}
\newcommand{\lemvvun}{18.5}
\newcommand{\lemvvdeux}{18.8}
\newcommand{\vvgen}{18.4}
\newcommand{\tbcun}{10.8}
\newcommand{\evhomun}{11.6}
\newcommand{\evhomdeux}{11.8}

\newcommand{\BS}{SPC}
\newcommand{\StrBS}{StrSPC}
\newcommand{\tens}{\otimes}

\newcommand{\jj}{j}
\newcommand{\inv}{inv}
\newcommand{\eqi}{\bullet}

\newcommand{\unitb}{{\it j}}
\newcommand{\compb}{{\it c}}
\newcommand{\assb}{{\it \alpha}}
\newcommand{\rb}{{\it \rho}}
\newcommand{\lb}{{\it \lambda}}

\newcommand{\unit}{u}
\newcommand{\comp}{c}
\newcommand{\assp}{\alpha}
\newcommand{\rp}{\rho}
\newcommand{\lp}{\lambda}
\newcommand{\comh}{c'} 
\newcommand{\assh}{\alpha'}
\newcommand{\rh}{\rho'}
\newcommand{\lh}{\lambda'}

\newcommand{\rhp}{{r'}}
\newcommand{\lhp}{{l'}}

\newcommand{\TCcompu}{{c^1}}
\newcommand{\TCcompd}{{c^2}}
\newcommand{\TCcompt}{{c^3}}
\newcommand{\rhpp}{{u^1}}
\newcommand{\lhpp}{{u^2}}

\newcommand{\one}{1}

\newcommand{\unm}{1}
\newcommand{\assmu}{\tilde{\alpha}}
\newcommand{\rmu}{\tilde{\rho}}
\newcommand{\lmu}{\tilde{\lambda}}
\newcommand{\Hadd}{H_{+}}
\newcommand{\Hmu}{H_{\times}}

\newcommand{\mult}{\varphi}
\newcommand{\assmd}{\beta}
\newcommand{\lmd}{\gamma}
\newcommand{\assm}{\beta''}

\newcommand{\multh}{\varphi'}
\newcommand{\assmh}{\beta'}
\newcommand{\lmh}{\gamma''}
\newcommand{\lmz}{\gamma'}

\newcommand{\multl}{\underline{\varphi}}
\newcommand{\assml}{\underline{\beta}}
\newcommand{\lml}{\underline{\gamma}}

\newcommand{\M}{{\cal M}}
\newcommand{\N}{{\cal N}}
\newcommand{\mor}{H}
\newcommand{\deltal}{\underline{\delta}}
\newcommand{\deltah}{\delta'}
\newcommand{\deltad}{\delta}

\newtheorem{theorem}{Theorem}[section]
\newtheorem{definition}[theorem]{Definition}
\newtheorem{proposition}[theorem]{Proposition}
\newtheorem{lemma}[theorem]{Lemma}

\newtheorem{remark}[theorem]{Remark}

\newtheorem{tag}[theorem]{}

\title{Enrichments over symmetric Picard categories}
\author{
Vincent Schmitt\\
{\small email: vrrschmitt@yahoo.fr}
}

\begin{document}
\maketitle                      

\begin{abstract}
Categorical rings were introduced 
in \cite{JiPi07}, which we call 2-rings.
In these notes we present basic definitions
and results regarding 2-modules.
This is work in progress.
\end{abstract}

\begin{section}{Introduction}
Categorical rings were introduced by M.Jibladze and T.Pirashvili
in \cite{JiPi07}. We call those 2-rings. The present set
of notes contains basic results about 2-modules and this 
work is in progress.\\
Section \ref{prel} contains technical preliminaries,
namely reminders on symmetric Picard categories as well
convenient references to previous works. 
In particular the developments in \cite{Sch08} for symmetric 
monoidal categories transpose  well to symmetric Picard categories 
and a suitable tensor product can be defined for the latter.
Section \ref{BSCat} and \ref{BSCattens}
treat enrichments over symmetric Picard categories.
Those were introduced in \cite{Dup08} and defined 
by means of multilinear maps. We also define
them using the tensor product.
Section \ref{Exples} contains expected examples 
of 2-rings and 2-enrichments.
Section \ref{Mod} contains basic results regarding 
categories of $\A$-modules for a 2-ring $\A$.
In particular we show that $\A$-modules are
particular algebras for the endo-2-functor 
$\A \otimes -$ of the 2-category of symmetric 
Picard categories.
A large appendix contains the more technical developments.   
\end{section}

\begin{section}{Preliminaries}\label{prel}
A categorical group structure $(\A,\jj)$ consists
of a monoidal category $\A$
and an assignment for every object $a$ of $\A$ of
an object $a^{\eqi}$ with
an isomorphism $\jj_a: \unc \rightarrow a^{\eqi} \otimes a$,
($a^{\eqi}$ is an {\em inverse} of $a$).
We are concerned in this paper with {\em symmetric 
Picard categories} which are the categorical groups $(\A,\jj)$ 
for which $\A$ has a symmetric monoidal structure
and its underlying category is a groupoid.
$\BS$ denotes the 2-category with objects symmetric Picard categories,
arrows symmetric monoidal functors and 2-cells
monoidal natural transformations. There is a
forgetful 2-functor $\BS \rightarrow \SMC$
forgetting the group structure
where $\SMC$ denotes the 2-category with 
objects symmetric monoidal categories, 
arrows symmetric monoidal functors,
and monoidal natural transformations as 2-cells. 
The 2-categorical properties of $\BS$ are similar 
to those of $\SMC$, the latter 2-category has been studied
in different works in particular
in \cite{HyPo02} and in \cite{Sch08}. We refer the reader to
this last work for basic notations, and more elaborate results. 
In this first section, we describe briefly and compare the important 
properties of $\SMC$ and $\BS$.\\

The 2-category $\SMC$ admits an internal hom and
the same holds for $\BS$ which hom is inherited from $\SMC$.
The following was mentioned
in \cite{Dup08} with a rather concise explanation.
\begin{lemma}
Given any two objects in $\A$ and $\B$ in $\SMC$ 
with $\B$ being a symmetric Picard category, the internal 
hom $[\A,\B]$ in $\SMC$ admits a symmetric Picard structure
given pointwise by that of $\B$.
\end{lemma}
We present a proof that relies on coherence results for categorical 
groups from \cite{Lap83}.
(This work also contains references to earlier works on the topic.)
Let us recall these. For a categorical group $\A$
with family of isomorphisms $\jj_a: \un \rightarrow a^{\eqi} \otimes a$
for each object $a$ there is a unique way of extending the assignment 
$a \mapsto a^{\eqi}$ into a functor $\A^{op} \rightarrow \A$
that makes the $\jj_a$ natural in $a$.
It is an equivalence. 
We will write it ${(-)}^{\eqi}$ and 
write $f^{\eqi}: b^{\eqi} \rightarrow a^{\eqi}$ for the image
of any arrow $f: a \rightarrow b$ by this functor.
There is a coherence theorem stating 
that any pair of $a$ and $b$ of objects of $\A$  
there is at most one ``canonical''  arrow $a \rightarrow b$,
those canonical arrows  being the ones generated in an expected way 
from the canonical arrows from the monoidal structures and the
$\jj_a$'s.  Eventually Laplaza's paper also provides a combinatorial
description of free categorical groups.\\

Let us recall also the following known facts for any symmetric Picard 
categories $\A$ and $\B$.
For any symmetric monoidal functor 
$(F,F^2,F^0): \A \rightarrow \B$ the component $F^0$ is determined
by $F^2$. 
Actually a monoidal structure on a functor $F: \A \rightarrow \B$
is given by a natural $F^2_{a,b}: Fa \otimes Fb \rightarrow F(a \otimes b)$
in $\B$
satisfying the only axiom that 
{\small
$$
\xymatrix{
Fa \otimes (Fb \otimes Fc)
\ar[d]_{1 \otimes F^2_{b,c}}
\ar@{-}[r]^{\cong}
&
(Fa \otimes Fb) \cong Fc
\ar[d]^{F^2_{a,b} \otimes 1}
\\
Fa \otimes F(b \otimes c)
\ar[d]_{F^2_{a,b \otimes c}}
& 
F(a \otimes b) \otimes Fc
\ar[d]^{F^2_{a \otimes b,1}}
\\
F(a \otimes (b \otimes c)) 
\ar@{-}[r]_{F(\cong)}
&
F((a \otimes b) \otimes c)
}
$$
}
commutes for any objects $a$, $b$ and $c$ of $\A$.
Also any natural transformation 
$\sigma: F \rightarrow G: \A \rightarrow \B$ between monoidal
functors is monoidal if it satisfies the only axiom that 
the diagram in $\B$
{\small
$$\xymatrix{
Fa \otimes Fb 
\ar[r]^-{F^2_{a,b}} 
\ar[d]_{\sigma_a \otimes \sigma_{b}}
&
F(a \otimes b)
\ar[d]^{\sigma_{a \otimes b}}
\\
Ga \otimes Gb 
\ar[r]_-{G^2_{a,b}} 
&
G(a \otimes b)
}$$
}
commutes for any objects $a$,$b$ of $\A$.\\

Let us consider a symmetric Picard category $(\A,\jj)$.  
Since $\A$ is a groupoid, one
has a functor $\A^{op} \rightarrow \A$ which is the identity
on objects and sends arrows to their inverses.
The functor $\inv: \A \rightarrow \A$
is obtained by composing the previous functors and $(-)^{\eqi}$
above. Let us denote by $!$ the unique canonical arrow
between two objects of $\A$, when it exists.
\begin{lemma}\label{invmdal}
For any symmetric Picard category $(\A,\jj)$, the associated functor 
$\inv$ admits a symmetric monoidal structure
where $\inv^2$ has component 
$\inv^2_{a,b}: a^{\eqi} \otimes b^{\eqi} \rightarrow {(a \otimes b)}^{\eqi}$
in any $(a,b)$ the composite  
$\xymatrix{ a^{\eqi} \otimes b^{\eqi} 
\ar[r]^{\syc} 
&
b^{\eqi} \otimes a^{\eqi}
\ar[r]^{!}
&
{(a \otimes b)}^{\eqi}
}$    
(which is also 
$\xymatrix{ a^{\eqi} \otimes b^{\eqi} 
\ar[r]^{!} 
&
{(b \otimes a)}^{\eqi}
\ar[r]^{s^{\eqi}}
&
{(a \otimes b)}^{\eqi}
}$
according to Lemma \ref{syPic} in Appendix).
\end{lemma}
\pf 
See Appendix \ref{pfinvmdal}.
\epf

\begin{tag}\label{monatgp}
For any symmetric Picard category $(\A,\jj)$
one has a monoidal natural isomorphism
$\un \rightarrow \inv \Box id: \A \rightarrow \A$
which component in any object $a$ is
$\jj_a: \un \rightarrow a^{\eqi} \otimes a$.
\end{tag}
\pf See Appendix \ref{pfmonatgp}. \epf

Now given objects $\A$ and $\B$ in $\SMC$ 
with $(\B,\jj)$ symmetric Picard, a group structure is obtained 
on the hom $[\A,\B]$ in $\SMC$, which is a groupoid,
as follows.
The strict symmetric monoidal functor
$[\A,-]: [\B,\B] \rightarrow [[\A,\B],[\A,\B]]$
sends the monoidal transformation
$\jj: \un \rightarrow inv \Box id: \B \rightarrow \B$
of \ref{monatgp} to a monoidal transformation
$$\xymatrix{
\un 
\ar@{=}[r] 
&
[\A,\un]
\ar[r]^-{[\A,\jj]}
&
[\A, inv \Box id]
\ar@{=}[r]
&
[\A,inv] \Box [\A,id]
\ar@{=}[r]
&
[\A,\inv] \Box id
}: [\A,\B] \rightarrow [\A,\B].
$$
which we define as the $\jj$ on $[\A,\B]$.
This is to say that $F^{\eqi}$ for 
any symmetric monoidal $F$ is the composite 
$$\xymatrix{ 
\A \ar[r]^F 
&
\B \ar[r]^{\inv} 
&
\B }.$$
and the natural isomorphisms $\jj_F: \un \cong F^{\eqi} \Box F$
are pointwise ${({\jj}_F)}_a: \un \rightarrow (Fa)^{\eqi} \otimes Fa$.
Then for any monoidal $\sigma: F \rightarrow G: \A \rightarrow 
\B$ one has that
${\sigma}^{\eqi}: G^{\eqi} \rightarrow F^{\eqi}$ is pointwise
${(\sigma_a)}^{\eqi}: {(Ga)}^{\eqi} \rightarrow {(Fa)}^{\eqi}$.\\ 
We will always consider that this is the chosen
group structure on $[\A,\B]$ when considered
as an object of $\BS$. This structure is determined by that of $\B$.\\  

One has a notion of {\em strictness} for arrows in $\BS$,
which is {\em different} from the notion of strictness in $\SMC$.
Let us consider any symmetric Picard categories $\A$ and $\B$.
For a symmetric monoidal functor $F: \A \rightarrow \B$
one has the natural isomorphism 
\begin{tag}\label{isofuntag}
$$F(a^{\eqi}) \cong_a {(Fa)}^{\eqi}.$$
\end{tag}
defined precisely in Appendix-\ref{defisofuntag}.
A symmetric monoidal functor $F: \A \rightarrow \B$ 
is called a {\em strict} 
arrow in $\BS$ when it preserves strictly the monoidal 
structure and moreover preserves strictly the isomorphisms 
$\jj$ meaning that the natural isomorphism \ref{isofuntag}
is an identity or equivalently that for any object $a$ in $\A$,
$F$ sends $a^{\eqi}$ to ${(Fa)}^{\eqi}$ and 
$\jj_a : \un \rightarrow a^{\eqi} \otimes a$
to $\jj_{Fa}: \un \rightarrow {(Fa)}^{\eqi} \otimes Fa$.
We write $\StrBS$ for the sub-2-category of $\BS$
with same objects, strict arrows and 2-cells inherited
from $\BS$.\\

The results from \cite{Sch08} regarding $\SMC$
transpose rather straightforwardly to $\BS$ as follows.\\

Given an arbitrary symmetric Picard category $\C$,
it happens that the isomorphism $\Dual_{\A,\B,\C}:
[\A,[\B,\C]] \rightarrow [\B,[\A,\C]]$ defined
in chapter 6 is also a strict arrow in $\BS$. 
For any arrow $F: \A \rightarrow [\B,C]$
in $\BS$, its image $F^*$ by $\Dual$, called it ``dual'', 
is strict in $\BS$ if and only for any object $b$ of $\B$, 
the arrow $F^*(b): \A \rightarrow \C$ is strict in $\BS$.
For any objects $\A$,$\B$ and $\C$ in $\BS$,
the arrow ${[\A,-]}_{\B,\C}:[\B,\C] \rightarrow [[\A,\B],[\A,\C]]$
defined in chapter 8 is strict in $\BS$.\\

The hom 2-functor of $\SMC$ defined in chapter 9 induces 
by restriction a hom 2-functor $\BS^{op} \times \BS \rightarrow \BS$
for $\BS$. The statements regarding the 2-naturality of $\Dual$
(chapter 10) and the evaluation functors (chapter 11)
still hold when replacing formally $\SMC$ by $\BS$. In particular 
the evaluation functors are strict arrows in $\BS$.\\ 

Similarly to the case of the 2-category $\SMC$, one has a tensor product
in $\BS$. For any symmetric Picard categories $\A$ and $\B$,
their tensor $\A \otimes \B$ satisfies the universal 
property of the existence of a 2-natural isomorphism
\begin{tag}\label{punitens}  
$$\BS(\A,[\B,\C]) \cong_{\C} \StrBS(\A \otimes \B, \C)$$
\end{tag}
between 2-functors $\StrBS \rightarrow \Cat$ in the argument $\C$.
Note that the 2-naturality in question involves only {\em strict} 
morphisms in $\BS$.
We briefly sketch a description of the above tensor $\A \tens \B$ 
by generator and relations. It is similar 
to that given with more details in \cite{Sch08} for the tensor 
in $\SMC$.\\

We consider a graph $\Gra$ with vertices
the terms of the free $\{\un, {(-)}^{\eqi}, \otimes \}$-algebra over 
the set $\Obj(\A) \times \Obj(\B)$, i.e. they are words 
of the formal language containing all pairs $(a,b)$ -- 
which we write $a \otimes b$ -- for objects $a$ of $\A$
and $b$ of $\B$, the one-symbol word $\un$, the words
$X^{\eqi}$ for any vertex $X$, and 
$X \otimes Y$ for any vertices $X$ and $Y$. 
The set of edges of $\Gra$ consists of:\\
- The ``canonical'' edges for the symmetric monoidal 
structure which are the 
$\ac_{X,Y,Z}: X \otimes (Y \otimes Z) \rightarrow (X \otimes Y) \otimes Z$,      $\rc_X: X \otimes \un \rightarrow X$,
$\lc_X: \un \otimes X \rightarrow X$,
$\syc_{X,Y}: X \otimes Y \rightarrow Y \otimes X$
for all vertices $X$,$Y$,$Z$;\\
- Edges $\jj_X : \un \rightarrow {X}^{\eqi} \otimes X$, one for
each vertex $X$;\\  
- Edges
$\gamma_{a,a',b}: (a \otimes b) \otimes (a' \otimes b) 
\rightarrow  (a \otimes a') \otimes b$
and
$\delta_{a,b,b'}: (a \otimes b) \otimes (a \otimes b')
\rightarrow 
a \otimes (b \otimes b')$
indexed by objects $a$,$a'$ of $\A$ and $b$,$b'$ of $\B$;\\
- Edges $a \otimes f: a \otimes b \rightarrow a \otimes b'$
indexed by 
objects $a$ of $\A$ and arrows $f: b \rightarrow b'$
of $\B$;\\ 
- Edges $f \otimes b: a \otimes b \rightarrow a' \otimes b$
indexed by objects $b$ of $\B$ and arrows $f: a \rightarrow a'$
of $\A$;\\
- Edges $X \otimes p: X \otimes Y \rightarrow X \otimes Z$ 
and $p \otimes X: Y \otimes X \rightarrow Z \otimes X$
for any vertex $X$ and any edge $p: Y \rightarrow Z$;\\  
with the convention that edges above with different names 
are different.\\

Let us consider $\FGd(\Gra)$ the free groupoid on $\Gra$, i.e its arrows
are mere concatenations of edges of $\Gra$ and their formal inverses.
For any vertex $X$, one has two graph endomorphisms of $\Gra$,
namely $X \otimes -$ and $- \otimes X$ sending respectively 
an arbitrary edge $f:Y \rightarrow Z$ to 
$X \otimes Y \rightarrow X \otimes Z$, resp.
$Y \otimes X \rightarrow Z \otimes X$. These two
extend uniquely to endofunctors of $\FGd(\Gra)$ and 
we extend the notation $X \otimes f$ and $f \otimes X$ to denote
the images of arrows of $\FGd(\Gra)$ by these functors.\\

The tensor $\A \otimes \B$ is the quotient of $\FGd(\Gra)$ 
by the congruence generated by the following relations $\sim$
on its arrows from \ref{cong6} to \ref{cong12} below.\\

For all edges $\xymatrix{X \ar[r]^t & Y}$ and
$\xymatrix{Z \ar[r]^s & W}$ of $\Gra$,
\begin{tag} \label{cong6}
$$\xymatrix{
&
X \otimes W
\ar[rd]^{t \otimes W}
&
\\
X \otimes Z 
\ar[ru]^{X \otimes s}
\ar[rd]_{t \otimes Z}
&
\sim
&
Y \otimes W. 
\\
&
Y \otimes Z
\ar[ru]_{Y \otimes s}
&
}$$
\end{tag}

\begin{tag}
\noindent Relations giving the coherence conditions
for $\ac$, $\rc$, $\lc$ and $\syc$ in $\A \otimes \B$. 
\end{tag}
These are the following.\\ 
- For any vertices $X$, $Y$, $Z$ and $T$,
$$\xymatrix{ X \otimes (Y \otimes (Z \otimes T)) 
\ar@{}[rrd]|{\sim}
\ar[r]^{\ac} 
\ar[d]_{1 \otimes \ac} & 
(X \otimes Y) \otimes (Z \otimes T) 
\ar[r]^{\ac}  & 
 ( ( X \otimes Y) \otimes Z) \otimes T \\
X \otimes ( ( Y \otimes Z) \otimes T) 
\ar[rr]_{\ac} & &
( X \otimes ( Y \otimes Z )) \otimes T. 
\ar[u]_{\ac \otimes 1} }$$
- For any vertices $X$ and $Y$,
$$\xymatrix{ X \otimes ( \un \otimes Y) \ar[rr]^{\ac} 
\ar[rd]_{1 \otimes \lc} &  \ar@{}[d]|{\sim} &  
(X \otimes \un) \otimes Y \ar[ld]^{\rc \otimes 1} \\
& X \otimes Y. &  }$$
- For any vertex $X$,
$$\xymatrix{ X \otimes \un \ar[rr]^{\syc} \ar[rd]_{\rc} 
& \ar@{}[d]|{\sim} & \un \otimes X 
\ar[ld]^{\lc} \\ 
& X. & }$$ 
- For any vertices $X$, $Y$ and $Z$,
$$\xymatrix{
X \otimes ( Y \otimes Z) \ar[r]^{\ac} 
\ar[d]_{1 \otimes s} & 
( X \otimes Y ) \otimes Z \ar[r]^{\syc} \ar@{}[d]|{\sim} &
Z \otimes ( X \otimes Y) \ar[d]^{\ac}\\
X \otimes ( Z \otimes Y) \ar[r]_{\ac} & 
(X \otimes Z) \otimes Y & (Z \otimes X) \otimes Y. 
  \ar[l]^{\syc \otimes 1} 
}$$

\begin{tag}\label{cong9}
Relations for the naturalities of $\ac$, $\rc$, $\lc$,
and $\syc$ in $\A \otimes \B$.\end{tag}
For instance, one has 
for any edge $f:X \rightarrow X'$ of $\Gra$,
and any vertices $Y$ and $Z$, 
$$\xymatrix@C=3pc{ 
X \otimes ( Y \otimes Z) 
\ar[r]^{\ac_{X,Y,Z}}
\ar@{}[rd]^{\sim}
\ar[d]_{f \otimes 1} &
(X \otimes Y) \otimes Z \ar[d]^{(f \otimes 1) \otimes 1}\\
X' \otimes ( Y \otimes Z) \ar[r]_{\ac_{X',Y,Z}}  &
(X' \otimes Y) \otimes Z.}
$$
We will not write here the other relations.
There are two more for the naturalities of $\ac_{X,Y,Z}$ in $Y$ and $Z$, 
one for that of $\lc_X$ in $X$,
one for that of $\rc_X$ in $X$
and two for those of $\syc_{X,Y}$
in $X$ and $Y$.\\ 

For any object $a$ in $\A$ and any arrows
$\xymatrix{b \ar[r]^{f}  & b' \ar[r]^{g} & b''}$
in $\B$, 
\begin{tag}\label{cong131}
$$\xymatrix{
a \otimes b 
\ar[rr]^{a \otimes (g \circ f)}
\ar[rd]_{a \otimes f}
&
\ar@{}[d]|{\sim}
&
a \otimes b''\\
&
a \otimes b'
\ar[ru]_{a \otimes g}
&
}$$
\end{tag}

For any object $b$ in $\B$ and any arrows 
$\xymatrix{a \ar[r]^{f}  & a' \ar[r]^{g} & a''}$
in $\A$,
\begin{tag}\label{cong133}
$$\xymatrix{
a \otimes b 
\ar[rr]^{(g \circ f) \otimes b}
\ar[rd]_{f \otimes b}
&
\ar@{}[d]|{\sim}
&
a'' \otimes b\\
&
a' \otimes b
\ar[ru]_{g \otimes b}
&
}$$
\end{tag}

For any objects $a$ in $\A$ and $b$ in $\B$,
\begin{tag}\label{cong132}
$a \otimes {id}_b \sim {id}_{a \otimes b}$
\end{tag}
and
\begin{tag}\label{cong134}
$ {id}_a \otimes b \sim  {id}_{a \otimes b}$. 
\end{tag} 
where ${id}_b$, ${id}_a$ and ${id}_{a \otimes b}$ above
are the identities respectively at $b$ in $\B$, at $a$ in $\A$ 
and at $a \otimes b$ in $\FGd(\Gra)$.\\

For any arrows $f: a \rightarrow a'$ in $\A$ and
$g: b \rightarrow b'$ in $\B$, 
\begin{tag}\label{cong135}
$$\xymatrix{a \otimes b \ar[d]_{a \otimes g} 
\ar[r]^{f \otimes b}
\ar@{}[rd]|{\sim} & 
a' \otimes b \ar[d]^{a' \otimes g}\\
a \otimes b' \ar[r]_{f \otimes b'}& a' \otimes  b'.}$$
\end{tag}

\begin{tag}\label{cong10}
Relations for the ``naturalities'' 
of $\gamma_{a,a',b}$ in $a$, $a'$ and $b$
and $\delta_{a, b, b'}$ in $a$, $b$ and $b'$.
\end{tag}
For instance by the relations for the ``naturality'' of $\gamma_{a,a',b}$
in $b$ it is meant that for any objects $a,a'$ in $\A$ and any 
arrow $g: b \rightarrow b'$ in $\B$, 
$$\xymatrix{ (a \otimes b) \otimes (a' \otimes b) 
\ar[d]_{ (1 \otimes g) \otimes (1 \otimes g) }
\ar[r]^-{\gamma_{a,a',b}} 
\ar@{}[rd]|{\sim}
& (a \otimes a') \otimes b 
\ar[d]^{1 \otimes g}\\
(a \otimes b')  \otimes (a' \otimes b') 
\ar[r]_-{\gamma_{a, a', b'}} 
& 
(a \otimes a') \otimes b'.
}$$  
We will not write explicitly now the five 
other relations.\\
 
For any objects $a$ in $\A$ and $b$, $b'$, $b''$ in $\B$,
\begin{tag}\label{cong113}
$$\xymatrix{ 
(a \otimes b) \otimes ( (a \otimes b') \otimes (a \otimes b'') )
\ar[r]^{\ac} \ar[d]_{ 1 \otimes \delta_{a,b',b''} }
\ar@{}[rdd]|{\sim}
& 
( (a \otimes b) \otimes (a \otimes b') ) \otimes (a \otimes b'')
\ar[d]^{ \delta_{a,b,b'} \otimes 1 }
\\
(a \otimes b) \otimes (a \otimes (b' \otimes b''))
\ar[d]_{ \delta_{a,b, b' \otimes b''} } & 
(a \otimes (b \otimes b')) \otimes (a \otimes b'')
\ar[d]^{ \delta_{a, b \otimes b' ,  b''} }\\
a \otimes (b \otimes (b' \otimes b'')) 
\ar[r]_{ 1 \otimes \ac_{ b, b',  b'' } } &
a \otimes ((b \otimes b') \otimes b''). 
}$$
\end{tag}

For any objects $a$ in $\A$ and $b$, $b'$ in $\B$,
\begin{tag}\label{cong115}
$$\xymatrix{
(a \otimes b) \otimes (a \otimes b') \ar[r]^-{ \delta_{a,b,b'} }
\ar[d]_{ \syc_{a \otimes b,a \otimes b'} } 
\ar@{}[rd]|{\sim} & 
a \otimes (b \otimes b') 
\ar[d]^{1 \otimes \syc_{b,b'}}\\
(a \otimes b') \otimes (a \otimes b) \ar[r]_-{\delta_{a,b',b}} & 
a \otimes (b' \otimes b).}$$
\end{tag}

For any objects $a$, $a'$, $a''$ in $\A$ and $b$ in $\B$,
\begin{tag}\label{cong143}
$$\xymatrix{ 
(a \otimes b) \otimes ( (a' \otimes b ) \otimes (a'' \otimes b ))
\ar[r]^{\ac} \ar[d]_{1 \otimes \gamma_{a',a'',b}}
\ar@{}[rdd]|{\sim}
& 
( (a \otimes b) \otimes (a' \otimes b) ) \otimes (a'' \otimes b)
\ar[d]^{ \gamma_{a,a',b} \otimes 1 }\\
(a \otimes b) \otimes ( (a' \otimes a'') \otimes b ))
\ar[d]_{ \gamma_{a,a' \otimes a'', b} }& 
((a \otimes a') \otimes b) \otimes (a'' \otimes b)
\ar[d]^{\gamma_{a \otimes a', a'' , b}}\\
(a \otimes (a' \otimes a'')) \otimes b 
\ar[r]_{ \ac \otimes 1 } &
( (a \otimes a') \otimes a'') \otimes b.
}$$
\end{tag}

For any objects $a$, $a'$ in $\A$ and $b$ in $\B$,
\begin{tag}\label{cong145}
$$\xymatrix{
(a \otimes b) \otimes (a' \otimes b) 
\ar[r]^-{\gamma_{a,a',b}} 
\ar[d]_{s_{a \otimes b, a' \otimes b}}
\ar@{}[rd]|{\sim} & 
(a \otimes a') \otimes b 
\ar[d]^{s_{a, a'} \otimes 1}\\
(a' \otimes b) \otimes (a \otimes b) 
\ar[r]_-{\gamma_{a',a,b}} & (a' \otimes a) \otimes b
.}$$  
\end{tag}

For any objects $a$,$a'$ in $\A$ and $b$,$b'$ in $\B$,
\begin{tag}\label{cong19}
$$\xymatrix{
((a \otimes b) \otimes (a \otimes b'))
\otimes 
((a' \otimes b) \otimes (a' \otimes b'))
\ar@{-}[r] 
\ar[d]_{\delta_{a,b,b'} \otimes \delta_{a',b,b'} }
\ar@{}[rd]|{\sim} 
&
((a \otimes b) \otimes (a' \otimes b))
\otimes 
((a \otimes b') \otimes (a' \otimes b'))
\ar[d]^{\gamma_{a,a',b} \otimes \gamma_{a,a',b'}}\\
(a \otimes (b \otimes b')) \otimes (a' \otimes (b \otimes b')) 
\ar[d]_{ \gamma_{a,a', b \otimes b'} } &
((a \otimes a') \otimes b) 
\otimes 
((a \otimes a') \otimes b') \ar[ld]^{\delta_{a \otimes a', b, b'}}\\
(a \otimes a') \otimes (b \otimes b'). &
 }$$
where the top arrow is the concatenation
$$\xymatrix@C=5pc{
\ar[r]^{ \ac_{X \otimes Y,Z,T} }
&
\ar[r]^{ \ac_{X,Y,Z}^{-1} \otimes 1 }
&
\ar[r]^{ (1 \otimes \syc_{Y,Z}) \otimes T }
&
\ar[r]^{ {\ac_{X,Z,Y}} \otimes T }
&
\ar[r]^{ {\ac_{X \otimes Z,Y,T}}^{-1} }  
&
}$$
with the $\ac^{-1}$ being the formal inverses 
of edges $\ac$ and $X$, $Y$, $Z$ and $T$
standing respectively for $a \otimes b$,
$a \otimes b'$, $a' \otimes b$ and $a' \otimes b'$.
\end{tag}

\begin{tag}\label{cong12}
Expansions of all relations above by iterations of
$X \otimes -$ and $- \otimes X$ for all vertices $X$.
\end{tag}
Which means precisely that the set of relations $\sim$ is the 
smallest set of relations on arrows of $\FGd(\Gra)$ 
containing the previous relations
(\ref{cong6} to \ref{cong19})
and satisfying the closure properties that
for any relation $f \sim g : Y \rightarrow Z$ that it contains
and any vertex $X$, 
it contains also the relations
\begin{center}
$X \otimes f \sim  X \otimes g: X \otimes Y \rightarrow X \otimes Z$
\end{center}
and 
\begin{center}
$f \otimes X \sim  g \otimes X: Y \otimes X \rightarrow Z \otimes X$.
\end{center} 

The proofs that the above category $\A \otimes \B$ is
a well defined symmetric Picard category and that it satisfies 
the universal property \ref{punitens}, are similar to 
those in \cite{Sch08} for the well definition
and universal property of the tensor product in $\SMC$.  
We will therefore not replicate them.
For any objects $\A$, $\B$ and $\C$ in $\BS$,
one obtains an adjunction 
$$\Ext \dashv \Res: [\A \otimes \B,\C] \rightarrow [\A,[\B, \C]]$$ 
in the 2-category $\BS$ (in this case it is an equivalence)
with $\Res \circ \Ext = 1$ and
where the arrows $\Res_{\A,\B,\C}$ are strict.
The isomorphism \ref{punitens} becomes 2-natural 
in $\A$ and $\B$ for a unique
tensor 2-functor $\BS \times \BS \rightarrow \BS$.\\

There exists a free symmetric Picard category on one 
generator, which we shall write $\un$, that differs obviously
from the ``unit'' for $\SMC$ written $\un$ and defined 
in \cite{Sch08}-chapter 18,
but that has a very similar presentation by generators and relation. 
The only differences are the following.
Its set objects is now the free $\{\un, \otimes, {(-)}^{\eqi} \}$-algebra 
over one generator $\gen$.
It is a quotient of free {\em groupoid}  generated by the graph
containing the usual canonical edges for the symmetric monoidal
structures (the $\ac_{X,Y,Z}$, $\rc_{X}$, $\lc_X$, $\syc_{X,Y}$)
and moreover containing one edge $\jj_X: \un \rightarrow X^{\eqi} \otimes X$
for each object $X$. The relations on this free groupoid
defining $\un$ are just those for the naturalities of the 
collections $\ac$, $\rc$, $\lc$
and $\syc$, those expressing the coherence axioms for the 
symmetric monoidal structure, and relations expressing the bifunctoriality
of $-\otimes-: \un \times \un \rightarrow \un$.
The universal property defining $\un$ is that for any symmetric
Picard category $\A$, there exits a unique strict arrow 
$\vv: \un \rightarrow [\A,\A]$ such that $\vv(\gen)$ is 
the identity arrow at $\A \rightarrow \A$ with its strict structure
in $\BS$.
One obtains with similar proofs, similar results. Namely:
\begin{tag}\label{unimonadj}
For any symmetric Picard category $\A$, the dual 
$\vs: \A \rightarrow [\unc,\A]$ of $\vv: \unc \rightarrow [\A,\A]$
has right adjoint in $\BS$ the evaluation at $\gen$ functor
$\ev_{\gen}: [\unc,\A] \rightarrow \A$.
\end{tag}

From this one can exhibit 
a kind of ``symmetric monoidal closed 2-structure''
on $\BS$ 
in the same way as done in \cite{Sch08} chapters
19, 20 and 21 for $\SMC$.
Namely one can define the canonical arrows 
$A'_{\A,\B,\C}: (\A \otimes \B) \otimes \C 
\rightarrow \A \otimes (\B \otimes \C)$,
$R'_{\A}: \A \rightarrow \A \otimes \un$,
$L'_{\A}: \A \rightarrow \un \otimes \A$
and $S_{\A,\B}: \A \otimes \B \rightarrow \B \otimes \A$
in this case with respective inverse equivalences
$A_{\A,\B,\C}$, $R_{\A}$, $L_{\A}$ and $S_{\B,\A}$
and satisfying the (strict) coherence axioms given
in chapter 20. Eventually this 2-categorical structure
induces a symmetric monoidal closed structure on $\BS/\sim$ where
$\sim$ denotes the congruence generated by the 2-cells
of $\BS$. 
\end{section}

\begin{section}{$\BS$-categories, $\BS$-functors and
$\BS$-natural transformations}\label{BSCat}
$\BS$-categories, $\BS$-functors and $\BS$-natural 
transformations have been considered by M.Dupont
in his thesis \cite{Dup08}. They are respectively
bicategories with homs in $\BS$, and
pseudo-functors and pseudo-natural transformations
with linear components, i.e with arrows and 2-cells
in the 2-category $\BS$ rather than in $\Cat$.
As such they obviously form a 2-category
with a forgetful 2-functor to the 2-category 
of bicategories, pseudo-functors and pseudo natural 
transformations.\\

We start by recalling these notions which are actually 
slight (enriched) variations of the usual notions of 
bicategories, pseudo-functors and pseudo-natural 
transformations.\\

By a n-linear natural transformation $\sigma$ between n-linear maps
$F$ and $G$, written
$\sigma: F \rightarrow G: \A_{1} \times ... \times \A_n \rightarrow \B$
we will mean a 2-cell of $[\A_1,[\A_2,...,[\A_n,\B]...]$.
Multi-linear maps and multilinear natural transformation
compose in an evident way, which can be justified by the 2-natural 
isomorphism $\Dual_{\A,\B,\C}: [\A,[\B,\C]] \cong [\B,[\A,\C]]$ 
and the existence of a forgetful 2-functor 
$\BS \rightarrow \Cat$ (see \cite{Sch08} chapter 6).\\

An {\em enrichment} $(\A, \compb, \unitb, \assb, \rb, \lb)$, 
over $\BS$ also named an {\em $\BS$-category} 
and which we might sometimes denote simply by $\A$,
consists of the following data:\\
- A small set with elements $x$,$y$,$z$... called the {\em objects}
of $\A$.\\
- A map $\A$ sending any pair $x$,$y$ of objects to an object 
$\A(x,y)$ of $\BS$ sometimes also written $\A_{x,y}$ for convenience
and called the {\em hom} of $x$ and $y$.\\
- A collection of bilinear maps 
$\compb_{x,y,z}: \A_{y,z} \times \A_{x,y} \rightarrow \A_{x,z}$, 
the {\em composition} maps
indexed by objects $x$,$y$,$z$ of $\A$ and we write 
$g \circ f$ for $\compb_{x,y,z}(g,f)$ for any objects 
$g$ of $\A_{y,z}$ 
and $f$ of $\A_{x,y}$ and $\tau * \sigma$ for 
$\compb_{x,y,z}(\tau,\sigma)$ for any arrows
$\tau$ of  $\A_{y,z}$ and $\sigma$ of $\A_{x,y}$.\\ 
- A collection of objects $\one_x$ of $\A_{x,x}$ indexed by objects
$x$ of $\A$;\\ 
- Collections of natural transformations $\assb_{x,y,z,t}$, which 
are {\em trilinear}, indexed by objects $x$,$y$,$z$,$t$,
and $\rb_{x,y}$ and
$\lb_{x,y}$ both {\em linear}, as follows
\begin{tag}\label{assb}
$${(\assb_{x,y,z,t})}_{h,g,f}: 
h \circ (g \circ f) \rightarrow h \circ (g \circ f)$$
\end{tag}
which lies in $\A_{x,t}$
for objects $h$ in $\A_{z,t}$, $g$ in $\A_{y,z}$ and $f$ in $\A_{x,y}$
\begin{tag}\label{rb}
$${(\rb_{x,y})}_f: f \rightarrow f \circ \one_x$$
\end{tag}
which lies in $\A_{x,y}$ 
for objects $f$ in $\A_{x,y}$;
\begin{tag}\label{lb}
$${(\lb_{x,y})}_f:  f \rightarrow \one_y \circ f$$
\end{tag} 
which lies in $\A_{x,y}$ for objects $f$ in $\A_{x,y}$;
and those are subjects to the coherence axioms
\ref{cohassb} and \ref{cohunassb} below.
\begin{tag}\label{cohassb}
For any objects $x$, $y$, $t$ and $u$ of $\A$, 
any objects $f$ of $\A_{x,y}$,
$g$ of $\A_{y,z}$,
$h$ of $\A_{z,t}$
and $k$ of $\A_{t,u}$,
the diagram in $\A_{x,u}$
$$\xymatrix{
k \circ (h \circ (g \circ f))) 
\ar[rr]^-{{(\assb_{x,z,t,u})}_{k,h,g \circ f}}
\ar[d]_{k * {(\assb_{x,y,z,t})}_{h,g,f}  }
&
&
(k \circ h) \circ (g \circ f)
\ar[d]^{{(\assb_{x,y,z,u})}_{k \circ h, g,f}}
\\
k \circ ( (h \circ g) \circ f)
\ar[r]_-{{(\assb_{x,y,t,u})}_{k, h \circ g, f}} 
& 
(k \circ  (h \circ g)) \circ f
\ar[r]_-{{(\assb_{y,z,t,u})}_{k, h, g} * f}
& 
((k \circ h) \circ g) \circ f
}$$
commutes.
\end{tag}
\begin{tag}\label{cohunassb}
For any objects $x$, $y$, $z$ in $\A$ and any 
objects $f$ of $\A(x,y)$ and $g$ of $\A(y,z)$
the diagram in $\A(x,y)$
$$
\xymatrix{
& g \circ f 
\ar[ld]_{g * {(\lb_{x,y})}_f }
\ar[rd]^{{(\rb_{y,z})}_g * f}
& \\
g \circ (\one_y \circ f)
\ar[rr]_{{(\assb_{x,y,y,z})}_{g,\one_y,f}}
& 
& 
(g \circ \one_y) \circ f
}
$$
commutes.
\end{tag}

For any $\BS$-category $\A$, its underlying 
bicategory is denoted $\A^0$.\\

Given two arbitrary $\BS$-categories $\A$ and $\B$, 
a {\em $\BS$-functor} $F: \A \rightarrow \B$ 
consists of the following data.\\
- A map $F$ sending objects of $\A$ to objects
of $\B$;\\
- Arrows $F_{x,y}: \A_{x,y} \rightarrow \B_{Fx,Fy}$ in $\BS$
for each pair of objects $x$,$y$ in $\A$.\\
- A collection of {\em bilinear} natural transformations
$F^2_{x,y,z}$ indexed by objects $x$,$y$,$z$ of $\A$ 
with components
$${(F^2_{x,y,z})}_{g,f}: F_{y,z}(g) \circ F_{x,y}(f) 
\rightarrow F_{x,z}(g \circ f)$$
in $\B_{Fx,Fz}$ for objects 
$g$ in $\A_{y,z}$ and $f$ in $\A_{x,y}$.\\
- A collection of arrows 
$F^0_x: 1_{Fx} \rightarrow F_{x,x}(1_x)$ in $\B_{Fx,Fx}$ 
indexed by objects $x$ of $\A$.\\
Those are subjects to the coherence axioms  
\ref{cohfun1b}, \ref{cohfun2b} 
and \ref{cohfun3b} below.
\begin{tag}\label{cohfun1b}
For any objects $x$,$y$,$z$,$t$ of $\A$,
and any objects $f$ of $\A_{x,y}$,
$g$ of $\A_{y,z}$ and $h$ of $\A_{z,t}$ the diagram 
in $\B(Fx,Ft)$
$$
\xymatrix@C=4pc{
F_{z,t}h \circ (F_{y,z}g \circ F_{x,y}f) 
\ar[r]^{  {(\assb_{Fx,Fy,Fz,Ft})}_{Fh,Fg,Ff} }
\ar[d]_{1 * {(F^2_{x,y,z})_{g,f}} }
& 
(F_{z,t}h \circ F_{y,z}g) \circ F_{x,y}f
\ar[d]^{ {(F^2_{y,z,t})}_{h,g} * 1 }
\\
F_{z,t}h \circ F_{x,z}(g \circ f)
\ar[d]_{{(F^2_{x,z,t})}_{h,g \circ f} }
& 
F_{y,t}(h \circ g) \circ F_{x,y}f
\ar[d]^{{ (F^2_{x,y,t}) }_{h \circ g,f} }
\\
F_{x,t}(h \circ (g \circ f))
\ar[r]_{  F_{x,t}(  {(\assb_{x,y,z,t})}_{h,f,g} )  }
& 
F_{x,t}((h \circ g) \circ f) 
}
$$
commutes.
\end{tag}
\begin{tag}\label{cohfun2b}
For any objects $x$,$y$ of $\A$ and
any object $f$ of $\A_{x,y}$, 
the diagram in $\B_{Fx,Fy}$
$$\xymatrix@C=3pc{
Ff 
\ar[d]_{ F( {\rb_{x,y}}_f ) }
\ar[r]^{{ (\rb_{Fx,Fy}) }_{Ff}}  
&
Ff \circ \one_{Fx}
\ar[d]^{1 * F^0_x}
\\
F(f \circ \one_x)
&
Ff \circ F(\one_x)
\ar[l]^{ { (F^2_{x,x,y}) }_{f,\one_x} }
}
$$
commutes.
\end{tag}
\begin{tag}\label{cohfun3b}
For any objects $x$,$y$ of $\A$
and any object $f$ of $\A_{x,y}$,
the diagram in $\B_{Fx,Fy}$ 
$$\xymatrix@C=3pc{
Ff 
\ar[d]_{ F( {\lb_{x,y}}_f ) }
\ar[r]^{{(\lb_{Fx,Fy})}_{Ff}}  
&
\one_{Fy} \circ Ff
\ar[d]^{F^0_y * 1}
\\
F(\one_y \circ f)
&
F(\one_y) \circ Ff
\ar[l]^{ {(F^2_{x,y,y})}_{\one_y,f} }
}
$$
commutes.
\end{tag}

Given two $\BS$-functors $F,G: \A \rightarrow \B$
a {\em $\BS$-natural transformation $(\sigma,\kappa)$} 
consists in a family of arrows $\sigma_x$ of $\B_{Fx,Gx}$
indexed by objects $x$ of $\A$ together with 
a collection of linear natural transformations 
$\kappa_{x,y}$ indexed by objects $x$ and $y$ as follows
$${(\kappa_{x,y})}_f: 
Gf \circ \sigma_x \rightarrow \sigma_y \circ Ff$$
lies in $\B_{Fx,Gy}$ for objects $f$ of $\A_{x,y}$,
and these satisfy the coherence axioms \ref{cohnat1b}
and \ref{cohnat2b} below.
\begin{tag}\label{cohnat1b}
For any objects $f$ in $\A_{x,y}$ and $g$ in $\A_{y,z}$, 
the diagram in $\B_{Fx,Gz}$  
$$
\xymatrix@C=3pc{
(Gg \circ Gf) \circ \sigma_x
\ar[r]^{ {(G^2_{x,y,z})}_{g,f} * \sigma_x  }
&
G(g \circ f) \circ \sigma_x
\ar[r]^{ {(\kappa_{x,z})}_{ g \circ f } }
&
\sigma_z \circ F(g \circ f)
\\
Gg \circ (Gf \circ \sigma_x)
\ar[u]^{ { ( \assb_{Fx,Gx,Gy,Gz} ) }_{Gg,Gf,\sigma_x} }
\ar[d]_{1 * {(\kappa_{x,y})}_f }
& 
&
\sigma_z \circ (Fg \circ Ff)
\ar[d]^{ {(\assb_{Fx,Fx,Fz,Gz})}_{\sigma_z,Fg,Ff} }
\ar[u]_{ 1*{(F^2_{x,y,z})}_{g,f} }
\\
Gg \circ (\sigma_y \circ Ff)
\ar[r]_{ {(\assb_{Fx,Fy,Gy,Gz})}_{Gg,\sigma_y,Ff} }
&
(Gg \circ \sigma_y) \circ Ff
\ar[r]_{ {(\kappa_{y,z})}_g * 1}
&
(\sigma_z \circ Fg) \circ Ff
}
$$
commutes. 
\end{tag}
\begin{tag}\label{cohnat2b}
For any object $x$ of $\A$,
the diagram in $\B_{Fx,Gx}$
$$
\xymatrix{
& \sigma_x 
\ar[ld]_{{(\rb_{Fx,Gx})}_{\sigma_x}}
\ar[rd]^{{(\lb_{Fx,Gx})}_{\sigma_x}}
& 
\\
\sigma_x \circ \one_{Fx} 
\ar[d]_{\sigma_x * F^0_x}
& 
& 
\one_{Gx} \circ \sigma_x
\ar[d]^{G^0_x * \sigma_x} 
\\
\sigma_x \circ F(\one_x) 
& 
& 
G(\one_x) \circ \sigma_x
\ar[ll]^{ { ( \kappa_{x,x} ) }_{\one_x} }
}
$$
commutes.
\end{tag}

We want to give alternative definitions of $\BS$-categories,
$\BS$-functors and $\BS$-natural transformations
by means of commuting diagrams in $\BS$ in a first 
instance without using the tensor.
They are obtained by replacing multilinear maps and multilinear
natural transformations from the previous definition
by corresponding arrows and 2-cells in $\BS$.
This yields the following.\\

One can define a $\BS$-category 
as a collection of objects $\A$, 
with homs $\A_{x,y}$ in $\BS$ as before, with 
collections of arrows\\
- $\A(x,-)_{y,z}: \A_{y,z} \rightarrow [\A_{x,y},\A_{x,z}]$
in $\BS$, indexed by objects $x$, $y$, $z$ of $\A$
with dual $\A_{-,y}: \A_{x,y} \rightarrow [\A_{y,z},\A_{x,z}]$
written $\A(-,y)$;\\
- $\unit_x : \unc \rightarrow \A_{x,x}$ indexed 
by $x$, which are strict;\\ 
and collections of 2-cells:\\
- $\assh_{x,y,z,t}$ in $\BS$,
indexed by objects $x$, $y$, $z$ and $t$ of $\BS$ as follows
\begin{tag}\label{assh}
{\small
$$
\xymatrix{
&
\A_{z,t}
\ar[rd]^{\A(y,-)}
\ar[ld]_{\A(x,-)}
&
\\
[\A_{x,z},\A_{x,t}]
\ar[d]_{[\A_{x,y},-]}
&
&
[\A_{y,z},\A_{y,t}]
\ar[d]^{[1,\A(x,-)]}
\\
[[\A_{x,y},\A_{x,z}],[\A_{x,y},\A_{x,t}]]
\ar[rr]_{[\A(x,-),1]}
\ar@{=>}[rru]^{\assh}
&
&
[\A_{y,z},[\A_{x,y},\A_{x,t}]]
}
$$
}
\end{tag}
- $\rh_{x,y}$ and $\lh_{x,y}$ indexed by objects $x$ and $y$ 
respectively as follows
\begin{tag}\label{rh}
$\rh_{x,y}:$ 
{\small
$$
\xymatrix{
\A_{x,y} 
\ar[r]^-{\A(x,-)}
\ar@{-}[dd]_{id}
&
[\A_{x,x},\A_{x,y}]
\ar[dd]^{[\unit_x,1]}
\\
\ar@{=>}[r]
&
\\
\A_{x,y}
&
[\unc,\A_{x,y}]
\ar[l]^{\ev_{\gen}}
}
$$
}
\end{tag}
and 
\begin{tag}\label{lh}
$\lh_{x,y}:$ 
{\small
$$
\xymatrix{
\A_{x,y} 
\ar[r]^-{\A(-,y)}
\ar@{-}[dd]_{id}
&
[\A_{y,y},\A_{x,y}]
\ar[dd]^{[\unit_y,1]}
\\
\ar@{=>}[r]
&
\\
\A_{x,y}
&
[\unc,\A_{x,y}]
\ar[l]^{\ev_{\gen}}
}
$$
}
\end{tag}
those satisfying coherence axioms \ref{cohassh} and \ref{cohunassh}
below
\begin{tag}\label{cohassh}
The 2-cells in $\BS$
{\tiny
$$\xymatrix@=1pc@C=3pc{
\A_{t,u} 
\ar[r]^-{\A(x,-)} 
\ar@{-}[d]_{id}
\ar@<50ex>@{=>}[d]|{\assh_{x,z,t,u}}
&
\txt{
$[\A_{x,t},$\\ 
$\A_{x,u}]$
} 
\ar[r]^-{[\A_{x,z},-]}
&
\txt{
$[[\A_{x,z},\A_{x,t}],$\\
$[\A_{x,z},\A_{x,u}]]$
}
\ar[r]^-{[\A(x,-),1]}
&
\txt{
$[\A_{z,t},$ \\
$[\A_{x,z},\A_{x,u}]]$
}
\ar@{-}[d]_{id}
\ar[r]^-{[1,[\A_{x,y},-]]}
\ar@{}[rrd]|{=}
&
\txt{
$[\A_{z,t},$ \\
$[ [\A_{x,y},\A_{x,z}],$\\ 
$[\A_{x,y},\A_{x,u}] ]]$
}
\ar[r]^-{[1,[\A(x,-),1]]}
&
\txt{
$[\A_{z,t},$\\
$[\A_{y,z},$\\
$[\A_{x,y},\A_{x,u}]]]$
}
\ar@{-}[d]^{id}
\\
\A_{t,u}
\ar[rr]^-{\A(z,-)}
\ar@{-}[d]_{id}
\ar@{}[rrd]|{=}
&
&
[\A_{z,t},\A_{z,u}]
\ar[r]^-{[1,\A(x,-)]}
\ar@{-}[d]_{id}
\ar@<60ex>@{=>}[d]_{[1,\assh_{x,y,z,u}]}
&
\txt{
$[\A_{z,t},$\\
$[\A_{x,z},\A_{x,u}]]$
}
\ar[r]^-{[1, [\A_{x,y},-]]}
&
\txt{
$[ \A_{z,t},$\\
$[ [\A_{x,y},\A_{x,z}],$\\ 
$[\A_{x,y},\A_{x,u}] ] ]$
}
\ar[r]^-{[1,[\A(x,-),1]]}
&
\txt{
$[\A_{z,t},$\\ 
$[\A_{y,z},$\\ 
$[\A_{x,y}, \A_{x,u}]]]$
}
\ar@{-}[d]^{id}
\\
\A_{t,u}
\ar[rr]_-{\A(z,-)}
&
&
[\A_{z,t},\A_{z,u}]
\ar[rr]_-{[1, \A(y,-)]}
&
&
\txt{
$[\A_{z,t},$\\
$[\A_{y,z},\A_{y,u}]]$
}
\ar[r]_-{[1,[1,\A(x,-)]]}
&
\txt{
$[\A_{z,t},$\\ 
$[\A_{y,z},$\\
$[\A_{x,y},\A_{x,u}]]]$
}
}
$$
}
and
{\tiny 
$$
\xymatrix@=1pc@C=4pc{
\A_{t,u}
\ar@{}[drrr]|{=}
\ar[r]^-{\A(x,-)} 
\ar@{-}[d]_{id}
& 
\txt{
$[\A_{x,t},$\\
$\A_{x,u}]$
}
\ar[r]^-{[\A_{x,y},-]}
&
\txt{
$[[\A_{x,y},$\\
$\A_{x,t}],$\\
$[\A_{x,y},$\\
$\A_{x,u}]]$
}
\ar[r]^-{[\A_{y,z},-]}
&
\txt{
$[ [\A_{y,z},$\\
$[\A_{x,y},$\\
$\A_{x,t}] ],$\\  
$[\A_{y,z},$\\
$[\A_{x,y},$\\
$\A_{x,u}]] ]$
}
\ar@<40ex>@{=>}[d]_{[\assh_{x,y,z,t},1]}
\ar[r]^-{[[\A(x,-),1],1]}
\ar@{-}[d]_{id}
& 
\txt{
$[[\A_{x,y},$\\
$\A_{x,z}],$\\
$[\A_{x,y},$\\
$\A_{x,t}]],$\\ 
$[\A_{y,z},$\\
$[\A_{x,y},$\\
$\A_{x,u}]]]$
}
\ar[r]^-{[ [\A_{x,y},-],1] }
& 
\txt{
$[[\A_{x,z},$\\
$\A_{x,t}],$\\ 
$[\A_{y,z},$\\
$[\A_{x,y},$\\
$\A_{x,u}]]]$
}
\ar[d]|{[\A(x,-),1]}
\\
\A_{t,u}
\ar@{-}[d]_{id}
\ar@{}[rrd]|{=} 
\ar[r]^-{\A(x,-)}
& 
\txt{$[\A_{x,t},$\\
 $\A_{x,u}]$}
\ar[r]^-{[\A_{x,y},-]}
&
\txt{
$[[\A_{x,y},$\\
$\A_{x,t}],$\\
$[\A_{x,y},$\\
$\A_{x,u}]]$
}
\ar[r]^-{[\A_{y,z},-]}
\ar@{-}[d]_{id}
\ar@{}[drr]|{=}
&
\txt{
$[
[\A_{y,z},$\\
$[\A_{x,y},$\\
$\A_{x,t}]],$\\ 
$[\A_{y,z},$\\
$[\A_{x,y},$\\
$\A_{x,u}]]]$
}
 \ar[r]^-{[[1,\A(x,-)],1]}
&
\txt{
$[[\A_{y,z},$\\
$\A_{y,t}],$\\
$[\A_{y,z},$\\
$[\A_{x,y},$\\
$\A_{x,u}]] ]$
}
\ar[r]^-{[\A(y,-),1]}
\ar@{-}[d]_{id}
\ar@{}[rd]|{=}
&
\txt{
$[\A_{z,t},$\\ 
$[\A_{y,z},$\\
$[\A_{x,y},$\\
$\A_{x,u}]] ]$
}
\ar@{-}[d]^{id}
\\
\A_{t,u} 
\ar[r]^-{\A(x,-)}
\ar@{-}[d]_{id}
\ar@<48ex>@{=>}[d]^{\assh_{x,y,t,u}}
& 
\txt{
$[\A_{x,t},$\\
$\A_{x,u}]$
}
\ar[r]^-{[\A_{x,y},-]}
&
\txt{
$[[\A_{x,y},$\\
$\A_{x,t}],$\\
$[\A_{x,y},$\\
$\A_{x,u}]]$
}
\ar[r]^-{[\A(x,-),1]}
&
\txt{
$[\A_{y,t},$\\
$[\A_{x,y},$\\
$\A_{x,u}]]$
}
\ar[r]^-{[\A_{y,z},-]}
\ar@{-}[d]_{id}
\ar@{}[rrd]|{=}
&
\txt{
$[[\A_{y,z},$\\
$\A_{y,t}],$\\
$[\A_{y,z},$\\ 
$[\A_{x,y},$\\
$\A_{x,u}]] ]$
}
\ar[r]^-{[\A(y,-),1]}
&
\txt{
$[\A_{z,t},$\\
$[\A_{y,z},$\\
$[\A_{x,y},$\\
$\A_{x,u}]] ]$
}
\ar@{-}[d]^{id}
\\
\A_{t,u}
\ar[rr]^-{\A(y,-)}
\ar@{-}[d]_{id}
\ar@{}[rrd]|{=}
&
&
\txt{
$[\A_{y,t},$\\
$\A_{y,u}]$
}
\ar@{-}[d]_{id}
\ar[r]^-{[1,\A(x,-)]}
\ar@{}[rrd]|{=}
&
\txt{
$[\A_{y,t},$\\
$[\A_{x,y},$\\
$\A_{x,u}]]$\\
}
\ar[r]^-{[\A_{y,z},-]}
&
\txt{
$[[\A_{y,z},$\\
$\A_{y,t}],$\\
$[\A_{y,z},$\\
$[\A_{x,y},$\\
$\A_{x,u}]]]$
}
\ar@{-}[d]_{id}
\ar[r]^-{[\A(y,-),1]}
\ar@{}[rd]|{=}
&
\txt{
$[\A_{z,t},$\\
$[\A_{y,z},$\\
$[\A_{x,y},$\\
$\A_{x,u}]]]$
}
\ar@{-}[d]^{id}
\\
\A_{t,u}
\ar[rr]^-{\A(y,-)}
\ar@{-}[d]_{id}
\ar@{}[rrrd]|{=}
&
&
\txt{
$[\A_{y,t},$\\
$\A_{y,u}]$
}
\ar[r]^-{[\A_{y,z},-]}
&
\txt{
$[[\A_{y,z},$\\ 
$\A_{y,t}],$\\
$[\A_{y,z},$\\
$\A_{y,u}]]$
}
\ar@{-}[d]_{id}
\ar[r]^-{[1,[1,\A(x,-)]]}
\ar@{}[rrd]|{=}
&
\txt{
$[[\A_{y,z},$\\
$\A_{y,t}],$\\ 
$[\A_{y,z},$\\
$[\A_{x,y},$\\
$\A_{x,u}]]]$
}
\ar[r]^-{[\A(y,-),1]}
&
\txt{
$[\A_{z,t},$\\
$[\A_{y,z},$\\
$[\A_{x,y},$\\
$\A_{x,u}]]]$
}
\ar@{-}[d]^{id}
\\
\A_{t,u}
\ar@{-}[d]_{id}
\ar[rr]^-{\A(y,-)}
\ar@<70ex>@{=>}[d]^{\assh_{y,z,t,u}}
&
&
\txt{
$[\A_{y,t},$\\
$\A_{y,u}]$
}
\ar[r]^{[\A_{y,z},-]}
&
\txt{
$[[\A_{y,z},$\\
$\A_{y,t}],$\\
$[\A_{y,z},$\\
$\A_{y,u}]]$
}
\ar[r]^-{[\A(y,-),1]}
&
\txt{
$[\A_{z,t},$\\ 
$[\A_{y,z},$\\
$\A_{y,u}]]$\\
}
\ar@{-}[d]_{id}
\ar[r]^-{[1,[1,\A(x,-)]]}
\ar@{}[rd]|{=}
&
\txt{
$[\A_{z,t},$\\
$[\A_{y,z},$\\
$[\A_{x,y},$\\
$\A_{x,u}]]]$
}
\ar@{-}[d]_{id}
\\
\A_{t,u}
\ar[rr]_-{\A(z,-)}
&
&
\txt{
$[\A_{z,t},$\\
$\A_{z,u}]$
}
\ar[rr]_-{[1,\A(y,-)]}
&
&
\txt{
$[\A_{z,t},$\\
$[\A_{y,z},$\\
$\A_{y,u}]]$
}
\ar[r]_-{[1,[1, \A(x,-)]]}
&
\txt{
$[\A_{z,t},$\\
$[\A_{y,z},$\\
$[\A_{x,y},$\\
$\A_{x,u}]]]$
}
}
$$
}
are equal.
\end{tag}
\begin{tag}\label{cohunassh}
The 2-cell 
$$
\xymatrix{
\A_{y,z} 
\ar[d]_{\A(y,-)}
\ar@{-}[rr]^{id}
&
\ar@{=>}[d]^{\rh}
&
\A_{y,z}
\ar[r]^-{\A(x,-)}
&
[\A_{x,y},\A_{x,z}]
\\
[\A_{y,y},\A_{y,z}] 
\ar[rr]_-{[\unit_y,1]}
&
&
[\unc,\A_{y,z}]
\ar[u]_{\ev_{\gen}}
}
$$
is equal to    
$$\xymatrix{
\A_{y,z}
\ar[d]_{\A(x,-)}
\ar@{-}[r]^{id}
\ar@{}[rd]|{=}
&
\A_{y,z}
\ar@{-}[r]^{id}
\ar[d]^{\A(x,-)}
\ar@{}[rd]|{=}
&
\A_{y,z}
\ar[d]^{\A(x,-)}
\ar@{-}[r]^{id}
&
\A_{y,z}
\ar[d]^{\A(x,-)}
\ar@{}[rd]|{=}
\ar@{-}[r]^{id}
&
\A_{y,z}
\ar[d]^{\A(x,-)}
\\
[\A_{x,y},\A_{x,z}]
\ar@{-}[dddd]_{id}
\ar@{-}[r]^{id}
&
[\A_{x,y},\A_{x,z}]
\ar[d]^{[\ev_{\gen},1]}
\ar@{-}[r]^{id}
&
[\A_{x,y},\A_{x,z}]
\ar[d]_{[\A_{x,y},-]}
\ar@<-5ex>@{=>}[r]^{\assh_{x,y,y,z}}
&
[\A_{y,y},\A_{y,z}]
\ar[dd]^{[1,\A(x,-)]}
\ar@{-}[r]^{id}
&
[\A_{y,y},\A_{y,z}]
\ar[d]^{[\unit_y,1]}
\\
&
\txt{
$[[\unc,\A_{x,y}],$\\
$\A_{x,z}]$
}
\ar[d]^{[[\unit_y,1],1]}
&
\txt{
$[[\A_{x,y},\A_{x,y}],$\\
$[\A_{x,y},\A_{x,z}]]$
}
\ar[d]^{[\A(x,-),1]}
&
&
[\unc,\A_{y,z}]
\ar[d]^{\ev_{\gen}}
\\
\ar@{=>}[r]^-{[\lh_{x,y},1]}
&
\txt{
$[[\A_{y,y},\A_{x,y}],$\\
$\A_{x,z}]$
}
\ar[dd]_{[\A(-,y),1]}
\ar@{}[r]|{=_{(I)}}
&
\txt{
$[\A_{y,y},$\\
$[\A_{x,y},\A_{x,z}]]$
}
\ar[d]^{[\unit_y,1]}
\ar@{}[rdd]|{=}
\ar@{-}[r]^{id}
&
\txt{
$[\A_{y,y},$\\
$[\A_{x,y},\A_{x,z}]]$
}
\ar[d]^{[\unit_y,1]}
\ar@{}[r]|{ \;\;\;=_{(II)} }
&
\A_{y,z}
\ar[dd]^{\A(x,-)}
\\
&
&
\txt{
$[\unc,$\\
$[\A_{x,y},\A_{x,z}]]$
}
\ar[d]^{\ev_{\gen}}
&
\txt{
$[\unc,$\\
$[\A_{x,y},\A_{x,z}]]$
}
\ar[d]^{\ev_{\gen}}
&
\\
[\A_{x,y},\A_{x,z}]
\ar@{-}[r]_{id}
&
[\A_{x,y},\A_{x,z}]
\ar@{-}[r]_{id}
&
[\A_{x,y},\A_{x,z}]
\ar@{-}[r]_{id}
&
[\A_{x,y},\A_{x,z}]
\ar@{-}[r]_{id}
&
[\A_{x,y},\A_{x,z}].
}$$
\end{tag}
Note: Equality $(I)$ in the first of the pastings 
above is established in \ref{ptequal1} in Appendix.
The equality $(II)$ results straightforwardly from 
Lemma \cite{Sch08}-\deuxnateva.\\

Let us justify the equivalence with the previous
definition of $\BS$-category.
For objects $x$, $y$ and $z$ 
the arrows of $\BS$
$$\A(x,-): \A(y,z) \rightarrow [\A(x,y),\A(x,z)]$$
correspond to the bilinear 
$$\compb_{x,y,z}: \A(y,z) \times \A(x,y) \rightarrow \A(x,z)$$
and for objects $x$, the objects $\one_x$ in $\A_{x,x}$ correspond 
to strict the strict arrows
$\unit_x: \unc \rightarrow \A_{x,x}$.
The trilinear 2-cells $\assb_{x,y,z,t}$ correspond to
the 2-cells $\assh_{x,y,z,t}$
and the linear natural transformations
$\rb_{x,y}$ and $\lb_{x,y}$ correspond 
respectively to 2-cells $\rh_{x,y}$
and $\lh_{x,y}$ in $\BS$.\\
For data as above, since the forgetful 
functors $\BS(X,Y) \rightarrow \Cat(X,Y)$
are faithful, the equalities of 2-cells of Axiom \ref{cohassh}
given below 
are equivalent to the equality of natural transformations
of Axioms \ref{cohassb}, and similarly Axiom \ref{cohunassh} 
given below and \ref{cohunassb} are equivalent.\\

Given two $\BS$-categories, a $\BS$-functor 
$\A \rightarrow \B$ consists of 
a map $F$ sending objects of $\A$ to objects of $\B$,
with a collection of arrows in $\BS$ 
$F_{x,y}: \A_{x,y} \rightarrow \B_{Fx,Fy}$ indexed 
by objects $x$ and $y$ of $\A$ and collections of 2-cells 
in $\BS$:\\
- ${F'}^2_{x,y,z}$ indexed by objects $x$, $y$, $z$ 
and as follows 
\begin{tag}\label{Fp2}
{\small
$$
\xymatrix{
\A_{y,z}
\ar[r]^-{F_{y,z}}
\ar[d]_{\A(x,-)}
&
\B_{Fy,Fz}
\ar[r]^-{\B(Fx,-)}
&
[\B_{Fx,Fy}, \B_{Fx,Fz}]
\ar[d]^{[F_{x,y},1]}
\ar@{=>}[lld]_{{F'}^2_{x,y,z}}
\\
[\A_{x,y},\A_{x,z}]
\ar[rr]_{[1,F_{x,z}]}
& & 
[\A_{x,y},\B_{Fx,Fz}]
.}
$$
}
\end{tag}
- $F^0_x$ indexed by objects $x$ as follows 
\begin{tag}\label{F0}
{\small 
$$
\xymatrix{
\unc 
\ar[rr]^{\unit_x} 
\ar[rrdd]_{\unit_{Fx}}
&  
& 
\A_{x,x}
\ar[dd]^{F_{x,x}}
\\
&
\ar@{=>}[ru]^{F^0_x}
&
\\
&
&
\B_{Fx,Fx}
}
$$
}
\end{tag}
Those  satisfy the coherence conditions \ref{cohfun1h}, 
\ref{cohfun2h} and \ref{cohfun3h} below.\\
 
\begin{tag}\label{cohfun1h}
The 2-cells in $\BS$
{\tiny
$$\xymatrix@C=2pc{
& & 
[\B_{Fx,Fz} ,\B_{Fx,Ft}]
\ar[r]^-{[\B_{Fx,Fy},-]}
\ar@<4ex>@{=>}[d]^{\assh_{Fx,Fy,Fz,Ft}}
&
\txt{
$[[\B_{Fx,Fy},\B_{Fx,Fz}],$\\
$[\B_{Fx,Fy},\B_{Fx,Ft}]]$
}
\ar[d]^{[\B(Fx,-),1]}
\\
\A_{z,t} 
\ar[d]_{\A(y,-)}
\ar[r]^-{F_{z,t}}
&
\B_{Fz,Ft}
\ar[ru]^{\B(Fx,-)}
\ar[r]^-{\B(Fy,-)}
&
[\B_{Fy,Fz},\B_{Fy,Ft}]
\ar@{=>}[lld]_{{F'}^2_{y,z,t}}
\ar[r]^-{[1,\B(Fx,-)]}
\ar[d]^{[F_{y,z},1]}
\ar@{}[rd]|{=}
&
\txt{
$[\B_{Fy,Fz},$\\
$[\B_{Fx,Fy},\B_{Fx,Ft}]]$
}
\ar[d]^{[F_{y,z},1]}
\\
[\A_{y,z},\A_{y,t}]
\ar[rr]^-{[1,F_{y,t}]}
\ar[rrd]_{[1,\B(Fx,-)]}
& & 
[\A_{y,z},\B_{Fy,Ft}]
\ar[r]^-{[1,\B(Fx,-)]}
\ar@{=>}[d]^{[1,{F'}^2_{x,y,t}]}
&
\txt{
$[\A_{y,z},$\\
$[\B_{Fx,Fy},\B_{Fx,Ft}]]$
}
\ar[r]^-{[1,[F_{x,y},1]]}
&
\txt{
$[\A_{y,z},$\\
$[\A_{x,y},\B_{Fx,Ft}]]$
}\\
& & 
\txt{
$[\A_{y,z},$\\
$[\A_{x,y},\A_{x,t}]]$
}
\ar[rru]_{[1,[1,F_{x,t}]]}
}$$
}
and
{\tiny 
$$\xymatrix@C=2pc{
\B_{Fz,Ft}
\ar[r]^-{\B(Fx,-)}
&
[\B_{Fx,Fz},\B_{Fx,Ft}]
\ar[d]^{[F_{x,z},1]}
\ar[r]^-{[\A_{x,y},-]}
\ar@{=>}[ddl]_{{F'}^2_{x,z,t}}
\ar@{}[rd]|{=}
&
\txt{
$[[\A_{x,y},\B_{Fx,Fz}],$\\
$[\A_{x,y},\B_{Fx,Ft}]]$
}
\ar[r]^-{[[F_{x,y},1],1]}
\ar[d]^{[[1,F_{x,z}],1]}
&
\txt{
$[[\B_{Fx,Fy},\B_{Fx,Fz}],$\\
$[\A_{x,y},\B_{Fx,Ft}]]$
}
\ar[r]^-{[\B(Fx,-),1]}
\ar@{=>}[d]^{[{F'}^2_{x,y,z},1]}
&
\txt{
$[\B_{Fy,Fz},$\\
$[\A_{x,y},\B_{Fx,Ft}]]$
}
\ar[d]^{[F_{y,z},1]}
\\
& 
[\A_{y,z},\B_{Fx,Ft}]
\ar[r]^-{[\A_{x,y},-]}
\ar@{}[rd]|{=}
&
\txt{
$[[\A_{x,y},\A_{x,z}],$\\
$[\A_{x,y},\B_{Fx,Ft}]]$
}
\ar[rr]_-{[\A(x,-),1]}
& 
& 
\txt{
$[\A_{y,z},$\\
$[\A_{x,y},\B_{Fx,Ft}]]$
}
\\
\A_{z,t}
\ar[r]^-{\A(x,-)}
\ar[uu]^{F_{z,t}}
\ar[drr]_{\A(y,-)}
&
[\A_{x,z},\A_{x,t}]
\ar[r]^-{[\A_{x,y},-]}
\ar[u]_{[1,F_{x,t}]}
&
\txt{
$[[\A_{x,y},\A_{x,z}],$\\
$[\A_{x,y},\A_{x,t}]]$
}
\ar[rr]^{[\A(x,-),1]}
\ar[u]_{[1,[1,F_{x,t}]]}
\ar@{}[rru]|{=}
\ar@{=>}[d]|{\assh_{x,y,z,t}}
&
&
\txt{
$[\A_{y,z},$\\
$[\A_{x,y},\A_{x,t}]]$
}
\ar[u]_{[1,[1,F_{x,t}]]}
\\
&
&
[\A_{y,z},\A_{y,t}]
\ar[rru]_{[1,\A(x,-)]}
}
$$
}
are equal.
\end{tag}

\begin{tag}\label{cohfun2h}
The 2-cells
{\tiny 
$$
\xymatrix{
\A_{x,y}
\ar[d]_{\A(x,-)} 
\ar[rr]^{F_{x,y}} 
& &  
\B_{Fx,Fy}
\ar[d]|{\B(Fx,-)}
\ar@{=>}[lld]_{{F'}^2_{x,x,y}}
\ar@/^70pt/[dddd]^1
\\
[\A_{x,x},\A_{x,y}]
\ar[dd]_{[\unit_x,1]}
\ar[rd]^{[1,F_{x,y}]}
&
&
[\B_{Fx,Fx},\B_{Fx,Fy}]
\ar[dd]^>>{[\unit_{Fx},1]}
\ar[ld]_{[F_{x,x},1]}
\\
&
[\A_{x,x},\B_{Fx,Fy}]
\ar[rd]_{[\unit_x,1]}
&
\ar@{=>}[l]_-{[F^0_x,1]}
& 
\ar@{=>}[l]_{\rh_{Fx,Fy}}
\\
[\unc,\A_{x,y}]
\ar@{}[rrd]|{=}
\ar@{}[ru]|{=}
\ar[d]_{\ev_{\gen}}
\ar[rr]^{[1,F_{x,y}]}
& & 
[\unc,\B_{Fx,Fy}]
\ar[d]^{\ev_{\gen}}
\\
\A_{x,y}
\ar[rr]_{F_{x,y}}
& & 
\B_{Fx,Fy}
}
$$
}
and
{\tiny
$$
\xymatrix{
\A_{x,y}
\ar@{-}[rr]^{id}
\ar[d]_{\A(x,-)}
&
\ar@{=>}[d]^{\rh_{x,y}}
&
\A_{x,y}
\ar[r]^{F_{x,y}}
& 
\B_{Fx,Fy}
\\
[\A_{x,x}, \A_{x,y}]
\ar[rr]_{[\unit_x,1]}
& 
& 
[\unc,\A_{x,y}]
\ar[u]_{\ev_{\gen}}
}
$$
}
are equal.
\end{tag}
\begin{tag}\label{cohfun3h}
The 2-cells in $\BS$
{\tiny 
$$
\xymatrix{
\A_{x,y}
\ar[d]_{\A(-,y)} 
\ar[rr]^{F_{x,y}} 
& &  
\B_{Fx,Fy}
\ar[d]|{\B(-,Fy)}
\ar@{=>}[lld]_{{({F'}^2_{x,y,y})}^*}
\ar@/^70pt/[dddd]^1
\\
[\A_{y,y},\A_{x,y}]
\ar[dd]_{[\unit_y,1]}
\ar[rd]^{[1,F_{x,y}]}
&
&
[\B_{Fy,Fy},\B_{Fx,Fy}]
\ar[dd]^>>{[\unit_{Fy},1]}
\ar[ld]_{[F_{x,y},1]}
\\
&
[\A_{y,y},\B_{Fx,Fy}]
\ar[rd]_{[\unit_y,1]}
&
\ar@{=>}[l]_-{[F^0_y,1]}
& 
\ar@{=>}[l]_{\lh_{Fx,Fy}}
\\
[\unc,\A_{x,y}]
\ar@{}[rrd]|{=}
\ar@{}[ru]|{=}
\ar[d]_{\ev_{\gen}}
\ar[rr]^{[1,F_{x,y}]}
& & 
[\unc,\B_{Fx,Fy}]
\ar[d]^{\ev_{\gen}}
\\
\A_{x,y}
\ar[rr]_{F_{x,y}}
& & 
\B_{Fx,Fy}
}
$$
}
and
{\tiny
$$
\xymatrix{
\A_{x,y}
\ar@{-}[rr]^{id}
\ar[d]_{\A(-,y)}
&
\ar@{=>}[d]^{\lh_{x,y}}
&
\A_{x,y}
\ar[r]^{F_{x,y}}
& 
\B_{Fx,Fy}
\\
[\A_{y,y}, \A_{x,y}]
\ar[rr]_{[\unit_y,1]}
& 
& 
[\unc,\A_{x,y}]
\ar[u]_{\ev_{\gen}}
}
$$
}
are equal.
\end{tag}

One can also define the $\BS$-natural transformations,
in a similar way. For this purpose, we need some notation.\\

Let $\A$ be an arbitrary $\BS$-category.
To give a {\em strict} arrow 
$\unc \rightarrow \A(x,y)$ is equivalent 
to give an arrow $x \rightarrow y$ of 
the underlying bicategory $\A^0$ and 
we might confuse the two.
We therefore define for any object $z$ of $\A$
and any arrow
$f: x \rightarrow y$ of $\A^0$ the 
arrows
\begin{tag}\label{Af1} 
$\A(f,1)$ as the composite arrow in $\BS$
$$
\xymatrix{
\A_{x,y} 
\ar[r]^-{\A(x,-)}
&
[\A_{x,y},\A_{x,z}]
\ar[r]^-{[F,1]}
&
[\unc,\A_{x,z}]
\ar[r]^-{\ev_{\gen}}
&
\A_{x,z}
}
$$
\end{tag}
and
\begin{tag}\label{A1f}
$\A(1,f)$ as 
$$
\xymatrix{
\A_{z,x} 
\ar[r]^-{\A(-,y)}
&
[\A_{x,y}, \A_{z,y}]
\ar[r]^-{[F,1]}
&
[\unc,\A_{z,y}]
\ar[r]^{\ev_{\gen}}
&
\A_{z,y}.
}
$$
\end{tag}

Given any arrows 
$\xymatrix{
x
\ar[r]^f 
& 
y 
\ar[r]^g
& 
z
\ar[r]^h 
& t}$ 
in $\A^0$, we define the two cells
\begin{tag}
${\TCcompu}_{f,y,z}$ 
$$
\xymatrix@C=3pc{
\A_{z,t} 
\ar[r]^-{\A(y,-)}
\ar[d]_{\A(x,-)}
&
[\A_{y,z},\A_{y,t}]
\ar[d]^{[1,\A(f,1)]}
\\
[\A_{x,z},\A_{x,t}]
\ar[r]_-{[\A(f,1),1]}
\ar@{=>}[ru]
&
[\A_{y,z},\A_{x,t}]
}
$$
\end{tag}
\begin{tag}
$\TCcompd_{x,g,t}$ 
$$
\xymatrix{
\A_{z,t}
\ar[d]_{\A(g,1)}
\ar[r]^-{\A(x,-)}
&
[\A_{x,z}, \A_{x,t}]
\ar[d]^{[\A(1,g),1]}
\ar@{=>}[ld]
\\
\A_{y,t}
\ar[r]_-{\A(x,-)}
&
[\A_{x,y}, \A_{x,t}]
}
$$
\end{tag}
and eventually
\begin{tag}
$\TCcompt_{x,y,h}$ 
$$
\xymatrix@C=3pc{
\A_{y,z}
\ar[r]^-{\A(1,h)}
\ar[d]_{\A(x,-)}
&
\A_{y,t}
\ar[d]^{\A(x,-)}
\\
[\A_{x,y},\A_{x,z}]
\ar[r]_-{[1,\A(1,h)]}
\ar@{=>}[ru]
&
[\A_{x,y},\A_{x,t}]
}
$$
\end{tag}
which are obtained 
from the trilinear natural transformation \ref{assb}
$${(\assb_{x,y,z,t})}_{h,g,f}: h \circ (g \circ f) 
\rightarrow (h \circ g) \circ f$$
by fixing one of its argument. For $\TCcompu$, $\TCcompd$ and $\TCcompt$
fix respectively $f$,$g$ and $h$.\\

Formally $\TCcompu_{f,z,t}$ is the composite 
$$\xymatrix{
\A_{z,t} \ar@{=>}[r]^-{\assh_{x,y,z,t}} 
&
[\A_{y,z},[\A_{x,y}, \A_{x,t}]]
\ar[r]^-{[1,[f,1]]}
&
[\A_{y,z}, [\unc, \A_{x,t}]]
\ar[r]^-{[1,\ev_{\gen}]}
&
[\A_{y,z},\A_{x,t}]}$$
(see \ref{defTCcomp1fyz} in Appendix), 
$\TCcompd_{x,g,t}$ it is the composite 
$$
\xymatrix{
\A_{z,t}
\ar[r]^-{\assh_{x,y,z,t}}
&
[\A_{y,z},[\A_{x,y},\A_{x,t}]]
\ar[r]^{[g,1]}
& 
[\unc,[\A_{x,y},\A_{x,t}]]
\ar[r]^{\ev_{\gen}}
&
[\A_{x,y},\A_{x,t}]
}
$$
and $\TCcompt_{x,y,h}$ is the image by $\ev_{\gen}$ of the composite
$$
\xymatrix{
\unc 
\ar[r]^-h
&
\A_{z,t}
\ar@{=>}[r]^-{\assh_{x,y,z,t}}
&
[\A_{y,z},[\A_{x,y},\A_{x,t}]]
}
$$
(see \ref{defTCcomp3xyh} in Appendix).
By definition all 2-cells $\TCcompu$, $\TCcompd$ and $\TCcompt$ are
identities when $\A$ is strict.\\

One has also for any objects $x$ and $y$ 
of $\A$ the 2-cell in $\BS$
\begin{tag}
$\rhp_{x,y}:$
{\small
$$
\xymatrix{
\unc 
\ar[rr]^-{\unit_x}
\ar[rrdd]_{\vv} 
& & 
\A_{x,x}
\ar[dd]^{\A(-,y)}\\
& \ar@{=>}[ru] & \\
& & [\A_{x,y},\A_{x,y}]},
$$
}
\end{tag}
which according to Remarks \ref{evgen2cel}
is determined by its value in $\gen$ 
which is the linear natural transformation
$${ \rb_{x,y} }_f: f \circ \one_x \rightarrow f$$
of \ref{rb},
and corresponds by the bijection \ref{2c2}/\ref{2c4} in Appendix
to the 2-cell $\rh_{x,y}$
{\small
$$
\xymatrix{
\A_{x,y} 
\ar[rr]^-{\A(x,-)}
\ar@{-}[dd]_{id}
&
&
[\A_{x,x}, \A_{x,y}]
\ar[dd]^{[\unit_x,1]}
\\
\ar@{=>}[rr]
& & 
\\
\A_{x,y}
&
& 
[\unc,\A_{x,y}]
\ar[ll]^{\ev_{\gen}}.
}
$$
}
Similarly for any objects $x$ and $y$ of $\A$
one has the 2-cell in $\BS$
\begin{tag}
$\lhp_{x,y}:$
{\small
$$\xymatrix{
\unc
\ar[rr]^{\unit_y}
\ar[rrdd]_{\vv}
&
&
\A_{y,y}
\ar[dd]^{\A(x,-)}
\\
& 
\ar@{=>}[ru]
& 
\\
& & 
[\A_{x,y},\A_{x,y}]
}$$
}
\end{tag}
that corresponds to the linear natural transformation
$${(\lb_{x,y})}_f : \one_y \circ f \rightarrow f$$ 
of \ref{lb} and which corresponds by the bijection 
\ref{2c2}/\ref{2c4} in Appendix to the 2-cell $\lh_{x,y}$
{\small
$$
\xymatrix{
\A_{x,y} 
\ar[rr]^{\A(-,y)}
\ar@{-}[dd]_{id}
&
&
[\A_{y,y}, \A_{x,y}]
\ar[dd]^{[\unit_y,1]}
\\
\ar@{=>}[rr]& & 
\\
\A_{x,y} 
&
&
[\unc,\A_{x,y}].
\ar[ll]^{\ev_{\gen}} 
}
$$
}

Given any strict arrow $\tilde{f}: \unc \rightarrow \A_{x,y}$,
with corresponding arrow $f: x \rightarrow y$ in $\A^0$, 
one has the 2-cell in $\BS$
\begin{tag}\label{defrhpp}
$\rhpp_f:$
{\small
$$
\xymatrix{
\unc 
\ar[rr]^-{\unit_y}
\ar[rrdd]_{\tilde{f}}
&
&
\A_{y,y}
\ar[dd]^{\A(f,1)}
\\
& \ar@{=>}[ru] &  
\\
&
& 
[\A_{x,y},\A_{x,y}]
}
$$ 
}
\end{tag}
which corresponds to the {\em 2-cell} 
$\rb_f : f \rightarrow \one_y \circ f: x \rightarrow y$ 
of $\A^0$. It is the pasting 
{\small
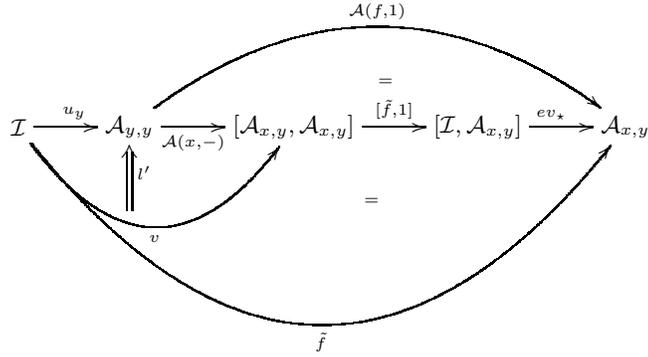
$$
\xymatrix{
& & \ar@{}[rd]|{=}
\\
\unc
\ar[r]^{\unit_y}
\ar@/_38pt/[rr]_-{\vv}
\ar@/_75pt/[rrrr]_-{\tilde{f}}
&
\A_{y,y}
\ar[r]_-{\A(x,-)}
\ar@/^38pt/[rrr]^-{\A(f,1)}
&
[\A_{x,y}, \A_{x,y}]
\ar[r]^-{[\tilde{f},1]}
&
[\unc,\A_{x,y}]
\ar[r]^{\ev_{\gen}}
&
\A_{x,y}
\\
& \ar@{=>}[u]_{\lhp} & 
\\
& & \ar@{}[uur]^{=}
}
$$
}
where the bottom identity 2-cell above is established
in Lemma \ref{evstrict} in Appendix.
Similarly one has the 2-cell in $\BS$
\begin{tag}\label{deflhpp}
{\small
$\lhpp_f:$
$$
\xymatrix{
\unc 
\ar[rr]^-{\unit_x}
\ar@{-}[rrdd]_{\tilde{f}}
&
&
\A_{x,x}
\ar[dd]^{\A(1,f)}
\\
& \ar@{=>}[ru] &  
\\
&
& 
[\A_{x,y},\A_{x,y}]
}
$$ 
}
\end{tag} 
that corresponds the 2-cell
$$ {( {\lb}_{x,y} ) }_f:  f \rightarrow  f \circ \one_x $$
and is the pasting 
{\small 
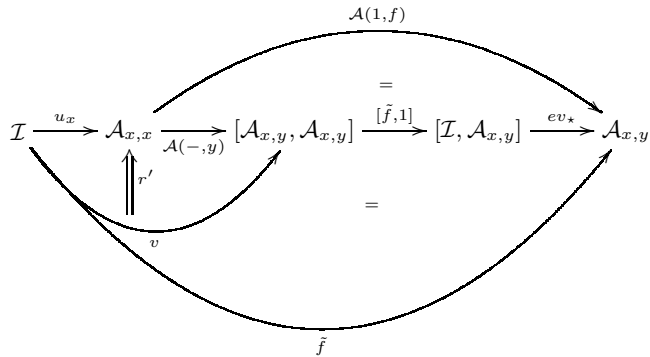
$$
\xymatrix{
& & \ar@{}[rd]|{=}
\\
\unc
\ar[r]^-{\unit_x}
\ar@/_38pt/[rr]_-{\vv}
\ar@/_75pt/[rrrr]_-{\tilde{f}}
&
\A_{x,x}
\ar[r]_-{\A(-,y)}
\ar@/^38pt/[rrr]^-{\A(1,f)}
&
[\A_{x,y}, \A_{x,y}]
\ar[r]^-{[\tilde{f},1]}
&
[\unc,\A_{x,y}]
\ar[r]^-{\ev_{\gen}}
&
\A_{x,y}
\\
& \ar@{=>}[u]_{\rhp} & 
\\
& & \ar@{}[uur]^{=}}
$$
}

Given two $\BS$-functors $F,G: \A \rightarrow \B$,
a {\em $\BS$-natural transformation}
$(\sigma, \kappa): F \rightarrow G: \A \rightarrow \B$
consists of a collection of {\em strict} arrows 
$\sigma_x: \unc \rightarrow \B(Fx,Gx)$
(or 1-cells $\sigma_x : Fx \rightarrow Gx$ in $\B_0$),  
indexed by objects $x$ of $\A$
together with a collection of 2-cells $\kappa_{x,y}$ in $\BS$
for objects $x$,$y$ of $\A$ as follows
\begin{tag}\label{kappa}
$$
\xymatrix{
\A_{x,y}
\ar[r]^{F_{x,y}}
\ar[d]_{G_{x,y}}
& 
\B_{Fx,Fy}
\ar[d]^{\B(1,\sigma_y)}
\\
\B_{Gx,Gy}
\ar@{=>}[ru]
\ar[r]_{\B(\sigma_x,1)}
& 
\B_{Fx,Gy}
}
$$
\end{tag}
and that satisfies the two coherence conditions 
\ref{cohnat1}, \ref{cohnat2} and
below. 
\begin{tag}\label{cohnat1}
For any object $x$, $y$ and $z$ in $\A$, the 2-cells
$\Xi_1$, $\Xi_2$, $\Xi_3$, $\Xi_4$, $\Xi_5$, $\Xi_6$, $\Xi_7$ and $\Xi_8$  
below satisfy the equality
$$\Xi_2 \circ \Xi_1 = 
\Xi_8 \circ {(\Xi_7)}^{-1} \circ \Xi_6 \circ \Xi_5 \circ \Xi_4 \circ 
{(\Xi_3)}^{-1}.$$
$\Xi_1$
is 
$$\xymatrix{
\A_{y,z}
\ar[d]_-{\A(x,-)}
\ar[r]^-{G_{y,z}}
& 
\B_{Gy,Gz}
\ar@{=>}[d]^{{G'}^2_{x,y,z}}
\ar[r]^-{\B(Gx,-)}
&
[\B_{Gx,Gy},\B_{Gx,Gz}]
\ar[d]^{[G_{x,y},1]}
\\
[\A_{x,y},\A_{x,z}]
\ar[rr]_-{[1,G_{x,z}]}
& & 
[\A_{x,y},\B_{Gx,Gz}]
\ar[r]_-{[1,\B(\sigma_x,1)]}
&
[\A_{x,y},\B_{Fx,Gz}]
}
$$
$\Xi_2$ is
$$
\xymatrix{
& & [\A_{x,y},\B_{Gx,Gz}]
\ar[rd]^{[1, \B(\sigma_x,1)]}
\ar@{=>}[dd]^{[1,\kappa_{x,z}]}
\\
\A_{y,z} 
\ar[r]^-{\A(x,-)}
&
[\A_{x,y},\A_{x,z}]
\ar[ru]^{[1,G_{x,z} ]}
\ar[rd]_{[1,F_{x,z}] }
& & 
[\A_{x,y},\B_{Fx,Gz}]
\\
& & 
[\A_{x,y}, \B_{Fx,Fz}]
\ar[ru]_{[1, \B(1,\sigma_z)]}
}
$$
$\Xi_3$ is
$$
\xymatrix{
& & 
[\B_{Gx,Gy}, \B_{Gx,Gz}]
\ar[rd]^{[1,\B(\sigma_x,1)]}
\\
\A_{y,z}
\ar[r]^-{G_{y,z}}
&
\B_{Gy,Gz}
\ar[ru]^{\B(Gx,-)}
\ar[rd]_{\B(Fx,-)}
&
&
[\B_{Gx,Gy}, \B_{Fx,Gz}]
\ar[r]^-{[G_{x,y},1]}
&
[\A_{x,y}, \B_{Fx,Gz}]
\\
& 
& 
[\B_{Fx,Gx}, \B_{Fx,Gz}]
\ar[ru]_{[\B(\sigma_x,1),1]}
\ar@{=>}[uu]|{ \TCcompu_{ \B(\sigma_x,1) , y, z } }
}
$$
$\Xi_4$ is 
$$
\xymatrix{
& & & [\B_{Gx,Gy}, \B_{Fx,Gz}]
\ar@{=>}[dd]^{[\kappa_{x,y},1]}
\ar[rd]^{[G_{x,y},1]}
\\
\A_{y,z}
\ar[r]^-{G_{y,z}}
&
\B_{Gy,Gz}
\ar[r]^-{\B(Fx,-)}
&
[\B_{Fx,Gy }, \B_{Fx,Gz }]
\ar[ru]^{[\B(\sigma_x,1),1]}
\ar[rd]_{[\B(1, \sigma_y),1]}
& & 
[\A_{x,y}, \B_{Fx,Gz}]
\\
& & & 
[\B_{Fx,Fy},\B_{Fx,Gz}]
\ar[ru]_{[F_{x,y},1]}
}$$
$\Xi_5$ is
$$
\xymatrix{
& & 
[\B_{Fx,Gy}, \B_{Fx,Gz}]
\ar[rd]^{[\B(1,\sigma_y),1]}
\ar@{=>}[dd]|{\TCcompd_{Fx, \sigma_y, Gz}}
\\
\A_{y,z}
\ar[r]^-{G_{y,z}}
&
\B_{Gy,Gz}
\ar[rd]_{\B(\sigma_y,1)}
\ar[ru]^{\B(Fx,-)}
&
&
[\B_{Fx,Fy},\B_{Fx,Gz}]
\ar[r]^-{[F_{x,y},1]}
&
[\A_{x,y}, \B_{Fx,Gz}] 
\\
& & 
\B_{Fy,Gz}
\ar[ru]_{\B(Fx,-)}
}
$$
$\Xi_6$ is
$$
\xymatrix{
& 
\B_{Gy,Gz}
\ar[rd]^{\B(\sigma_y,1)}
\ar@{=>}[dd]^{\kappa_{y,z}}
&
\\
\A_{y,z}
\ar[ru]^{G_{y,z}}
\ar[rd]_{F_{y,z}}
& & 
\B_{Fy,Gz}
\ar[r]^-{\B(Fx,-)}
&
[\B_{Fx,Fy},\B_{Fx,Gz}]
\ar[r]^-{[F_{x,y},1]}
&
[\A_{x,y},\B_{Fx,Gz}]
\\
& 
\B_{Fy,Fz}
\ar[ru]_{\B(1,\sigma_z)}
}
$$
$\Xi_7$ is
$$\xymatrix{
& & \B_{Fy,Gz}
\ar[rd]^{\B(Fx,-)}
\\
\A_{y,z}
\ar[r]^-{F_{y,z}}
&
\B_{Fy,Fz}
\ar[ru]^{\B(1,\sigma_z)}
\ar[rd]_{\B(Fx,-)}
& & 
[\B_{Fx,Fy},\B_{Fx,Gz}]
\ar[r]^-{[F_{x,y},1]}
&
[\A_{x,y}, \B_{Fx,Gz}]
\\
& & 
[\B_{Fx,Fy}, \B_{Fx,Fz}]
\ar[ru]_{[1,\B(1,\sigma_z)]}
\ar@{=>}[uu]|{\TCcompt_{Fx,Fy,\sigma_z}}
}
$$
$\Xi_8$ is
$$
\xymatrix{
\A_{y,z}
\ar[r]^-{F_{y,z}}
\ar[d]_{\A(x,-)}
&
\B_{Fy,Fz}
\ar[r]^-{\B(Fx,-)}
\ar@{=>}[d]^{{F'}^2_{x,y,z}}
&
[\B_{Fx,Fy},\B_{Fx,Fz}]
\ar[d]^{[F_{x,y},1]}
\\
[\A_{x,y},\A_{x,z}]
\ar[rr]_{[1,F]}
& & 
[\A_{x,y}, \B_{Fx,Fz}]
\ar[r]^-{[1,\B(1,\sigma_z)]}
&
[\A_{x,y}, \B_{Fx,Gz}]
}
$$
\end{tag}

\begin{tag}\label{cohnat2}
For any object $x$ of $\A$, the 2-cells
$$
\xymatrix{
& 
& 
\ar@{=>}[dd]|{\lhpp_{\sigma_x}}
& 
\\
& 
\ar@{=>}[d]^{F^0}
& & 
\\
\unc
\ar@/^38pt/[rr]^-{\unit}
\ar[r]^-{\unit}
\ar@/^70pt/[rrr]^-{\sigma_x}
&
\A_{x,x}
\ar[r]_-{F_{x,x}}
&
\B_{Fx,Fx}
\ar[r]_-{\B(1,\sigma_x)}
&
\B_{Fx,Gx}
}
$$
and 
$$
\xymatrix{
& 
&
\ar@{=>}[d]|{\rhpp_{\sigma_x}}
&
& 
\\
\unc 
\ar@{-}[d]_{id}
\ar[rr]^-{\unit_{Gx}}
\ar@/^38pt/[rrrr]^{\sigma_x}
&
\ar@{=>}[d]^{G^0}
& 
\B_{Gx,Gx}
\ar[rr]^-{\B(\sigma_x,1)}
\ar@{=>}[rrd]^{\kappa}
&
&
\B_{Fx,Gx}
\\
\unc
\ar[rr]_-{\unit_x}
& 
&
\A_{x,x}
\ar[u]^{G_{x,x}}
\ar[rr]_-{F_{x,x}}
&
&
\B_{Fx,Fx} 
\ar[u]_{\B(1,\sigma_x)}
}
$$
are equal.
\end{tag}

\end{section}

\begin{section}{$\BS$-categories via the tensor}\label{BSCattens}
In this section we give definitions of $\BS$-categories 
and $\BS$-functors that rely on the tensor product in $\BS$.\\

A {\em $\BS$-category} $(\A, \unit, \comp, \assp, \rp, \lp)$
consists of the following data:\\
- As before: a small set of objects with a map sending any 
pair $x$,$y$ of objects to an object $\A_{x,y}$ of $\BS$;\\
- Collections of {\em strict}
morphisms $\unit_x : \unc \rightarrow \A_{x,x}$ and
$\comp_{x,y,z}: \A_{y,z} \otimes \A_{x,y} \rightarrow \A_{x,z}$ 
indexed by objects of $\A$ with collections of
2-cells $\assp_{x,y,z,t}$, $\rp_{x,y}$ and $\lp_{x,y}$ in $\BS$ 
indexed by objects of $\A$ and as follows
\begin{tag}\label{assp}
$$\xymatrix{
(\A_{z,t} \otimes \A_{y,z}) \otimes \A_{x,y} 
\ar[rr]^{A'}
\ar[d]_{\comp_{y,z,t} \otimes 1}
& 
& 
\A_{z,t} \otimes (\A_{y,z} \otimes \A_{x,y})
\ar[d]^{1 \otimes \comp_{x,y,z}}
\ar@{=>}[lld]_{\assp_{x,y,z,t}}
&
\\
\A_{y,t} \otimes \A_{x,y} 
\ar[rd]_{\comp_{x,y,t}}
& 
& 
\A_{z,t} \otimes \A_{x,z}
\ar[ld]^{\comp_{x,z,t}}
\\
& 
\A_{x,t}
&
}
$$
\end{tag}
\begin{tag}\label{rp}
$$\xymatrix{
\A_{x,y} 
\ar[r]^{R'}
\ar@{-}[dd]_{id}
& 
\A_{x,y} \otimes \unc
\ar[dd]^{1 \otimes \unit_x}
\\
\ar@{=>}[r]^-{\rp_{x,y}}
&
\\
\A_{x,y}
& 
\A_{x,y} \otimes \A_{x,x} 
\ar[l]^-{\comp_{x,x,y}}
}
$$
\end{tag}
and
\begin{tag}\label{lp}
$$
\xymatrix{
\A_{x,y} 
\ar[r]^-{L'}
\ar@{-}[dd]_{id}
& 
\unc \otimes \A_{x,y}
\ar[dd]^{\unit_y \otimes 1}
\\
\ar@{=>}[r]^-{\lp_{x,y}}
&
\\
\A_{x,y}
& 
\A_{y,y} \otimes \A_{x,y} 
\ar[l]^-{\comp_{x,y,y}}
}
$$
\end{tag}
Those satisfy the coherence Axioms \ref{cohassp} and
\ref{cohunassp} below.
\begin{tag}\label{cohassp}
For any objects $x$,$y$,$z$,$t$ and $u$ of $\A$, the 2-cells
$$ 
\xymatrix@C=1pc{
((\A_{t,u} \A_{z,t}) \A_{y,z})  \A_{x,y} 
\ar[r]^{A'}
\ar[d]_{(\comp_{z,t,u} \otimes 1) \otimes 1}
\ar@{}[rd]|{=}
&
(\A_{t,u}  \A_{z,t})  (\A_{y,z}  \A_{x,y})
\ar@{-}[r]^{id}
\ar[d]^{\comp_{z,t,u} \otimes 1}
\ar@{}[rdd]|{=}
&
(\A_{t,u}  \A_{z,t}) (\A_{y,z} \A_{x,y})
\ar[r]^{A'}
\ar[d]^{1 \otimes \comp_{x,y,z}}
\ar@{}[rd]|{=}
&
\A_{t,u}  ( \A_{z,t} (\A_{y,z} \A_{x,y}))
\ar[d]^{1 \otimes (1 \otimes \comp_{x,y,z})}
\\
(\A_{z,u}  \A_{y,z}) \A_{x,y}
\ar[r]^{A'}
\ar[d]_{\comp_{y,z,t} \otimes 1}
&
\A_{z,u} (\A_{y,z}  \A_{x,y})
\ar[d]^{1 \otimes \comp_{x,y,z}}
\ar@{=>}[ldd]|{\assp_{x,y,z,u}}
&
(\A_{t,u} \A_{z,t}) \A_{x,z}
\ar[r]^{A'}
\ar[d]^{\comp_{z,t,u} \otimes 1}
& 
\A_{t,u} (\A_{z,t} \A_{x,y})
\ar[d]^{1 \otimes \comp_{x,y,t}}
\ar@{=>}[ldd]|{\assp_{x,z,t,u}}
\\
\A_{y,u}  \A_{x,y} 
\ar[d]_{\comp_{x,y,u}}
&
\A_{z,u} \A_{x,z}
\ar[d]_{\comp_{x,z,u}} 
\ar@{-}[r]^{id}
\ar@{}[rd]|{=}
& 
\A_{z,u} \A_{x,z}
\ar[d]_{\comp_{x,z,u}}
&
\A_{t,u} \A_{x,t}
\ar[d]^{\comp_{x,t,u}}\\
\A_{x,u}
\ar@{-}[r]_{id}
&
\A_{x,u}
\ar@{-}[r]_{id}
&
\A_{x,u}
\ar@{-}[r]_{id}
&
\A_{x,u}
}
$$
and 
$$
\xymatrix@C=1pc{
((\A_{t,u}  \A_{z,t})  \A_{y,z})  \A_{x,y}
\ar[r]^{A' \otimes 1}
\ar[d]_{(\comp_{z,t,u} \otimes 1) \otimes 1}
&
((\A_{t,u}  (\A_{z,t}  \A_{y,z}))  \A_{x,y}
\ar[r]^{A'}
\ar[d]|{(1 \otimes \comp_{y,z,t}) \otimes 1}
\ar@{}[rd]|{=}
\ar@{=>}[ldd]|{\assp_{y,z,t,u} \otimes 1}
&
\A_{t,u} ((\A_{z,t} \A_{y,z})  \A_{x,y})
\ar[r]^{1 \otimes A'}
\ar[d]|{1 \otimes (\comp_{y,z,t} \otimes 1)}
& 
\A_{t,u} ( \A_{z,t} (\A_{y,z} \A_{x,y}))
\ar[d]^{1 \otimes (1 \otimes \comp_{x,y,z})}
\ar@{=>}[ldd]|{1 \otimes \assp_{x,y,z,t}}
\\
(\A_{z,u}  \A_{y,z}) \A_{x,y}
\ar[d]_{\comp_{y,z,u} \otimes 1}
&
(\A_{t,u}  \A_{y,t}) \A_{x,y}
\ar[d]|{\comp_{y,t,u} \otimes 1}
\ar[r]^{A'}
&
\A_{t,u} (\A_{y,t} \A_{x,y})
\ar[d]|{1 \otimes \comp_{x,y,t}}
\ar@{=>}[ldd]|{\assp_{x,y,t,u}}
&
\A_{t,u}  (\A_{z,t} \A_{x,z})
\ar[d]^{1 \otimes \comp_{x,z,t}}
\\
\A_{y,u} \A_{x,y}
\ar[d]_{\comp_{x,y,u}}
\ar@{-}[r]^{id}
\ar@{}[rd]|{=}
&
\A_{y,u} \A_{x,y}
\ar[d]_{\comp_{x,y,u}}
& 
\A_{t,u} \A_{x,t}
\ar@{-}[r]^{id}
\ar[d]_{\comp_{x,t,u}}
\ar@{}[rd]|{=}
&
\A_{t,u}  \A_{x,t}
\ar[d]^{\comp_{x,t,u}}
\\
\A_{x,u}
\ar@{-}[r]_{id}
&
\A_{x,u}
\ar@{-}[r]_{id}
&
\A_{x,u}
\ar@{-}[r]_{id}
&
\A_{x,u}
}
$$
are equal.
Note that the domains of the above 2-cells are equal 
since $(1 \otimes A') \circ A' \circ (A' \otimes 1) = A' \circ A'$
by Lemma \cite{Sch08}-\waxiomdouze.
\end{tag}
\begin{tag}\label{cohunassp}
For any objects $x$, $y$, $z$ of $\A$,
the 2-cells\\  
$\Xi_1$ $=$
$$\xymatrix{
\A_{y,z} \otimes \A_{x,y}
\ar[d]_{R' \otimes 1} 
\ar@{-}[rr]^{id}
&
\ar@{=>}[d]^{\rp_{y,z} \otimes 1}
&
\A_{y,z} \otimes \A_{x,y}
\ar[r]^-{\comp_{x,y,z}}
&
\A_{x,z}
\\
(\A_{y,z} \otimes \unc) \otimes \A_{x,y}
\ar[rr]_{(1 \otimes \unit_y) \otimes 1}
&
&
(\A_{y,z} \otimes \A_{y,y}) \otimes \A_{x,y}
\ar[u]_{\comp_{y,y,z} \otimes 1} 
}
$$\\
$\Xi_2$ $=$ $$\xymatrix{
\A_{y,z} \otimes \A_{x,y}
\ar[d]_{1 \otimes L'} 
\ar@{-}[rr]^{id}
&
\ar@{=>}[d]^{1 \otimes \lp_{x,y}}
&
\A_{y,z} \otimes \A_{x,y}
\ar[r]^-{\comp_{x,y,z}}
&
\A_{x,z}
\\
\A_{y,z} \otimes (\unc \otimes \A_{x,y})
\ar[rr]_{1 \otimes (\unit_y \otimes 1)}
&
&
\A_{y,z} \otimes (\A_{y,y} \otimes \A_{x,y})
\ar[u]_{1 \otimes \comp_{x,y,y}} 
}
$$
and\\ 
$\Xi_3$ $=$
$$
\xymatrix@C=3pc{
& 
\txt{
$(\A_{y,z} \otimes \unc)$\\ 
$\otimes \A_{x,y}$
}
\ar[r]^-{(1 \otimes \unit_y) \otimes 1}
\ar[dd]^{A'}
&
\txt{
$(\A_{y,z} \otimes \A_{y,y})$\\
$\otimes \A_{x,y}$
}
\ar[r]^-{\comp_{y,y,z} \otimes 1}
\ar[dd]_{A'}
&
\A_{y,z} \otimes \A_{x,y}
\ar[rd]^{\comp_{x,y,z}}
\\
\A_{y,z} \otimes \A_{x,y}
\ar[ru]^{R' \otimes 1}
\ar[rd]_{1 \otimes L'}
\ar@{}[r]|-{=}
&
\ar@{}[r]|-{=}
&
&
&
\A_{x,z}
\\
& 
\txt{
$\A_{y,z} \otimes$\\ 
$(\unc \otimes \A_{x,y})$
}
\ar[r]_-{1 \otimes (\unit_y \otimes 1)}
&
\txt{
$\A_{y,z} \otimes$\\ 
$(\A_{y,y} \otimes \A_{x,y})$
}
\ar[r]_-{1 \otimes \comp_{x,y,y}}
\ar@{=>}[ruu]|{\assp_{x,y,y,z}}
&
\A_{y,z} \otimes \A_{x,y}
\ar[ru]_{\comp_{x,y,z}}
}
$$
satisfy the equality
$\Xi_1 = \Xi_3 * \Xi_2$
\end{tag}

Let us justify the equivalence of the definitions
of $\BS$-categories. We define the following
bijective correspondence between data involved the definitions.
Arrows ${\A(x,-)}_{y,z} : \A_{y,z} \rightarrow [\A_{x,y},\A_{x,z}]$
and $\comp_{x,y,z} : \A_{y,z} \otimes \A_{x,y} \rightarrow \A_{x,z}$
correspond via the adjunction \ref{punitens}.
By Lemma \cite{Sch08}-$\calcdeux$ one has a bijective correspondence
between 2-cells of the kind $\assp_{x,y,z,t}$
and 2-cells $\assh_{x,y,z,t}$, the later being images
by $\Res \circ \Res$ of the first ones.
The codomains of the 2-cells $\rp_{x,y}$ and $\rh_{x,y}$
are equal by Lemma \ref{calc51},
and these 2-cells correspond when are equal. 
The codomains of the 2-cells $\lp_{x,y}$ and $\lh_{x,y}$
are equal by Lemma \ref{calc6} and these 2-cells correspond 
when they are equal. 
For such corresponding data, the proofs of the equivalence 
of Axioms \ref{cohassp} and \ref{cohassh}
rely on the adjunction \ref{punitens}.
The 2-cells of Axiom \ref{cohassp} 
have images by $\Res \circ \Res \circ \Res$
the two 2-cells of Axioms \ref{cohassh}
and their common domain is a strict arrow 
with strict images by $\Res$ and $\Res \circ \Res$.
Computation details are in Appendix in \ref{pteqax1}.
The proof of the equivalence of Axioms \ref{cohunassp}
and Axioms \ref{cohunassh} is similar. The 2-cells 
$\Xi_1$ of Axiom \ref{cohunassp} has a strict domain
and its image by $\Res$ is $\rh$ whereas the 2-cell
$\Xi_3 \circ \Xi_2$ has image by $\Res$ the second 
2-cell of Axiom \ref{cohunassh}. Computation details
are in Appendix in \ref{pteqax2}.\\

We have an alternative definition for the $\BS$-functors
with the tensor in $\BS$.\\

Given two arbitrary $\BS$-categories $\A$ and $\B$, 
a {\em $\BS$-functor} $F: \A \rightarrow \B$ 
consists of the following data:\\
- A map $F$ sending objects of $\A$ to objects
of $\B$;\\
- For any objects $x$,$y$ of $\A$,
and arrow $F_{x,y}: \A(x,y) \rightarrow \B(Fx,Fy)$ in $\BS$;\\
- Collections of 2-cells of $\BS$:
the $F^2_{x,y}$, indexed by pair 
of objects $x$,$y$ of $\A$ and the $F^0_x$, indexed by objects
$x$ of $\A$, as follows
\begin{tag}\label{F2}
{\small
$$
\xymatrix@C=3pc{
\A_{y,z} \otimes \A_{x,y} 
\ar[r]^-{F_{y,z} \otimes F_{x,y}}
\ar[d]_{\comp}
&
\B_{Fy,Fz} \otimes \B_{Fx,Fy}
\ar[d]^{\comp}
\ar@{=>}[ld]^{F^2_{x,y,z}}
\\
\A_{x,z} 
\ar[r]_{F_{x,z}}
&
\B_{Fx,Fz}
}
$$
}
\end{tag}
and
{\small 
$$
\xymatrix{
\unc 
\ar[rr]^{\unit_x} 
\ar[rrdd]_{\unit_{Fx}}
&  
& 
\A_{x,x}
\ar[dd]^{F_{x,x}}
\\
&
\ar@{=>}[ru]^{F^0_x}
&
\\
&
&
\B_{Fx,Fx}
}
$$
}
and that satisfy the coherence conditions \ref{cohfun1p}, \ref{cohfun2p} 
and \ref{cohfun3p} below.\\
\begin{tag}\label{cohfun1p}
For any objects $x$,$y$,$z$,$t$ of $\A$, the 2-cells 
{\tiny
$$\xymatrix@C=4pc{
\txt{
$(\A_{z,t} \otimes \A_{y,z})$\\ 
$\otimes \A_{x,y}$
}
\ar[r]^{\comp \otimes 1} 
\ar[d]_{(F_{z,t} \otimes F_{y,z}) \otimes F_{x,y}}
&
\A_{y,t} \otimes \A_{x,y}
\ar[r]^{\comp}
\ar[d]|{F_{y,t} \otimes F_{x,y}}
& 
\A_{x,t}
\ar[d]^{F_{x,t}}
\\
\txt{
$(\B_{Fz,Ft} \otimes \B_{Fy,Fz})$\\
$\otimes \B_{Fx,Fy}$
}
\ar@{=>}[ru]^{F^2_{y,z,t} \otimes 1}
\ar[r]_{\comp \otimes 1}
\ar[d]_{A'} 
&
\B_{Fy,Ft} \otimes \B_{Fx,Fy}
\ar@{=>}[ru]^{F^2_{x,y,t}}
\ar[r]^{\comp}
&
\B_{Fx,Ft}
\ar@{-}[d]^{id}
\\
\txt{
$\B_{Fz,Ft} \otimes$\\
$(\B_{Fy,Fz} \otimes \B_{Fx,Fy})$
}
\ar[r]_{1 \otimes \comp}
\ar@{=>}[rru]_{\assp_{Fx,Fy,Fz,Ft}}
&
\B_{Fz,Ft} \otimes \B_{Fx,Fz}
\ar[r]_-{\comp}
&
\B_{Fx,Ft}
}
$$
}
and
{\tiny 
$$
\xymatrix@C=4pc{
\txt{
$(\A_{z,t} \otimes \A_{y,z})$\\ 
$\otimes \A_{x,y}$
}
\ar[r]^{\comp \otimes 1}
\ar[d]_{A'}
&
\A_{y,t} \otimes \A_{x,y}
\ar[r]^{\comp}
&
\A
\ar@{-}[d]^{id}
\\
\txt{
$\A_{z,t} \otimes$\\
$(\A_{y,z} \otimes \A_{x,y})$
}
\ar[d]_{F_{z,t} \otimes (F_{y,z} \otimes F_{x,y})}
\ar[r]_{1 \otimes \comp}
\ar@{=>}[rru]^{\assp_{x,y,z,t}}
&
\A_{z,t} \otimes \A_{x,z}
\ar[r]^{\comp}
\ar[d]|{F_{z,t} \otimes F_{x,z}}
&
\A_{x,t}
\ar[d]^{F_{x,t}}
\\
\txt{
$\B_{Fz,Ft} \otimes$\\
$(\B_{Fy,Fz} \otimes \B_{Fx,Fy})$
}
\ar[r]_{1 \otimes \comp}
\ar@{=>}[ru]_{1 \otimes F^2_{x,y,z}}
&
\B_{Fz,Ft} \otimes \B_{Fx,Fz}
\ar[r]_{\comp}
\ar@{=>}[ru]_{F^2_{x,z,t}}
&
\B_{Fx,Ft}
}
$$
}
are equal.
\end{tag}
\begin{tag}\label{cohfun2p}
For any objects $x$,$y$ of $\A$, 
the 2-cells
{\tiny
$$
\xymatrix@C=4pc{
\A_{x,y} 
\ar[r]^{F_{x,y}}
\ar[d]_{R'}
\ar@{}[rd]|{=}
&
\B_{Fx,Fy}
\ar[d]^{R'}
\ar@/^70pt/@{-}[ddd]^{id}
\\
\A_{x,y} \otimes \unc
\ar[d]_{1 \otimes \unit}
\ar[r]^{F_{x,y} \otimes 1}
&
\B_{Fx,Fy} \otimes \unc
\ar[d]^{1 \otimes \unit}
\ar@{=>}[ld]_{1 \otimes F^0_x}
&
\\
\A_{x,y} \otimes \A_{x,x}
\ar[d]_{\comp}
\ar[r]^{F_{x,y} \otimes F_{x,x}}
&
\B_{Fx,Fy} \otimes \B_{Fx,Fx}
\ar@{=>}[ld]_{F^2_{x,x,y}}
\ar[d]^{\comp}
\ar@{}[ru]|{\Leftarrow^{\rp_{Fx,Fy}}}
\\
\A_{x,y} 
\ar[r]_{F_{x,y}}
&
\B_{Fx,Fy}
}
$$
}
and
{\tiny
$$
\xymatrix{
\A_{x,y} 
\ar@{-}[rr]^{id}
\ar[d]_{R'}
&
\ar@{=>}[d]^{\rp_{x,y}}
&
\A_{x,y}
\ar[r]^F 
&
\B_{Fx,Fy}
\\
\A_{x,y} \otimes \unc
\ar[rr]_{1 \otimes \unit_x}
&
&
\A_{x,y} \otimes \A_{x,x}
\ar[u]_{\comp}
}
$$
}
are equal.
\end{tag}
\begin{tag}\label{cohfun3p}
For any objects $x$,$y$ of $\A$, the 2-cells 
{\tiny
$$\xymatrix@C=4pc{
\A_{x,y} 
\ar[r]^{F}
\ar[d]_{L'}
\ar@{}[rd]|{=}
&
\B_{Fx,Fy}
\ar@{-}@/^75pt/[ddd]^{id}
\ar[d]^{L'}
\\
\unc  \otimes \A_{x,y}
\ar[d]_{\unit_y \otimes 1}
\ar[r]^{1 \otimes F_{x,y}}
&
\unc \otimes \B_{Fx,Fy}
\ar[d]|{\unit_{Fy} \otimes 1}
\ar@{=>}[ld]_{ F^0_{x,y} \otimes 1}
&
\\
\A_{y,y} \otimes \A_{x,y}
\ar[d]_{\comp}
\ar[r]^{F_{y,y} \otimes F_{x,y}}
&
\B_{Fy,Fy} \otimes \B_{Fx,Fy}
\ar@{=>}[ld]_{F^2_{x,y,y}}
\ar[d]^{\comp}
\ar@{}[ru]_{\Leftarrow^{\lp_{Fx,Fy}}}
\\
\A 
\ar[r]_{F_{x,y}}
&
\B_{Fx,Fy}
}
$$
}
and
{\tiny 
$$
\xymatrix{
\A_{x,y} 
\ar@{-}[rr]^{id}
\ar[d]_{L'}
&
\ar@{=>}[d]^{\lp}
&
\A_{x,y}
\ar[r]^{F_{x,y}} 
&
\B_{Fx,Fy}
\\
\unc \otimes \A_{x,y}
\ar[rr]_{\unit_y \otimes 1}
&
&
\A_{y,y} \otimes \A_{x,y}
\ar[u]_{\comp}
}
$$
}
are equal.
\end{tag}

The adjunction \ref{punitens} 
gives a bijective correspondence between
2-cells $F^2$ as in \ref{F2} and
2-cells ${F'}^2$ as in \ref{Fp2}.
For such coresponding data, it turns out that
Axioms \ref{cohfun1h} and \ref{cohfun1p} are equivalent, 
this is proved in \ref{eqcofun1hp} in Appendix,
Axioms \ref{cohfun2h} and \ref{cohfun2p} are equivalent,
this is proved in \ref{eqcofun2hp} in Appendix, and
Axioms \ref{cohfun3h} and \ref{cohfun3p} are equivalent,
this is proved in \ref{eqcofun3hp} in Appendix.\\

We say that a $\BS$-category $\A$ as above,
is {\em strict} if and only if 
the arrows $\A(x,-)$ are strict in $\BS$
and the 2-cells $\assp$, $\rp$, $\lp$,
or equivalently the 2-cells
$\assh$, $\rh$ and $\lh$, are all identities.
We have also the notion of {\em strict} $\BS$-functors
$(F,F^0,F^2): \A \rightarrow \B$: they are the ones 
for which the components $F_{x,y}: \A_{x,y} \rightarrow
\B_{Fx,Fy}$ are strict arrows in $\BS$ and for which 
the 2-cells of the collections $F^0$ and ${F'}^2$
are identities. Note that for an $F$ as above with
$\A$ strict the ${F'}^2$ are identities if and only if 
the $F^2$ are.
\end{section}

\begin{section}{First examples}\label{Exples}
One-point enrichments are of particular interest and 
named {\em 2-rings.} As mentioned by M.Dupont in his thesis 
\cite{Dup08}, they are also the {\em categorical rings} defined by 
Jibladze and Pirashvili \cite{JiPi07}, those are also known to 
be the {\em Ann-categories} of \cite{Qu87}. Given a 2-ring $\A$, we shall
write simply $\A$ for the hom $\A(*,*)$ of its unique object
$*$. Therefore the formal definition of a 2-ring 
$(\A,\comp,\unit, \assp,\rp,\lp)$,
as a ``weak'' monoid in $\BS$, namely 
a Picard category $\A$ with a {\em multiplication} 
$\comp: \A \otimes \A \rightarrow \A$
and a {\em unit}
$\unit: \unc \rightarrow \A$ (which are strict arrows!) and
appropriate 2-cells $\assp$, $\rp$ and $\lp$, is obtained
by removing the subscripts $x,y,z,...$ from the definitions
of $\BS$-categories.
We shall use the alternative definition of a 2-ring
$(\A,\comh,\unit, \assh,\rh,\lh)$ not using the tensor
obtained similarly by forgetting subscripts and where multiplication 
$\comh: \A \rightarrow [\A,\A]$ denotes the unique arrow $\A_{*,-}$.\\

The following definition of 2-rings can be obtained 
by written explicitly all linearity conditions from
the definition of enriched categories and functors. 
It is equivalent and very close to that of Jibladze and Pirashvili
\cite{JiPi07} (for their categorical rings).
Detailed explanations that the 2-rings with their morphisms
in the sense below are just one-point $\BS$-categories with 
their functors is given in Appendix-\ref{pfeqdef2rgobj}
and \ref{pfeqdef2rgmor}.

\begin{definition}\label{defobjJP}
A {\em 2-ring}  consists of a symmetric Picard category
$(\A,\jj)$, where $\A$ is denoted additively
$(\A,+,0,\ac,\rc,\lc,\syc)$, 
together with a functor $\A \times \A \rightarrow \A$,
denoted by a multiplication ``.'', an object $\unm$ of 
$\A$ and natural
isomorphisms
$$\assmu_{a,b,c}: (a.b).c \rightarrow a.(b.c),$$
$$\rmu_a: a.\unm \rightarrow a,$$
$$\lmu_a: \unm.a \rightarrow a,$$
$$\underline{a}_{b,b'}: a.b + a.b' \rightarrow a.(b+b'),$$
$$\overline{b}_{a,a'}: a.b + a'.b \rightarrow (a+a').b,$$
such that the data $(\A, ., \one, \assmu, \rmu, \lmu)$ 
defines a monoidal structure on $\A$ and 
the diagrams below from \ref{amofun1} to \ref{bJP5}
commute for all possible objects of $\A$.\\

\begin{tag}\label{amofun1}
{\small
$$\xymatrix{
a.b + (a.b' + a.b'')
\ar[d]_{id + \underline{a}_{b',b''}}
\ar[r]^-{\ac}
&
(a.b + a.b') + a.b''
\ar[d]^{\underline{a}_{b,b'} + id}
\\
a.b + a.(b'+b'')
\ar[d]_{\underline{a}_{b,b'+b''} }
&
a.(b+b') + a.d
\ar[d]^{\underline{a}_{b+b',b''}}
\\
a.(b.(b'+b''))
\ar[r]_-{a.\ac}
&
a.((b+b')+b'')
}$$
}
\end{tag}

\begin{tag}\label{bmofun1}
{\small
$$\xymatrix{
a.b + (a'.b + a''.b)
\ar[d]_{id + \overline{b}_{a',a''}}
\ar[r]^-{\ac}
&
(a.b + a'.b) + a''.b
\ar[d]^{\overline{b}_{a,a'} + id}
\\
a.b + (a'+a'').b
\ar[d]_{\overline{b}_{a,a'+a''} }
&
(a+a').b + a''.b
\ar[d]^{\overline{b}_{a+a',a''}}
\\
(a+(a'+a'')).b
\ar[r]_-{\ac.b}
&
((a+a')+a'').b
}$$
}
\end{tag}

\begin{tag}\label{asymofun}
{\small
$$
\xymatrix{
a.b + a.b' 
\ar[r]^-{\underline{a}_{b,b'}}
\ar[d]_{\syc}
&
a.(b+b')
\ar[d]^{a.\syc}
\\
a.b' + a.b
\ar[r]_-{\underline{a}_{b',b}}
&
a.(b'+b).
}
$$
}
\end{tag}
\begin{tag}\label{bsymofun}
{\small
$$
\xymatrix{
a.b + a'.b 
\ar[r]^-{\overline{b}_{a,a'}}
\ar[d]_{\syc}
&
(a+a').b
\ar[d]^{\syc.b}
\\
a'.b + a.b
\ar[r]_-{\overline{b}_{a',a}}
&
(a'+a).b
}
$$
}
\end{tag}

\begin{tag}\label{lastaxiom}
{\small
$$\xymatrix{
(a.b + a.b')
+
(a'.b + a'.b')
\ar@{-}[rr]^-{\cong}
\ar[d]_{\underline{a}_{b,b'} + \underline{a'}_{b,b'}}
&
&
(a.b + a'.b)
+
(a'.b + a'.b')
\ar[d]^{\overline{b}_{a,a'} + \overline{b'}_{a,a'}}
\\
a.(b+b')
+
a'.(b+b')
\ar[rd]_{\overline{b+b'}_{a,a'}}
&
&
(a+a').b
+
(a+a').b'
\ar[ld]^{\underline{a+a'}_{b,b'}}
\\
&
(a+a').(b+b')
}$$
}
\end{tag}

\begin{tag}\label{bJP1}
{\small
$$
\xymatrix{
(a.b).c + (a'.b).c
\ar[r]^-{\overline{c}_{a.b,a'.b}}
\ar[d]_{\assmu_{a,b,c} + \assmu_{a',b,c}}
&
((a.b)+(a'.b)).c
\ar[r]^-{\overline{b}_{a,a'}.c}
&
((a+a').b).c
\ar[d]^{\assmu_{a+a',b,c}}\\
a.(b.c) + a'.(b.c)
\ar[rr]_-{\overline{b.c}_{a,a'}}
& 
&
(a+a').(b.c)
}
$$
}
\end{tag}
\begin{tag}\label{bJP2}
{\small
$$
\xymatrix{
(a.b).c + (a.b').c
\ar[r]^-{\overline{c}_{a.b,a.b'}}
\ar[d]_{\assmu_{a,b,c}+\assmu_{a,b',c}}
&
((a.b)+(a.b')).c
\ar[r]^-{\underline{a}_{b,b'}.c}
&
(a.(b+b')).c
\ar[d]^{\assmu_{a,b+b',c}}
\\
a.(b.c) + a.(b.c')
\ar[r]_-{\underline{a}_{b.c,b.c'}}
&
a.(b.c+b.c')
\ar[r]_-{a.\overline{c}_{b,b'}}
&
a.((b+b').c)
}
$$
}
\end{tag}
\begin{tag}
\label{bJP3}
{\small
$$
\xymatrix{
(a.b).c + (a.b).c' 
\ar[rr]^{\underline{a.b}_{c,c'}}
\ar[d]|{\assmu_{a,b,c} + \assmu_{a,b,c'}}
& &
(a.b).(c+c')
\ar[d]|{\assmu_{a,b,c+c'}}
\\
a.(b.c) + a.(b.c')
\ar[r]_{\underline{a}_{b.c,b.c'}}
&
a.(b.c+b.c')
\ar[r]_{a.\underline{b}_{c,c'}}
&
a.(b.(c+c'))
}
$$
}
\end{tag}
\begin{tag}\label{bJP4}
{\small
$$
\xymatrix{
(a+b).\unm
\ar[d]_{ \overline{\unm}_{a,b}}
\ar[r]^-{ \rmu_{a+b}}
& a+b  \\
a+b
\ar@{=}[ru]. 
}
$$
}
\end{tag}
\begin{tag}\label{bJP5}
{\small
$$
\xymatrix{
\unm.(a+b) 
\ar[d]_{\underline{\unm}_{a,b}}
\ar[r]^-{\lmu_{a+b}}
& a+b \\
a+b
\ar@{=}[ru] 
}
$$
}
\end{tag}
\end{definition}

Note that in the definition given in \cite{JiPi07} inverses
of maps $\underline{a}_{b,c}$ and $\overline{a}_{b,c}$
rather than the maps themselves are considered and
diagrams \ref{amofun1} and \ref{asymofun} are replaced by
\begin{tag}
{\small
$$
\xymatrix{
& 
a.(b+b')+a.(c+c')
\ar[rd]^{ \underline{a}_{b,b'}+ \underline{a}_{c,c'}}
& 
\\
a.((b+b')+(c+c'))
\ar[ru]^{ \underline{a}_{b+b',c+c'} }
\ar@{-}[d]_{a.\cong}
&
&
(a.b + a.b') + (a.c + a.c')
\ar@{-}[d]^{\cong}
\\
a.( (b+ c) + (b'+c')) 
\ar[rd]_{ \underline{a}_{b+c,b'+c'} }
&
&
(a.b +a.c)+(a.b'+a.c')
\\
&
a.(b+c)+a.(b'+c')
\ar[ru]_{ \underline{a}_{b,c} + \underline{a}_{b',c'} }
&
}
$$
}
\end{tag}
and similarly diagrams \ref{bmofun1} and \ref{bsymofun} are replaced by 
\begin{tag}
{\small
$$
\xymatrix{
& 
(a+a').c+ (b+b').c
\ar[rd]^{\overline{c}_{a,a'} + \overline{c}_{b,b'}}
& 
\\
((a+a')+(b+b')).c
\ar[ru]^{ \overline{c}_{a+a',b+b'} }
\ar@{-}[d]_{\cong.c}
&
&
(a.c + a'.c) + (b.c + b'.c)
\ar@{-}[d]^{\cong}
\\
( (a + b) + (a'+ b')).c 
\ar[rd]_{\overline{c}_{a+b,a'+b'}}
&
&
(a.c + b.c)+(a'.c + b'.c)
\\
&
(a + b).c + (a' + b').c'
\ar[ru]_{\overline{c}_{a,b} + \overline{c}_{a',b'}}
&
}
$$
}
\end{tag}
The definitions here and in \cite{JiPi07} are indeed equivalent 
(To see this use for instance Lemma \cite{Sch08}-\lemsymofundeux.)

\begin{definition}\label{defmorJP}
A morphism of categorical ring $\A \rightarrow \B$
consists of a functor $H:\A \rightarrow \B$
with a symmetric monoidal structure between the 
symmetric categorical groups 
$$\Hadd: (\A,+,0,\ac,\rc,\lc,\syc) \rightarrow 
(\B,+,0,\ac,\rc,\lc,\syc)$$
and a monoidal structure between the monoidal categories
$$\Hmu: (\A,.,1,\assmu,\rmu,\lmu) \rightarrow 
(\B,.,1,\assmu,\rmu,\lmu)$$\\
such that the following diagrams
\begin{tag}\label{diag1JP}
{\small
$$
\xymatrix{
& 
H(a).H(b)+H(a).H(b')
\ar[r]^-{ {\Hmu}^2_{a,b} + {\Hmu}^2_{a,b'} }
&
H(a.b)+H(a.b')
\ar[rd]^{{\Hadd}^2_{a.b,a.b'}}
\\
H(a).(H(b)+H(b'))
\ar[ru]^{ \underline{H(a)}_{ H(b),H(b') } }
\ar[rd]_{ H(a).{\Hadd^2}_{b,b'} }
& 
& 
&
H(a.b+a.b')
\\
&
H(a).H(b+b')
\ar[r]_-{ {\Hmu}^2_{a,b+b'} }
&
H(a.(b+b'))
\ar[ru]_{ H ( \underline{a}_{b,b'} ) }
}
$$
}
\end{tag}
and
\begin{tag}\label{diag2JP}
{\small
$$
\xymatrix{
& 
H(a).H(b)+H(a').H(b)
\ar[r]^-{{\Hmu}^2_{a,b} + {\Hmu}^2_{a',b}}
&
H(a.b)+H(a'.b)
\ar[rd]^{{\Hadd}^2_{a.b,a.b'} }\\
(H(a)+H(a')).H(b)
\ar[rd]_{{\Hadd}^2_{a,a'}.H(b)}
\ar[ru]^{\overline{H(b)}_{H(a),H(a')}}
& & & 
H(a.b+a'.b)\\
&
H(a+a').H(b)
\ar[r]_-{{\Hmu}^2_{a+a',b}}  
& 
H((a+a').b)
\ar[ru]_{H(\overline{b}_{a,a'})}
}
$$
}
\end{tag}
commute for all possible objects involved.  
\end{definition}

The following is a crucial example of $\BS$-category.
\begin{proposition}\label{PicsEn}
The 2-category $\BS$ gets strictly enriched over itself, i.e.
it admits a strict enriched structure as follows. The hom map 
sends any pair $\A$,$\B$ of objects to 
$[\A,\B]$, the composition maps are the  
$[\A,-]_{\B,\C}: [\B,\C] \rightarrow [[\A,\B],[\A,\C]]$ and
the unit arrows $\unit_{\A}$ are the $\vv: \unc \rightarrow [\A,\A]$.
\end{proposition}
\pf
See \ref{pfPicsEn} in Appendix.
\epf

Let us make the following remark about the terminology.
If $\BS'$ denotes just for the purpose 
of this explanation the enriched structure
of $\BS$ over itself then for any $\A$ in $\BS$
and any 1-cell $F:\B \rightarrow \C$ in $\BS$, 
one has that:\\
- ${\BS'(\A,-)}_{\B,\C}$ is $[\A,-]_{\B,\C}: [\B,\C] \rightarrow [[\A,\B],[\A,\C]]$;\\
- ${\BS'(-,\C)}_{\A,B}$ is 
${[-,\C]}_{\A,\B}: [\A,\B] \rightarrow [[\B,\C],[\A,\C]]$;\\
- $\BS'(1,F)$ is $[1,F]: [\A,\B] \rightarrow [\A,\C]$;\\
- $\BS'(F,1)$ is $[F,1]: [\B,\C] \rightarrow [\A,\C]$.\\ 
The first two points results from the definitions.
The other two points are the following lemma proved
in Appendix \ref{pflemtoto}.
\begin{lemma}\label{lemtoto}
For any $\A$ and any arrows $F: \B \rightarrow \C$ and 
$\tilde{F}: \unc \rightarrow [\B,\C]$ strict 
with $\ev_{\gen}(\tilde{F}) = F$ the diagrams
in $\BS$
$$
\xymatrix{
[\C,\A]
\ar[r]^-{[F,\A]}
\ar[d]_{R'}
&
[\B,\A]
\\
[\C,\A] \otimes \unc
\ar[r]_{1 \otimes \tilde{F}}
&
[\C,\A] \otimes [\B,\C]
\ar[u]_{\comp}
}
$$
and 
$$
\xymatrix{
[\A,\B]
\ar[r]^-{[\A,F]}
\ar[d]_{L'}
&
[\A,\C]
\\
\unc \otimes [\A,\B] 
\ar[r]_{\tilde{F} \otimes 1}
&
[\B,\C] \otimes [\A,\B]
\ar[u]_{\comp}
}
$$
both commute. 
\end{lemma}

Given any $\BS$-category $\A$, 2-rings are obtained
by restriction of $\A$ to anyone of its points.
The particular case of the strict enriched structure 
on $\BS$ yields a strict 2-ring structure on $[\A,\A]$
for any Picard category $\A$.\\

Another important example of 2-ring is provided by 
the unit $\unc$ of $\BS$. 
\begin{proposition}\label{unit2rg}
The unit $\unc$ of $\BS$ admits a strict 2-ring structure
with multiplication given by 
$L_{\unc}: \unc \otimes \unc \rightarrow \unc$
(or $\vv: \unc \rightarrow [\unc,\unc]$) 
and unit the identity at $\unc$. 
\end{proposition}
\pf 
See Appendix \ref{pfunit2rg}.
\epf

\end{section}

\begin{section}{Modules and their morphisms}\label{Mod}
Any 2-ring $\A$ yields a category $\A$-mod of $\A$-modules
and their morphisms.
Formally $\A$-mod is the category of $\BS$-functors $\A \rightarrow \BS$ 
and  $\BS$-natural transformations between them.
In this section we present alternative descriptions of $\A$-modules
and their morphisms. In particular we show that the category $\A$-mod 
is isomorphic to a category of $T$-algebras and their morphisms 
for the doctrine $T = \A \otimes -$ over $\BS$.
Eventually we prove in Proposition \ref{ImodeqBS} that the category 
$\unc-mod$ of modules over the unit 2-ring $\unc$ is equivalent to $\BS$.\\

In this section $\A$ stand for a 2-ring with multiplication
$\comh:\A \rightarrow [\A,\A]$/ $\comp: \A \otimes \A \rightarrow \A$
with unit $\unit: \unc \rightarrow \A$ and 
coherence 2-cell $\assh$/$\assp$ (\ref{assh}/\ref{assp}), 
$\rp$/$\rh$ (\ref{rp}/\ref{rh}) and $\lp$/$\lh$ (\ref{lp}/\ref{lh}).\\
 
Considering an arbitrary $\BS$-functor $F: \A \rightarrow \BS$,
let us write $\M$ for the object $F(\star)$ of $\BS$
image by
$F$ of the unique point $\star$ of $\A$,
$\multh$ for the arrow unique component 
$F_{\star,\star}: \A \rightarrow [\M,\M]$
of $F$,
$\assmh$ for the 2-cell ${F'}^2_{\star,\star,\star}$  
and $\lmz$ for the 2-cell $F^0_{\star}$.
Then one obtains the following first definition 
of $\A$-modules by rewriting the data \ref{Fp2}
and \ref{F0} and Axioms \ref{cohfun1h}, \ref{cohfun2h}
and \ref{cohfun3h} with these new notations.\\
 
A $\A$-module $\M = (\M,\multh,\assmh,\lmz)$ consists 
of the following data in $\BS$: an object $\M$,
with an arrow $\multh : \A \rightarrow [\M,\M]$,
called its {\em action}, 
and two 2-cells $\assmh$ and $\lmz$ as follows
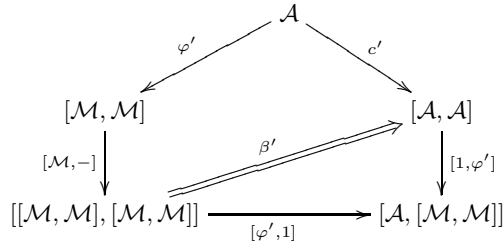
\begin{tag}\label{assmh}
{\small
$$
\xymatrix{
&
\A
\ar[rd]^{\comh}
\ar[ld]_{\multh}
&
\\
[\M,\M]
\ar[d]_{[\M,-]}
&
&
[\A,\A]
\ar[d]^{[1,\multh]}
\\
[[\M,\M],[\M,\M]]
\ar[rr]_{[\multh,1]}
\ar@{=>}[rru]^{\assmh}
&
&
[\A,[\M,\M]]
}
$$
}
\end{tag}
\begin{tag}\label{lmz}
{\small 
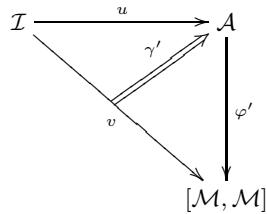
$$
\xymatrix{
\unc 
\ar[rr]^{\unit} 
\ar[rrdd]_{\vv}
&  
& 
\A
\ar[dd]^{\multh}
\\
&
\ar@{=>}[ru]^{\lmz}
&
\\
&
&
[\M,\M]
}
$$
}
\end{tag}
and those satisfy the coherence axioms
\ref{cohmdobj1h}, \ref{cohmdobj2h} and \ref{cohmdobj3h} below.
\begin{tag}\label{cohmdobj1h}
The 2-cells 
{\tiny
$$\xymatrix@C=3pc{
& & 
\txt{
$[[\M,\M],$\\
$[\M,\M]]$
}
\ar[r]^-{[[\M,\M],-]}
\ar@<4ex>@{=>}[d]^{\assh_{\M,\M,\M,\M}}
&
\txt{
$[[[\M,\M],$\\
$[\M,\M]],$\\
$[[\M,\M],$\\
$[\M,\M]]]$
}
\ar[d]^{[[\M,-],1]}
\\
\A
\ar[d]_{\comh}
\ar[r]^-{\multh}
&
[\M,\M]
\ar[ru]^{[\M,-]}
\ar[r]^-{[\M,-]}
&
\txt{
$[[\M,\M],$\\
$[\M,\M]]$
}
\ar@{=>}[lld]_{\assmh}
\ar[r]^-{[1,[\M,-]]}
\ar[d]^{[\multh,1]}
\ar@{}[rd]|{=}
&
\txt{
$[[\M,\M],$\\
$[[\M,\M],$\\
$[\M,\M]]]$
}
\ar[d]^{[\multh,1]}
\\
[\A,\A]
\ar[rr]^-{[1,\multh]}
\ar[rrd]_{[1,\comh]}
& & 
\txt{
$[\A,$\\
$[\M,\M]]$
}
\ar[r]^-{[1,[\M,-]]}
\ar@{=>}[d]^{[1,\assmh]}
&
\txt{
$[\A,$\\
$[[\M,\M],$\\
$[\M,\M]]]$
}
\ar[r]^-{[1,[\multh,1]]}
&
\txt{
$[\A,$\\
$[\A,$\\
$[\M,\M]]]$
}
\\
& & 
[\A,[\A,\A]]
\ar[rru]_{[1,[1,\multh]]}
}$$
}
and
{\tiny 
$$
\xymatrix@C=3pc{
[\M,\M]
\ar[r]^-{[\M,-]}
&
\txt{
$[[\M,\M],$\\
$[\M,\M]]$
}
\ar[d]^{[\multh,1]}
\ar[r]^-{[\A,-]}
\ar@{=>}[ddl]_{\assmh}
\ar@{}[rd]|{=}
&
\txt{
$[[\A,[\M,\M]],$\\
$[\A,[\M,\M]]]$
}
\ar[r]^-{[[\multh,1],1]}
\ar[d]^{[[1,\multh],1]}
&
\txt{
$[[[\M,\M],$\\
$[\M,\M]],$\\
$[\A,[\M,\M]]]$
}
\ar[r]^-{[[\M,-],1]}
\ar@{=>}[d]^{[\assmh,1]}
&
\txt{
$[[\M,\M],$\\
$[\A,[\M,\M]]]$
}
\ar[d]^{[\multh,1]}
\\
& 
[\A,[\M,\M]]
\ar[r]^-{[\A,-]}
\ar@{}[rd]|{=}
&
\txt{
$[[\A,\A],$\\
$[\A,[\M,\M]]]$
}
\ar[rr]_-{[\comh,1]}
& 
& 
\txt{
$[\A,[\A,$\\
$[\M,\M]]]$}
\\
\A
\ar[r]^-{\comh}
\ar[uu]^{\multh}
\ar[drr]_{\comh}
&
[\A,\A]
\ar[r]^-{[\A,-]}
\ar[u]_{[1,\multh]}
&
\txt{
$[[\A,\A],$\\
$[\A,\A]]$
}
\ar[rr]^{[\comh,1]}
\ar[u]_{[1,[1,\multh]]}
\ar@{}[rru]|{=}
\ar@{=>}[d]^{\assh}
&
&
[\A,[\A,\A]]
\ar[u]_{[1,[1,\multh]]}
\\
&
&
[\A,\A]
\ar[rru]_{[1,\comh]}
}
$$
}
are equal.
\end{tag}
\begin{tag}\label{cohmdobj2h}
The 2-cells in $\BS$
{\tiny 
$$
\xymatrix{
\A
\ar[d]_{\comh} 
\ar[rr]^{\multh} 
& &  
[\M,\M]
\ar[d]|{[\M,-]}
\ar@{=>}[lld]_{\assmh}
\ar@/^70pt/[dddd]^1
\\
[\A,\A]
\ar[dd]_{[\unit,1]}
\ar[rd]^{[1,\multh]}
&
&
[[\M,\M],[\M,\M]]
\ar[dd]^{[\vv,1]}
\ar[ld]_{[\multh,1]}
\\
&
[\A,[\M,\M]]
\ar[rd]_{[\vv,1]}
&
\ar@{=>}[l]_-{[\lmz,1]}
& 
\ar@{}[l]_{=}
\\
[\unc,\A]
\ar@{}[rrd]|{=}
\ar@{}[ru]|{=}
\ar[d]_{\ev_{\gen}}
\ar[rr]^{[1,\multh]}
& & 
[\unc,[\M,\M]]
\ar[d]^{\ev_{\gen}}
\\
\A
\ar[rr]_{\multh}
& & 
[\M,\M]
}
$$
}
and
{\tiny
$$
\xymatrix{
\A
\ar@{-}[rr]^{id}
\ar[d]_{\comh}
&
\ar@{=>}[d]^{\rh}
&
\A
\ar[r]^{\multh}
& 
[\M,\M]
\\
[\A, \A]
\ar[rr]_{[\unit,1]}
& 
& 
[\unc,\A]
\ar[u]_{\ev_{\gen}}
}
$$
}
are equal.
\end{tag}
\begin{tag}\label{cohmdobj3h}
The 2-cells in $\BS$
{\tiny 
$$
\xymatrix{
\A
\ar[d]_{{\comh}^*} 
\ar[rr]^{\multh} 
& &  
[\M,\M]
\ar[d]|{[-,\M]}
\ar@{=>}[lld]_{{(\assmh)}^*}
\ar@/^70pt/[dddd]^1
\\
[\A,\A]
\ar[dd]_{[\unit,1]}
\ar[rd]^{[1,\multh]}
&
&
[[\M,\M],[\M,\M]]
\ar[dd]^{[\vv,1]}
\ar[ld]_{[\multh,1]}
\\
&
[\A,[\M,\M]]
\ar[rd]_{[\vv,1]}
&
\ar@{=>}[l]_-{[\lmz,1]}
& 
\ar@{}[l]_{=}
\\
[\unc,\A]
\ar@{}[rrd]|{=}
\ar@{}[ru]|{=}
\ar[d]_{\ev_{\gen}}
\ar[rr]^{[1,\multh]}
& & 
[\unc,[\M,\M]]
\ar[d]^{\ev_{\gen}}
\\
\A
\ar[rr]_{\multh}
& & 
[\M,\M]
}
$$
}
and
{\tiny 
$$
\xymatrix{
\A
\ar@{-}[rr]^{id}
\ar[d]_{\comh}
&
\ar@{=>}[d]^{\lh}
&
\A
\ar[r]^{\multh}
& 
[\M,\M]
\\
[\A, \A]
\ar[rr]_{[\unit,1]}
& 
& 
[\unc,\A]
\ar[u]_{\ev_{\gen}}
}
$$
}
are equal.
\end{tag}
We shall also denote by $\lmh$ for the 2-cell  
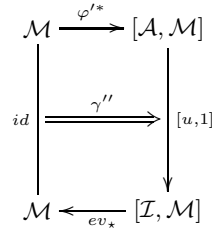
\begin{tag}\label{lmh}
{\small
$$
\xymatrix{
\M 
\ar[r]^-{\multh^*}
\ar@{-}[dd]_{id}
&
[\A,\M]
\ar[dd]^{[\unit,1]}
\\
\ar@{=>}[r]^{\lmh}
&
\\
\M
&
[\unc,\M]
\ar[l]^{\ev_{\gen}}
}
$$
}
\end{tag}
that corresponds to $\lmz$ via the bijection
\ref{2c2}/\ref{2c3}.\\

Consider now a $\BS$-natural transformation between
presheaves with domain a one point category $\A$
$$(\sigma,\kappa): F \rightarrow G: \A \rightarrow \BS.$$
Let us write 
$\M = (\M,\multh, \assmh, \lmh)$ and
$\N = (\N, \psi', \assmh,\lmh)$ for the two modules
corresponding respectively to $F$ and $G$.  
The $\BS$-natural transformation $\sigma$ has a unique 
component at the unique object $\gen$ of $\A$,
which is a strict arrow 
$\unc \rightarrow [\M,\N]$ in $\BS$ or equivalently an 
arrow $\mor: \M \rightarrow \N$.
The collection $\kappa$ consists of a unique 2-cell
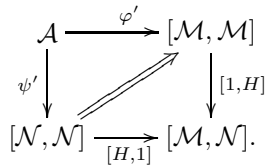
\begin{tag}\label{deltah}
$$
\xymatrix{
\A
\ar[r]^{\multh}
\ar[d]_{\psi'}
&
[\M,\M]
\ar[d]^{[1,\mor]}
\\
[\N,\N]
\ar[r]_{[\mor,1]}
\ar@{=>}[ru]
&
[\M,\N]
.}
$$
\end{tag}
which we name $\deltah$.
We obtain therefore the following definition
of morphism of $\A$-modules
by rewriting  Axioms \ref{cohnat1} and 
\ref{cohnat2} with these new notations.\\

A morphism of $\A$-module 
$(\mor,\deltah): (\M,\multh,
\assmh,\lmz) \rightarrow (\N,\psi',\assmh,\lmz)$
consists of an arrow $\mor:\M \rightarrow \N$ in $\BS$
with a 2-cell $\deltah$ as in \ref{deltah},
those satisfying Axioms \ref{cohmdmor1h}
and \ref{cohmdmor2h} below.
\begin{tag}\label{cohmdmor1h}
The 2-cells
$\Xi_1$, $\Xi_2$, $\Xi_3$, $\Xi_4$, $\Xi_5$, $\Xi_6$, $\Xi_7$ and $\Xi_8$  
in $\BS$ below satisfy the equality
$$\Xi_2 \circ \Xi_1 = 
\Xi_8 \circ \Xi_7 \circ \Xi_6 \circ \Xi_5 \circ \Xi_4 \circ \Xi_3.$$
$\Xi_1$ is
$$\xymatrix{
\A
\ar[d]_-{\comh}
\ar[r]^-{\psi'}
& 
[\N,\N]
\ar@{=>}[d]^{\assmh}
\ar[r]^-{[\N,-]}
&
[[\N,\N],[\N,\N]]
\ar[d]^{[\psi',1]}
\\
[\A,\A]
\ar[rr]_-{[1,\psi']}
& & 
[\A,[\N,\N]]
\ar[r]^-{[1,[\mor,1]]}
&
[\A,[\M,\N]]
}
$$
$\Xi_2$ is
$$
\xymatrix{
& & [\A,[\N,\N]]
\ar[rd]^{[1, [\mor,1]]}
\ar@{=>}[dd]^{[1,\deltah]}
\\
\A
\ar[r]^-{\comh}
&
[\A,\A]
\ar[ru]^{[1,\psi' ]}
\ar[rd]_{[1,\multh] }
& & 
[\A,[\M,\N]]
\\
& & 
[\A, [\M,\M]]
\ar[ru]_{[1, [1,\mor]]}
}
$$
$\Xi_3$ is the identity
$$
\xymatrix{
& & 
[[\N, \N],[\N,\N]]
\ar[rd]^{[1,[\mor,1]]}
\\
\A
\ar[r]^-{\psi'}
&
[\N,\N]
\ar[ru]^{[\N,-]}
\ar[rd]_{[\M,-]}
&
&
[[\N,\N], [\M,\N]]
\ar[r]^-{[\psi',1]}
&
[\A, [\M,\N]]
\\
& 
& 
[[\M,\N], [\M,\N]]
\ar[ru]_{[[\mor,1],1]}
\ar@{}[uu]|{=}
}
$$
$\Xi_4$ is
$$
\xymatrix{
& & & [[\N,\N],[\N,\N]]
\ar@{=>}[dd]^{[\deltah,1]}
\ar[rd]^{[\psi',1]}
\\
\A
\ar[r]^-{\psi'}
&
[\N,\N]
\ar[r]^-{[\M,-]}
&
[[\M,\N], [\M,\N]]
\ar[ru]^{[[\mor,1],1]}
\ar[rd]_{[[1,\mor],1]}
& & 
[\A, [\M,\N]]
\\
& & & 
[[\M,\M],[\M,\N]]
\ar[ru]_{[\multh,1]}
}$$
$\Xi_5$ is the identity 2-cell
$$
\xymatrix{
& & 
[[\M,\N], [\M,\N]]
\ar[rd]^{[[1,\mor],1]}
\ar@{}[dd]|{=}
\\
\A
\ar[r]^-{\psi'}
&
[\N,\N]
\ar[rd]_{[\mor,1]}
\ar[ru]^{[\M,-]}
&
&
[[\M,\M],[\M,\N]]
\ar[r]^-{[\multh,1]}
&
[\A, [\M,\N]] 
\\
& & 
[\M,\N]
\ar[ru]_{[\M,-]}
}
$$
$\Xi_6$ is
$$
\xymatrix{
& 
[\N,\N]
\ar[rd]^{[\mor,1]}
\ar@{=>}[dd]^{\deltah}
&
\\
\A
\ar[ru]^{\psi'}
\ar[rd]_{\multh}
& & 
[\M,\N]
\ar[r]^-{[\M,-]}
&
[[\M,\M],[\M,\N]]
\ar[r]^-{[\multh,1]}
&
[\A,[\M,\N]]
\\
& 
[\M,\M]
\ar[ru]_{[1,\mor]}
}
$$
$\Xi_7$ is the identity 2-cell
$$\xymatrix{
& & [\M,\N]
\ar[rd]^{[\M,-]}
\ar@{}[dd]|{=}
\\
\A
\ar[r]^-{\multh}
&
[\M,\M]
\ar[ru]^{[1,\mor]}
\ar[rd]_{[\M,-]}
& & 
[[\M,\M], [\M,\N]]
\ar[r]^-{[\multh,1]}
&
[\A, [\M,\N]]
\\
& & 
[[\M,\M], [\M,\M]]
\ar[ru]_{[1,[1,\mor]]}
}
$$
and $\Xi_8$ is
$$
\xymatrix{
\A
\ar[r]^-{\multh}
\ar[d]_{\comh}
&
[\M,\M]
\ar[r]^-{[\M,-]}
\ar@{=>}[d]^{\assmh}
&
[[\M,\M],[\M,\M]]
\ar[d]^{[\multh,1]}
\\
[\A,\A]
\ar[rr]_{[1,\multh]}
& & 
[\A, [\M,\M]]
\ar[r]^-{[1,[1,\mor]]}
&
[\A, [\M,\N]].
}
$$
\end{tag}
\begin{tag}\label{cohmdmor2h}
The 2-cells in $\BS$
$$
\xymatrix{
& 
& 
\ar@{=>}[dd]^{\lhpp_{\mor}}
& 
\\
& 
\ar@{=>}[d]^{\lmz}
& & 
\\
\unc
\ar@/^38pt/[rr]^-{\vv}
\ar[r]^-{\unit}
\ar@/^70pt/[rrr]^-{\mor}
&
\A
\ar[r]_-{\multh}
&
[\M,\M]
\ar[r]_-{[1,\mor]}
&
[\M,\N]
}
$$
and 
$$
\xymatrix{
& 
&
\ar@{}[d]|{=}
&
& 
\\
\unc 
\ar@{-}[d]_{id}
\ar[rr]^-{\vv}
\ar@/^38pt/[rrrr]^{\mor}
&
\ar@{=>}[d]^{\lmz}
& 
[\N,\N]
\ar[rr]^-{[\mor,1]}
\ar@{=>}[rrd]^{\deltah}
&
&
[\M,\N]
\\
\unc
\ar[rr]_-{\unit}
& 
&
\A
\ar[u]^{\psi'}
\ar[rr]_-{\multh}
&
&
[\M,\M] 
\ar[u]_{[1,\mor]}
}
$$
are equal.
\end{tag}

From these first definitions of $\A$-modules
and $\A$-module morphisms, 
one obtains immediately the following simple
other definitions involving multilinear maps and multilinear 
natural transformations.\\

A $\A$-module $(\M, \multl, \assml, \lml)$ consists
of a bilinear map $\multl: \A \times \M \rightarrow \M$,
which we write as a multiplication $\multl(a,m) = a . m$, 
with two natural transformations,  
$\assml$ which is trilinear, and $\lml$ which is linear, as follows.
\begin{tag}
$$\assml_{a_1,a_2,m}: a_1. (a_2 . m) \rightarrow (a_1. a_2). m$$
lies in $\M$, for $a_1$, $a_2$ objects of $\A$ and $m$ object of $\M$. 
\end{tag}
\begin{tag}
$$\lml_m: m \rightarrow \one_{\A} . m  $$ lies in $\M$
for $m$ object of $\M$.
\end{tag}
Those satisfy the coherence conditions 
\ref{cohmdobj1b},
\ref{cohmdobj2b} and 
\ref{cohmdobj3b} below.
\begin{tag}\label{cohmdobj1b}
For any objects $a_1$, $a_2$, $a_3$ of $\A$ and 
$m$ of $\M$ the diagram in $\M$ 
$$
\xymatrix{
a_1 . (a_2 . (a_3 . m)))
\ar[d]_{ a_1 . \assml_{a_2,a_3,m} }
\ar[r]^{ \assml_{a_1,a_2,a_3 . m} }
&
(a_1 . a_2) . (a_3 . a)
\ar[dd]^{ \assml_{a_1 . a_2,a_3,m} }
\\
a_1 . ((a_2 . a_3) . m)  
\ar[d]_{ \assml_{a_1, a_2 . a_3,m } }
&
\\
(a_1 . (a_2 . a_3)) . m
\ar[r]_{\assb_{a_1,a_2,a_3} . m}
& 
((a_1 . a_2) . a_3) . m
}
$$
commutes. 
\end{tag}
\begin{tag}\label{cohmdobj2b}
For any objects $a$ of $\A$ and $m$ of $\M$
the diagram in $\M$
$$
\xymatrix{
a .  m
\ar[rr]^{a . \lml_m}
\ar[ddr]_{\rb_a  . m}
& & 
a . (\one_{\A} . m)
\ar[ddl]^{\assml_{a,\one_{\A},m }}
\\
\\
& 
(a . \one_{\A}) . m
& 
}
$$
commutes.
\end{tag}
\begin{tag}\label{cohmdobj3b}
For any objects $a$ of $\A$ and $m$ of $\M$
the diagram in $\M$
$$
\xymatrix{
a . m
\ar[ddr]_{\lb_a . m}
\ar[rr]^{\lml_{a . m}}
& & 
\one_{\A} . (a . m)
\ar[ddl]^{\assml_{\one_{\A},a,m}}
\\
\\
& 
(\one_{\A} . a) . m
}
$$
commutes.
\end{tag}

A $\A$-module morphism $(\mor,\deltal) : \M \rightarrow \N$ consists
of an arrow $\mor :\M \rightarrow \N$ in $\BS$ with a 
bilinear natural transformation
\begin{tag}\label{deltal} 
$$\deltal_{a,m}: a. \mor(m) \rightarrow \mor(a.m)$$
which lies in $\N$ for $a$ object of $\A$ and 
$m$ object of $\M$
\end{tag} 
which satisfy the Axioms
\ref{cohmdmor1b}
and \ref{cohmdmor2b} below.
\begin{tag}\label{cohmdmor1b}
For any object $m$ of $\M$ the following diagram
in $\N$
$$
\xymatrix{
\mor m 
\ar[rr]^{\mor (\lml_m)}
\ar[rdd]_{\lml_{ \mor (m)}}
& & 
\mor (\gen .  m)
\\
\\
& 
\gen . \mor m
\ar[ruu]_{\deltal_{\gen,m}}
& 
}
$$
commutes.
\end{tag}
\begin{tag}\label{cohmdmor2b}
For any objects $a_1$, $a_2$ in $\A$ and $m$ in $\M$,
the following diagram in $\N$ 
$$
\xymatrix{
a_1 . (a_2 . \mor m)
\ar[rr]^-{\assml_{a_1,a_2,m}}
\ar[d]_{a_1 . \deltal_{a_2,m}}
& & 
(a_1 . a_2) . \mor m
\ar[d]^{\deltal_{a_1. a_2,m}}
\\
a_1 . \mor(a_2 . m)
\ar[r]_{\deltal_{a_1, a_2 . m}}
&  
\mor(a_1 . (a_2 . m))
\ar[r]_{\mor(\assml_{a_1,a_2,m})}
&
\mor((a_1 . a_2) . m)
}
$$
commutes.
\end{tag}
Note that Axiom \ref{cohmdmor2b} is obtained from
Axiom \ref{cohmdmor2h} by evaluation at the generator 
$\gen$ since the component in $\gen$ of ${\lhpp_{\mor}}$ is an 
identity. By Remark \ref{evgen2cel} Axioms \ref{cohmdmor2h}
and \ref{cohmdmor2b} are equivalent.\\

Eventually we give definitions of $\A$-modules and
their morphisms using the tensor in $\BS$. This will
show that in which sense $\A$-modules occur as algebras for the 
doctrine $\A \otimes -$ over $\BS$.\\  

A $\A$-module $(\M, \mult, \assmd, \lmd)$ consists of 
a {\em strict} arrow $\mult: \A \otimes \M \rightarrow \M$ in $\BS$
in with 2-cells $\assmd$ and  $\lmd$ in $\BS$ as follows\\
\begin{tag}\label{assmd}
$$\xymatrix{
(\A \otimes \A) \otimes \M
\ar[rr]^-{A'}
\ar[d]_{\comp \otimes 1}
& 
& 
\A \otimes ( \A  \otimes \M )
\ar[d]^{1 \otimes \mult}
\ar@{=>}[lld]_{\assmd}
&
\\
\A \otimes \M 
\ar[rd]_{\mult}
& 
& 
\A \otimes \M
\ar[ld]^{\mult}
\\
& 
\M
&
}
$$
\end{tag}
and
\begin{tag}\label{lmd}
$$
\xymatrix{
\M
\ar[r]^-{L'}
\ar@{-}[dd]_{id}
& 
\unc \otimes \M
\ar[dd]^{\unit \otimes 1}
\\
\ar@{=>}[r]^{ \lmd }
&
\\
\M
& 
\A \otimes \M
\ar[l]^-{\mult}
}
$$
\end{tag}
satisfying the coherence conditions \ref{cohmdobj1p},
\ref{cohmdobj2p} and \ref{cohmdobj3p} below.
\begin{tag}\label{cohmdobj1p}
The 2-cells
{\tiny
$$ 
\xymatrix{
((\A \A) \A) \M 
\ar[r]^{A'}
\ar[d]_{(\comp \otimes 1) \otimes 1}
\ar@{}[rd]|{=}
&
(\A \A)  (\A \M)
\ar@{-}[r]^{id}
\ar[d]^{\comp \otimes 1}
\ar@{}[rdd]|{=}
&
(\A \A) (\A \M)
\ar[r]^{A'}
\ar[d]^{1 \otimes \mult}
\ar@{}[rd]|{=}
&
\A ( \A (\A \M))
\ar[d]^{1 \otimes (1 \otimes \mult)}
\\
(\A \A)  \M
\ar[r]^{A'}
\ar[d]_{\comp \otimes 1}
&
\A  (\A \M)
\ar[d]^{1 \otimes \mult}
\ar@{=>}[ldd]_{\assmd}
&
(\A \A)  \M
\ar[r]^{A'}
\ar[d]^{\comp \otimes 1}
& 
\A  (\A  \M)
\ar[d]^{1 \otimes \mult}
\ar@{=>}[ldd]_{\assmd}
\\
\A  \M
\ar[d]_\mult
&
\A \M
\ar[d]_\mult
\ar@{-}[r]_{id}
\ar@{}[rd]|{=}
& 
\A \M
\ar[d]_\mult
&
\A  \M
\ar[d]^\mult\\
\M
\ar@{-}[r]_{id}
&
\M
\ar@{-}[r]_{id}
&
\M
\ar@{-}[r]_{id}
&
\M
}
$$
}
and 
{\tiny
$$
\xymatrix{
((\A \A)  \A) \M
\ar[r]^{A' \otimes 1}
\ar[d]_{(\comp \otimes 1) \otimes 1}
&
((\A  (\A  \A))  \M
\ar[r]^{A'}
\ar[d]|{(1 \otimes \comp) \otimes 1}
\ar@{}[rd]|{=}
\ar@{=>}[ldd]|{\assp \otimes 1}
&
\A  ((\A \A) \M)
\ar[r]^{1 \otimes A'}
\ar[d]|{1 \otimes (\comp \otimes 1)}
& 
\A  ( \A  (\A  \M))
\ar[d]^{1 \otimes (1 \otimes \mult)}
\ar@{=>}[ldd]_{1 \otimes \assmd}
\\
(\A  \A) \M
\ar[d]_{\comp \otimes 1}
&
(\A  \A)  \M
\ar[d]|{\comp \otimes 1}
\ar[r]^{A'}
&
\A  (\A  \M)
\ar[d]|{1 \otimes \mult}
\ar@{=>}[ldd]_{\assmd}
&
\A  (\A  \M)
\ar[d]^{1 \otimes \mult}
\\
\A  \M
\ar[d]_\mult
\ar@{-}[r]^{id}
\ar@{}[rd]|{=}
&
\A  \M
\ar[d]_\mult
& 
\A  \M
\ar@{-}[r]^{id}
\ar[d]_{\mult}
\ar@{}[rd]|{=}
&
\A  \M
\ar[d]^{\mult}
\\
\M
\ar@{-}[r]_{id}
&
\M
\ar@{-}[r]_{id}
&
\M
\ar@{-}[r]_{id}
&
\M
}
$$
}
are equal.
\end{tag}

\begin{tag}\label{cohmdobj2p}
The 2-cells  
$\Xi_1$ $=$
{\small
$$\xymatrix{
\A \otimes  \M
\ar[d]_{R' \otimes 1} 
\ar@{-}[rr]^{id}
&
\ar@{=>}[d]^{\rp \otimes 1}
&
\A \otimes \M
\ar[r]^-{\mult}
&
\M
\\
(\A \otimes \unc) \otimes \M
\ar[rr]_{(1 \otimes \unit) \otimes 1}
&
&
(\A \otimes \A) \otimes \M
\ar[u]_{\comp \otimes 1} 
}
$$
}\\
$\Xi_2$ $=$ 
{\small
$$\xymatrix{
\A \otimes \M
\ar[d]_{1 \otimes L'} 
\ar@{-}[rr]^{id}
&
\ar@{=>}[d]^{1 \otimes \lmd}
&
\A \otimes \M
\ar[r]^-{\mult}
&
\M
\\
\A \otimes (\unc \otimes \M)
\ar[rr]_{1 \otimes (\unit \otimes 1)}
&
&
\A \otimes (\A \otimes \M)
\ar[u]_{1 \otimes \mult} 
}$$
}
and\\ 
$\Xi_3$ $=$
{\tiny
$$
\xymatrix@C=3pc{
& 
(\A \otimes \unc) \otimes \M
\ar[r]^-{(1 \otimes \unit) \otimes 1}
\ar[dd]^{A'}
&
(\A \otimes \A) \otimes \M
\ar[r]^-{\comp \otimes 1}
\ar[dd]_{A'}
&
\A \otimes \M
\ar[rd]^{\mult}
\\
\A \otimes \M
\ar[ru]^{R' \otimes 1}
\ar[rd]_{1 \otimes L'}
\ar@{}[r]|-{=}
&
\ar@{}[r]|-{=}
&
&
&
\M
\\
& 
\A \otimes (\unc \otimes \M)
\ar[r]_-{1 \otimes (\unit \otimes 1)}
&
\A \otimes (\A \otimes \M)
\ar[r]_-{1 \otimes \mult}
\ar@{=>}[ruu]^{\assmd}
&
\A \otimes \M
\ar[ru]_{\mult}
}
$$
}
satisfy the equality
$\Xi_1 = \Xi_3 \circ \Xi_2$.
\end{tag}
\begin{tag}\label{cohmdobj3p}
The 2-cells
$\Xi_1$ $=$
{\tiny
$$\xymatrix{
\A \otimes \M
\ar[d]_{L' \otimes 1}
\ar@{-}[rr]^{id}
&
\ar@{=>}[d]^{\lp \otimes 1} 
& 
\A \otimes  \M
\ar[rr]^{\varphi}
& & 
\M
\\ 
(\unc \otimes \A) \otimes \M
\ar[rr]_-{(\unit  \otimes 1) \otimes 1}
& 
& 
(\A \otimes \A) \otimes \M
\ar[u]_-{\comp \otimes 1}
}$$
}\\
$\Xi_2$ $=$
{\tiny
$$\xymatrix{
\A \otimes \M
\ar@{}[rrd]|{=}
\ar[rr]^{\mult}
\ar[d]_{L'}
& & 
\M
\ar[d]_{L'}
\ar@{-}[rr]^{id}
&
\ar@{=>}[d]^{\lmd}
&
\M
\\
\unc \otimes (\A \otimes \M)
\ar[rr]_{1 \otimes \mult}
& & 
\unc \otimes \M
\ar[rr]_{\unit \otimes 1}
& 
&
\A \otimes \M
\ar[u]_{\mult}
}
$$}\\
$\Xi_3$ $=$
{\tiny
$$
\xymatrix{
\A \otimes \M 
\ar[rr]^{L'}
\ar[rd]_{L' \otimes 1}
&
&
\unc \otimes (\A \otimes \M)
\ar[r]^{\unit \otimes 1}
&
\A \otimes (\A \otimes \M)
\ar[r]^{1 \otimes \mult}
&
\A \otimes \M
\ar[r]^{\mult}
&
\M
\\
& 
(\unc \otimes \A) \otimes \M
\ar[ru]_{A'}
\ar@{=>}[u]_{\theta}
}
$$
}
where the ``canonical'' 2-cell $\theta$ is defined in Appendix 
in \ref{def2ctheta} and has image by $\Res$ an identity,
and\\
$\Xi_4$ $=$
{\tiny
$$\xymatrix{
\A \otimes \M
\ar[r]^-{L' \otimes 1}
&
(\unc \otimes \A) \otimes \M
\ar[r]^-{(\unit \otimes 1) \otimes 1}
\ar[dd]_{A'}
\ar@{}[rdd]|{=}
&
(\A \otimes \A) \otimes \M
\ar[r]^-{\mult \otimes 1}
\ar[dd]_{A'}
&
\A \otimes \M
\ar[rd]^-{\varphi}
\\
& & & 
&
\M
\\
&
\unc \otimes (\A \otimes \M)
\ar[r]_{\unit \otimes 1}
 & 
\A \otimes (\A \otimes \M)
\ar[r]_-{1 \otimes \varphi}
\ar@{=>}[ruu]^{\assmd}
&
\A \otimes \M
\ar[ru]_-{\varphi}
}
$$}\\
satisfy the equality
$\Xi_1 = \Xi_4 \circ {(\Xi_3)}^{-1} \circ \Xi_2$.
\end{tag}

With the previous definition of $\A$-modules,
a morphism of $\A$-modules $(\mor,\deltad): 
(\M,\mult,\assmd, \lmd) \rightarrow (\N,\psi,\assmd,\lmd)$ 
consists of 
an arrow $\mor:\M \rightarrow \N$  
with a 2-cell $\deltad$ in $\BS$
\begin{tag}\label{deltad}
$$
\xymatrix{
\A \otimes \M 
\ar[r]^{1 \otimes \mor}
\ar[d]_{\mult}
& 
\A \otimes \N
\ar@{=>}[ld]
\ar[d]^{\psi}
\\
\M 
\ar[r]_{\mor}
& 
\N
}
$$
\end{tag}
that satisfy Axioms \ref{cohmdmor1p} and \ref{cohmdmor2p} below.
\begin{tag}\label{cohmdmor1p}
The 2-cells
$$\xymatrix{
&
\A  (\A  \M)
\ar[r]^{1 \otimes (1 \otimes \mor)}
&
\A  (\A  \N)
\ar[dd]^{1 \otimes \psi}
\ar@{=>}[lddd]^{\assmd}
\\
(\A  \A)  \M
\ar[r]^{1 \otimes \mor}
\ar[ru]^{A'}
\ar[d]_{\comp \otimes 1}
&
(\A  \A) \N
\ar[d]^{\comp \otimes 1}
\ar[ru]_{A'}
\ar@{}[u]|{=}
\ar@{}[ld]|{=}
\\
\A \M
\ar[d]_{\mult}
\ar[r]^-{1 \otimes \mor}
& 
\A  \N
\ar@{=>}[ld]_{\deltad}
\ar[d]_{\psi}
& 
\A  \N
\ar[ld]^{\psi}
\\
\M
\ar[r]_\mor
&
\N
}$$
and
$$
\xymatrix{
(\A  \A)  \M
\ar[d]_{\comp \otimes 1}
\ar[r]^-{A'}
&
\A  (\A  \M)
\ar[d]|{1 \otimes \mult}
\ar[r]^-{1 \otimes (1 \otimes \mor)}
\ar@{=>}[ldd]_{\assmd}
&
\A  (\A  \N)
\ar[d]^{1 \otimes \psi}
\ar@{=>}[ld]|{1 \otimes \deltad}
\\
\A  \M
\ar[d]_{\mult}
&
\A  \M
\ar[d]|{\mult}
\ar[r]^-{1 \otimes \mor}
&
\A \N
\ar@{=>}[ld]^{\deltad}
\ar[d]^{\psi}
\\
\M 
\ar@{-}[r]_-{id}
&
\M
\ar[r]_-{\mor}
& 
\N
}
$$
are equal.
\end{tag}
\begin{tag}\label{cohmdmor2p}
The 2-cells 
$$
\xymatrix{
& 
\ar@<8ex>@{=>}[d]^{\lmd}
\\
\M
\ar[r]_-{L'}
\ar@/^38pt/[rrr]^1
&
\unc \otimes \M
\ar[r]_-{\unit \otimes 1}
&
\A \otimes \M
\ar[r]_-{\varphi}
&
\M
\ar[r]_-{\mor}
&
\N
}
$$
and 
$$
\xymatrix{
&
\ar@<6ex>@{=>}[d]^{\lmd} 
& 
&
\\
\N
\ar[r]^-{L'}
\ar@/^38pt/[rrr]^1
&
\unc \otimes \N
\ar[r]_-{\unit \otimes 1}
&
\A \otimes \N
\ar[r]^{\psi}
\ar@{=>}[rd]^{\deltad}
&
\N
\\
\M
\ar[u]^{\mor}
\ar[r]_{L'}
\ar@{}[ru]|{=}
&
\unc \otimes \M
\ar[r]_-{\unit \otimes 1}
\ar[u]|{1 \otimes \mor}
\ar@{}[ru]|{=}
&
\A \otimes \M
\ar[r]_-{\mult}
\ar[u]|{1 \otimes \mor}
&
\M
\ar[u]_{\mor}
}
$$
are equal.
\end{tag}

To justify these new definitions let us consider again 
an arbitrary $\BS$-functor $F: \A \rightarrow \BS$
which defines a $\A$-module $\M$ with multiplication
$\multh: \A \rightarrow [\M,\M]$ and 2-cells
$\assmh$ as in \ref{assmh} and $\lmz$ as in \ref{lmz}.
The multiplication corresponds by adjunction 
\ref{punitens} to a strict arrow 
$\mult: \A \otimes \M \rightarrow \M$.
According to the adjunction \ref{punitens}, 
Lemmas \ref{calc1} and \cite{Sch08}-\calcdeux, 
the map $\Res \circ \Res$ defines a bijection between the sets 
of 2-cells of the following kinds
$$
\xymatrix{
(\A\A)\M 
\ar[rr]^-{A'}
\ar[d]_{\comp \otimes 1}
& 
& 
\A(\A\M)
\ar[d]^{1 \otimes \mult}
\ar@{=>}[lld]
\\
\A\M
\ar[rd]_{\mult}
& & 
\A\M
\ar[ld]^{\mult}
\\
& 
\M
&
}
$$
and
$$
\xymatrix{
\A 
\ar[r]^{\multh}
\ar[d]_{\comh}
& 
[\M,\M]
\ar[r]^-{[\M,-]}
& 
[[\M, \M],[\M,\M]]
\ar[d]^{[\multh,1]}
\ar@{=>}[lld]
\\
[\A,\A]
& 
&
[\A,[\M,\M]]
\ar[ll]^{[1,\multh]}
}
$$
and one has a 2-cell $\assmd$ corresponding to 
$\assmh$/${F'}^2_{*,*,*}$ via the above bijection.
Note then that $\assmd$ has image by $\Res$ the 2-cell
$F^2_{\gen,\gen,\gen}$ which we also write
$\assm$. 
One obtains a 2-cell $\lmd$ as in \ref{lmd}
that is equal to the 2-cell $\lmh$ according Lemma \ref{calc6}.\\
 
We are going to check that the points \ref{ptm1}, \ref{ptm2},
\ref{ptm3}, \ref{ptm4}, \ref{ptm5} and \ref{ptm6} below hold
for a $\BS$-functor $F$ and the related data as above.
Since the arrows domains of the 2-cells
of the equalities of 
Axioms \ref{cohmdobj1p}, \ref{cohmdobj2p} and \ref{cohmdobj3p}
are strict, it will result from the adjunction \ref{punitens}
that these axioms are equivalent respectively to
Axioms \ref{cohfun1p}, \ref{cohfun2p} and \ref{cohfun3p}
for the $\BS$-functor $F$.
\begin{tag}\label{ptm1}
The 2-cell
{\small
$$\xymatrix{
(\A \otimes \A) \otimes \A
\ar[d]_{(\multh \otimes \multh) \otimes \multh}
\ar[r]^-{\comp \otimes 1}
& 
\A \otimes \A
\ar[d]|{\multh \otimes \multh}
\ar[r]^-{\comp}
& 
\A
\ar[d]^{\multh}
\\
([\M,\M] \otimes [\M,\M]) \otimes [\M,\M]
\ar[d]_{A'}
\ar[r]_-{\comp \otimes 1}
\ar@{=>}[ru]^{\assm \otimes 1}
&
[\M,\M] \otimes [\M,\M]
\ar[r]_-{\comp}
\ar@{=>}[ru]^{\assm}
&
[\M,\M]
\ar[d]^{id}
\\
[\M,\M] \otimes ([\M,\M] \otimes [\M,\M])
\ar[r]_-{1 \otimes \comp}
\ar@{}[rru]|{=}
&
[\M,\M] \otimes [\M,\M]
\ar[r]_-{\comp}
&
[\M,\M]
} 
$$
}
is the image by $\Res$
of the first of the 2-cells of Axiom \ref{cohmdobj1p}.
\end{tag}
\pf See \ref{pfptm1} in Appendix. \epf
\begin{tag}\label{ptm2}
The 2-cell 
{\small
$$
\xymatrix{
(\A\A)\A 
\ar[r]^{\comp \otimes 1}
\ar[d]_{A'}
&
\A \A
\ar[r]^{\comp}
&
\A
\ar@{-}[d]^{id}
\\
\A (\A\A)
\ar@{=>}[rru]^{\assp}
\ar[r]_{1 \otimes \comp}
\ar[d]_{\multh \otimes (\multh \otimes \multh)}
&
\A \A
\ar[d]|{\multh \otimes \multh}
\ar[r]^{\comp}
& 
\A
\ar[d]^{\multh}
\\
[\M,\M] 
\otimes 
([\M,\M] \otimes [\M,\M])
\ar@{=>}[ru]_{1 \otimes \assm}
\ar[r]_-{1 \otimes \comp}
&
[\M,\M] \otimes [\M,\M]
\ar@{=>}[ru]^{\assm}
\ar[r]_-{\comp}
&
[\M,\M]
}
$$ 
}
is the image by $\Res$ of the second 2-cell
of Axiom \ref{cohmdobj2p}
\end{tag}
\pf See \ref{pfptm2} in Appendix. \epf
\begin{tag}\label{ptm3}
The 2-cell
{\small
$$\xymatrix{
\A 
\ar[rr]^-{id}
\ar[d]_{R'}
& 
\ar@{=>}[d]^{\rp}
& 
\A 
\ar[rr]^-{\multh}
& & 
[\M,\M]
\\
\A \otimes \unc 
\ar[rr]_-{1 \otimes \unit}
&& 
\A \otimes \A
\ar[u]_{\comp}
}$$
}
is the image by $\Res$ of the 2-cell
$\Xi_1$ of Axiom \ref{cohmdobj2p}.
\end{tag}
\pf Immediate. \epf
\begin{tag}\label{ptm4}
The 2-cell
{\small
$$\xymatrix{
\A 
\ar@{}[rd]|{=}
\ar[r]^{\multh}
\ar[d]_{R'}
&
[\M,\M]
\ar[d]^{R'}
\ar@/^70pt/[ddd]^{id}
\\
\A \otimes \unc
\ar[d]_{1 \otimes \unit} 
\ar[r]^{\multh \otimes 1}
&
[\M,\M] \otimes \unc
\ar@{=>}[ld]_{1 \otimes \lmz}
\ar[d]^{1 \otimes \vv}
&
\\
\A \otimes \A
\ar[r]_-{\multh \otimes \multh}
\ar[d]_{\comp}
&
[\M,\M] \otimes [\M,\M]
\ar@{=>}[ld]_{\assm}
\ar[d]^{\comp}
\ar@{}[ru]|{=}
\\
\A 
\ar[r]_{\multh}
&
[\M,\M] 
}$$
}
is the image by $\Res$ of the 2-cell
$\Xi_3 \circ \Xi_2$ of Axiom \ref{cohmdobj2p}.
\end{tag}
\pf See \ref{pfptm4} in Appendix. \epf
\begin{tag}\label{ptm5}
The 2-cell
{\small
$$\xymatrix{
\A 
\ar[rr]^-{id}
\ar[d]_{L'}
& 
\ar@{=>}[d]^{\lp}
& 
\A 
\ar[rr]^-{\multh}
& & 
[\M,\M]
\\
\unc \otimes \A 
\ar[rr]_-{\unit \otimes 1}
&& 
\A \otimes \A
\ar[u]_{\comp}
}$$
}
is the image by $\Res$ of the 2-cell 
$\Xi_1$ of Axiom \ref{cohmdobj3p}.
\end{tag}
\pf immediate.
\epf
\begin{tag}\label{ptm6}
The 2-cell
{\small
$$\xymatrix{
\A 
\ar@{}[rd]|{=}
\ar[r]^{\multh}
\ar[d]_{L'}
&
[\M,\M]
\ar[d]^{L'}
\ar@/^70pt/[ddd]^{id}
\\
\unc \otimes \A
\ar[d]_{\unit \otimes 1} 
\ar[r]^{1 \otimes \multh}
&
\unc \otimes [\M,\M] 
\ar@{=>}[ld]_{\lmz \otimes 1}
\ar[d]^{\vv \otimes 1}
&
\\
\A \otimes \A
\ar[r]_-{\multh \otimes \multh}
\ar[d]_{\comp}
&
[\M,\M] \otimes [\M,\M]
\ar@{=>}[ld]_{\assm}
\ar@{}[ru]|{=}
\ar[d]^{\comp}
\\
\A 
\ar[r]_{\multh}
&
[\M,\M] 
}$$
}
is the image by $\Res$ of the 2-cell 
$\Xi_4 \circ {(\Xi_3)}^{-1} \circ \Xi_2$ 
of Axiom \ref{cohmdobj3p}.
\end{tag}
\pf See \ref{pfptm6} in Appendix. \epf

Let us consider now two modules
$\M$ and $\N$ with respective multiplications denoted by 
$\mult: \A \otimes \M \rightarrow \M$/$\multh: \A \rightarrow [\M,\M]$
and $\psi: \A \otimes \N \rightarrow \N$/$\psi': \A \rightarrow [\N,\N]$,
with sets of coherence two-cells for both written
$\assmd$/$\assm$/$\assmh$ and $\lmd$/$\lmh$/$\lmz$,
and an arrow $H: \M \rightarrow \N$ in $\BS$. 
A 2-cell $\deltad$ as in \ref{deltad}
corresponds by the adjunction \ref{punitens}
to a 2-cell $\deltah$.
For the above related data the equivalence of Axioms
\ref{cohmdmor1p} and \ref{cohmdmor1h} 
results form the two following points \ref{ptm7} and \ref{ptm8}
below whereas the equivalence \ref{cohmdmor2p} and \ref{cohmdmor2h} 
follows from Remark \ref{evgen2cel} and points
\ref{ptm9} and \ref{ptm10}. 
\begin{tag}\label{ptm7}
The first of the 2-cell 
of \ref{cohmdmor1p} has image by $\Res \circ \Res$
 the pasting $\Xi_2 \circ \Xi_1$ of \ref{cohmdmor1h}, 
it has a strict domain which has a strict image by $\Res$.
\end{tag}
\pf See \ref{pfptm7} in Appendix. \epf
\begin{tag}\label{ptm8}
The image by $\Res \circ \Res$ of the second 2-cell
of Axiom \ref{cohmdmor1p} is the pasting 
$\Xi_7 \circ \Xi_6 \circ \Xi_5
\circ \Xi_4 \circ \Xi_3 \circ \Xi_2 \circ \Xi_1$
of \ref{cohmdmor1h}.
\end{tag}
\pf See \ref{pfptm8} in Appendix. \epf 
\begin{tag}\label{ptm9}
The first 2-cell of Axiom \ref{cohmdmor2h}
has image by $\ev_{\gen}$ the first 2-cell of Axiom \ref{cohmdmor2p}
and has a strict domain.
\end{tag}
\pf See \ref{pfptm9} in Appendix. \epf
\begin{tag}\label{ptm10}
The second 2-cell of Axiom \ref{cohmdmor2h} has image 
by $\ev_{\gen}$ the second cell of Axiom \ref{cohmdmor2p}.
\end{tag}
\pf See \ref{pfptm10} in Appendix. \epf

For any 2-ring $\A$, a $\A$-module $\M$ is
said {\em strict} when the corresponding $\BS$-presheaf
is strict which is to say that its action
$\multh: \A \rightarrow [\M,\M]$ is strict 
and the 2-cells $\assmh$ and $\lmz$ are identities.
A few remarks are in order. Consider any $\A$-module $\M$.
If its 2-cell $\lmz$ is an identity then certainly
its 2-cell $\lmh$ is also an identity. Conversely if 
the action $\multh$ is strict then for any $m$ in $\M$, 
$\multh^*(m): \A \rightarrow \M$ is strict,
and the component $\epsilon_{\multh^*(m) \circ \unit}$ is an 
identity and from this fact one has that if $\lmh$ is strict
then also is $\lmz$.
If the 2-cells $\assmd$ are identities
then certainly are the 2-cells
$\assmh$.  Conversely if $\A$ is a strict 2-ring 
and the 2-cells $\assmh$ 
are identities then the 2-cells 
$\assmd$ also are.\\

\begin{proposition}\label{ImodeqBS}
The forgetful functor $\unc-Mod \rightarrow \BS$ is part of an
equivalence of categories. Its equivalence inverse
factors as 
$$\xymatrix{\BS \ar[r]^-{\cong} & \unc-Mod^s \ar[r]^-{inc} & \unc-Mod}$$
where the left functor is an isomorphism between $\BS$ 
and the full sub-category $\unc-Mod^s$ of $\unc-Mod$ generated
by the strict modules and $inc$ is the inclusion functor. 
\end{proposition}
\pf
One has a forgetful 2-functor 
$\unc-Mod \rightarrow \BS$ and the result follows
then from Lemmas \ref{objPic}, \ref{uniqmor}
and \ref{essurj} below.
\epf

\begin{lemma}\label{objPic}
Any object $\A$ of $\BS$ admits a unique
strict $\unc$-module structure, its
multiplication is given by the arrow
$\multh = \vv: \unc \rightarrow [\A,\A]$ 
(or equivalently $\mult = L_{\A}: \unc \otimes \A \rightarrow \A$). 
\end{lemma}
\pf
If $\A$ has a strict $\unc$-module structure 
with multiplication $\multh = \vv : \unc \rightarrow [\A,\A]$
the 2-cell $\lmh$ (\ref{lmh}) in this case being an identity one has 
that the composite 
$$
\xymatrix{
\A 
\ar[r]^-{{\multh}^*}
&
[\unc, \A]
\ar[r]^-{\ev_{\gen}}
&
\A
}
$$ in $\BS$
is necessarily the identity at $\A$.
Actually there is a unique arrow $f: \A \rightarrow [\unc,\A]$ in $\BS$
with strict images in $[\unc,\A]$ --  or equivalently such that $f^*$
is strict -- and such that the composite
$$\xymatrix{
\A
\ar[r]^-f
&
[\unc,\A]
\ar[r]^-{\ev_{\gen}}
&
\A
}$$ is the identity and this arrow is $\vs$.\\

Now to establish that $\multh = \vv: \unc \rightarrow [\A,\A]$ 
gives a strict $\unc$-module structure on $\A$,
it remains to check that one has in this case an identity 
2-cell $\assmh$ (\ref{assmh}) 
which is the commutativity of
the external diagram in the pasting below
$$
\xymatrix{
\unc 
\ar[rd]^{\vv}
\ar[r]^-{\vv}
\ar[dd]_{\vv}
&
[\A,\A]
\ar[d]^{[\A,-]}
\\
&
[[\A,\A],[\A,\A]]
\ar[d]^{[\vv,1]}
\\
[\unc,\unc]
\ar[r]_{[1,\vv]}
&  
[\unc,[\A,\A]]
}
$$ 
where the top right diagram commutes
according to Lemma \cite{Sch08}-\lemvvun and the bottom
left also does according to Lemma \cite{Sch08}-\lemvvdeux.
\epf

\begin{lemma}\label{uniqmor}
For any $\unc$-module $\A$
and any arrow $H: \A \rightarrow \B$ in $\BS$
there is a unique $\unc$-module morphism
from $\A$ to the strict $\unc$-module structure
on the symmetric Picard category $\B$ which underlying 
map in $\BS$ is $H$.
If the multiplication of $\A$ is given by 
$\multh: \unc \rightarrow [\A,\A]$ 
this morphism has 2-cell $\deltah$ as in \ref{deltah}
$$
\xymatrix{
\unc 
\ar[r]^-{\vv}
\ar[d]_{\multh}
&
[\B,\B]
\ar[d]^{[H,1]}
\ar@{=>}[ld]
\\
[\A,\A]
\ar[r]_-{[1,H]}
&
[\A,\B],
}
$$
which is determined by its value in $\gen$ by Remark \ref{evgen2cel} 
and such that
\begin{center} 
$(\deltah_{\gen})_a$  $=$ 
$\xymatrix{ 
\gen.Ha 
\ar@{-}[r]^{id}
&
Ha 
\ar[r]^-{H(\lml_a)}
&
H( \multh(\gen)(a) )}.$
\end{center}
\end{lemma}
\pf
Let us write $t \times a$ for $\multh(t)(a)$ for 
any objects $t$ of $\unc$ and $a$ of $\A$.\\
 
The coherence Axiom \ref{cohmdmor1b} for the pair $H$ and
$\deltah$, with corresponding bilinear 
$\deltal$ as in \ref{deltal},
amounts to the commutativity of the diagram in $\B$
{\small
$$
\xymatrix{
Ha 
\ar[rr]^{H(\lml_a)}
\ar@{-}[rdd]_{id}
& & 
H(\gen \times a)
\\
\\
& \gen . Ha
\ar[ruu]_{\deltal_{\gen,a}}.
}
$$ 
}
That the arrow  $H:\A \rightarrow \B$ in $\BS$ together 
with the 2-cell $\deltal$ defined by the condition above 
satisfies Axiom \ref{cohmdmor2b} amounts to the commutativity
of the diagram in $\B$
\begin{tag}\label{dicoass}
{\small
$$
\xymatrix@C=3pc{
t_1.(t_2.Ha)
\ar@{-}[rr]^{id}
\ar[d]_{t_1.\deltal_{t_2,a}}
& & 
(t_1.t_2).Ha
\ar[d]^{\deltal_{t_1.t_2,a}}
\\
t_1 . H(t_2 \times a)
\ar[r]_{\deltal_{t_1, t_2 \times a}}
& 
H( t_1 \times (t_2 \times a) )
\ar[r]_{H(\assml_{t_1,t_2,a})}
& 
H((t_1.t_2) \times a)
}
$$
}
\end{tag}
for all objects $t_1$, $t_2$ in $\unc$ 
and $a$ in $\A$.\\

We prove this last point by induction on the structure of the 
objects $t_1$ in $\unc$ for arbitrary objects $t_2$ and $a$.\\ 
 
For $t_1 = I$ diagram \ref{dicoass}
is the external diagram in the pasting
$$
\xymatrix{
I
\ar@{-}[rr]^{id}
\ar@{-}[dd]_{id} 
& 
& 
I
\ar[dd]^{{\deltal}_{I,a}}
\ar[ld]_{H^0}
\\
& 
H(I)
\ar[d]|{H({\multh^0}_{t_2 \times a})}
\ar[rd]^{H({\multh^0}_a)}
\\
I
\ar[ru]^{H^0}
\ar[r]_-{\deltal_{I,t_2 \times a}} 
& 
H( I \times (t_2 \times a) )
\ar[r]_-{H(\assml_{I,t_2,a})}
& 
H( I \times a) 
}
$$
in which all diagrams commute and in particular 
the bottom-right triangle since 
the natural transformation
$\assml_{-,t_2,a}: - \times (t_2 \times a) \rightarrow (- . t_2) \times a:
\unc \rightarrow \A$ is monoidal.\\  

For $t_1 = \gen$ diagram \ref{dicoass}
is the external diagram in the pasting
$$
\xymatrix{
\gen.( t_2.Ha )
\ar@{-}[rr]^{id}
\ar[d]_{\gen.\deltal_{t_2,a}} 
& 
& 
(\gen.t_2).Ha
\ar[d]^{\deltal_{\gen.t_2,a}}
\\
\gen.H(t_2 \times a)
\ar@{-}[d]_{id}
\ar@{-}[rr]^{id}
& 
& 
H((\gen.t_2) \times a)
\\
H(t_2 \times a)
\ar[rr]_{H(\lml_{t_2 \times a})} 
& 
&
H(\gen \times (t_2 \times a))
\ar[u]_{H(\assml_{\gen,t_2,a})}
}
$$
where the bottom diagram commutes
by the coherence Axiom \ref{cohmdobj3b} for 
the $\unc$-module $\A$.\\ 

For $t_1 = t_1' \otimes t_1''$ diagram \ref{dicoass}
is the external diagram in the pasting
$$
\xymatrix{
(t_1' \otimes t_1'') . (t_2 . Ha)
\ar@{-}[d]_{id}
\ar@{-}[rr]^{id}
& & 
((t_1' \otimes t_1'') . t_2).Ha
\ar@{-}[d]^{id}
\\
(t_1'.(t_2.Ha))
\otimes 
(t_1''. (t_2 . Ha))
\ar@{-}[rr]^{id}
\ar[d]_{t_1'.\deltal_{t_2,a} \otimes t_1''.\deltal_{t_2,a}}
& & 
((t_1'.t_2). Ha)
\otimes 
((t_1''.t_2) . Ha)
\ar[d]^{\deltal_{t_1'.t_2,a} \otimes \deltal_{t_1''.t_2,a}}
\\
(t_1'.H(t_2 \times a)) \otimes (t_1''.H(t_2 \times a))
\ar[d]_{\deltal_{t_1',t_2 \times a} \otimes \deltal_{t_1'', t_2 \times a} }
& & 
H((t_1'.t_2) \times a) \otimes H((t_1''.t_2) \times a))
\ar[d]^{H^2_{ (t_1'.t_2) \times a,  (t_1''.t_2) \times a }  }
\\
H(t_1' \times (t_2 \times a)) \otimes H(t_1'' \times (t_2 \times a))
\ar[d]_{H^2_{ t_1' \times (t_2 \times a), t_1'' \times (t_2 \times a) }}
& & 
H( ((t_1'.t_2) \times a) \otimes ((t_1''.t_2) \times a) )
\ar[d]^{ H(( {{\multh}^2_{t_1'.t_2,t_1''.t_2})}_a ) }
\\
H( (t_1' \times (t_2 \times a)) \otimes (t_1'' \times (t_2 \times a)) )
\ar[d]_{H(( {\multh^2_{t_1',t_1''})}_{t_2 \times a} )}
\ar[rru]^{H(\assml_{t_1',t_2,a}) \otimes H(\assml_{t_1'',t_2,a})}
& & 
H ( ((t_1'.t_2) \otimes (t_1''.t_2)) \times a )
\ar@{-}[d]^{id}
\\  
H( (t_1' \otimes t_1'') \times (t_2 \times a) )
\ar[rr]_{H( \assml_{t_1' \otimes t_1'',t_2,a} )}
& & 
H ( ((t_1' \otimes t_1'').t_2) \times  a )
}
$$
where the middle diagram is commutative if the 
diagram \ref{dicoass} commutes for the values $t_1 = t_1'$ 
and  $t_1 = t_1''$ 
and the bottom diagram commutes since
the natural transformation
$\assml_{-,t_2,a}: - \times (t_2 \times a) \rightarrow (- . t_2) \times a:
\unc \rightarrow \A$ is monoidal.\\  

For $t_1 = t^{\eqi}$, the diagram \ref{dicoass}
is 
$$
\xymatrix{
t^{\eqi}.(t_2.Ha)
\ar@{-}[rr]^{id}
\ar[d]_{ t^{\eqi}.\deltal_{t_2,a} }
& & 
(t^{\eqi}.t_2).Ha
\ar[d]^{\deltal_{t^{\eqi}.t_2,a}}
\\
t^{\eqi}.H(t_2 \times a)
\ar[r]_{\deltal_{ t^{\eqi}, t_2 \times a } }
&
H( t^{\eqi} \times (t_2 \times a) )
\ar[r]_{ H(\assml_{t^{\eqi},t_2,a }) }
&
H((t^{\eqi}.t_2) \times a)
}
$$
Note that according to Lemmas \ref{sigmap} and \ref{pcanco}
the arrow $\delta_{t^{\eqi},a}$ is
$$\xymatrix{t^{\eqi}.Ha 
\ar@{-}[r]
&
(t.Ha)^{\eqi}
&
{H(t \times a)}^{\eqi}
\ar@{-}[r]^-{\cong}
\ar[l]_-{{(\deltal_{t,a})}^{\eqi}}
&
H((t \times a)^{\eqi})
\ar@{-}[r]^-{H(\cong)}
&
H(t^{\eqi} \times a).}$$
Therefore the left-bottom leg rewrites
\begin{tabbing}
\=1. \={\small $\xymatrix{
t^{\eqi}.(t_2.Ha)
\ar[r]^-{ t^{\eqi}.\deltal_{t_2,a} }
&
t^{\eqi}.H(t_2 \times a)
\ar[r]^-{\deltal_{ t^{\eqi},t_2 \times a }}
&
H(t^{\eqi} \times (t_2 \times a))
\ar[r]^-{ H(\assml_{ t^{\eqi},t_2,a }) }
&
H((t^{\eqi}.t_2) \times a)
}$}\\
\>2. \>{\small $\xymatrix{
t^{\eqi}.(t_2.Ha)
\ar@{-}[r]^-{id}
&
{(t.(t_2.Ha))}^{\eqi}
&
{( t.H(t_2 \times a)) }^{\eqi}
\ar[l]_-{{( t.\deltal_{t_2,a} )}^{\eqi}}
&
{H(t \times (t_2 \times a))}^{\eqi}
\ar[l]_-{  { ( \deltal_{t,t_2 \times a} ) }^{\eqi}  }
\ar@{-}[r]^-{\cong}
&
H({(t \times (t_2 \times a))}^{\eqi} )\;...}$}\\
\> \>{\small $\xymatrix{
...
\ar@{-}[r]^-{H(\cong)} 
&
H(t^{\eqi} \times (t_2 \times a))
\ar[r]^-{ H(\assml_{ t^{\eqi},t_2,a }) }
&
H((t^{\eqi}.t_2) \times a)
}$}\\
\>3. \>{\small $\xymatrix{
t^{\eqi}.(t_2.Ha)
\ar@{-}[r]^-{id}
&
{(t.(t_2.Ha))}^{\eqi}
&
{( t.H(t_2 \times a)) }^{\eqi}
\ar[l]_-{{( t.\deltal_{t_2,a} )}^{\eqi}}
&
{H(t \times (t_2 \times a))}^{\eqi}
\ar[l]_-{  { ( \deltal_{t,t_2 \times a} ) }^{\eqi}  }
\ar@{-}[r]^-{\cong}
&
H({(t \times (t_2 \times a))}^{\eqi} )
\;...}$}\\
\> \>{\small $\xymatrix{...
 &
H( {((t .t_2) \times a )}^{\eqi} )
\ar@{-}[l]_-{ H( { \assml_{t,t_2,a} }^{\eqi} ) }
\ar@{-}[r]^-{ H(\cong) }
&
H((t^{\eqi}.t_2) \times a)
}$}\\
\>4. \>{\small $\xymatrix{
t^{\eqi}.(t_2.Ha)
\ar@{-}[r]^-{id}
&
{(t.(t_2.Ha))}^{\eqi}
&
{( t.H(t_2 \times a)) }^{\eqi}
\ar[l]_-{{( t.\deltal_{t_2,a} )}^{\eqi}}
&
{H(t \times (t_2 \times a))}^{\eqi}
\ar[l]_-{  { ( \deltal_{t,t_2 \times a} ) }^{\eqi}  }
&
{H( ((t.t_2) \times a) )}^{\eqi}\;...
\ar[l]_-{ {(H(\assml_{t,t_2,a}))}^{\eqi} }  }$}
\\
\> \>{\small
$\xymatrix{...
\ar@{-}[r]^-{\cong}
 &
H( {((t .t_2) \times a )}^{\eqi} )
\ar@{-}[r]^-{ H(\cong) }
&
H((t^{\eqi}.t_2) \times a)
}$}
\end{tabbing}
In the derivation above arrows 2. and 3. are equal due
to Lemma \ref{sigmap}
and arrows 3. and 4. are equal due to the naturality
of the isomorphism \ref{isofuntag}.\\

The top-right leg of the diagram above rewrites
\begin{tabbing}
\=1. \={\small
$\xymatrix{
t^{\eqi}.(t_2.Ha)
\ar@{-}[r]^-{id}
&
(t^{\eqi}.t_2).Ha
\ar[r]^-{\deltal_{t^{\eqi}.t_2,a}}
&
H(  ( t^{\eqi}.t_2) \times a )}$}
\\
\>2. \>{\small
$\xymatrix{
t^{\eqi}.(t_2.Ha)
\ar@{-}[r]^-{id}
&
(t^{\eqi}.t_2).Ha
\ar@{-}[r]^-{id}
&
{(t.t_2)}^{\eqi}.Ha
\ar[r]^-{\deltal_{{(t.t_2)}^{\eqi},a}}
&
H( {(t.t_2)}^{\eqi} \times a)
\ar@{-}[r]^-{id}
&
H(  ( t^{\eqi}.t_2) \times a )
}$}\\
\>3. \>{\small$\xymatrix{
t^{\eqi}.(t_2.Ha)
\ar@{-}[r]^-{id}
&
{(t.t_2)}^{\eqi}.Ha
\ar@{-}[r]^-{id}
&
{( (t.t_2).Ha )}^{\eqi}
&
{H((t.t_2) \times a)}^{\eqi}
\ar[l]_-{ {(\deltal_{t.t_2,a})}^{\eqi} }
\ar@{-}[r]^-{\cong}
&
H({ ((t.t_2) \times a) }^{\eqi} )\;...
}$}\\
\> \>
{\small $\xymatrix{
...
\ar@{-}[r]^-{H(\cong)}
&
H({(t.t_2)}^{\eqi} \times a)
\ar@{-}[r]^-{id}
&
H(  ( t^{\eqi}.t_2) \times a )
}$}
\end{tabbing}
Therefore the two legs above are equal
when diagram \ref{dicoass} commutes for $t_1 = t$.
\epf

\begin{lemma}\label{essurj}
For any $\unc$-module $\A$
the unique $\unc$-module morphism
from $\A$ to the strict $\unc$-module structure
on the symmetric Picard category $\A$ and which underlying 
map is the identity map $1_{\A}$ at $\A$ in $\BS$ 
-- given by Lemma \ref{uniqmor} -- is invertible in $\unc$-mod.
\end{lemma}
\pf
Let $\multh: \unc \rightarrow [\A,\A]$ denote the multiplication
of $\A$, and $t \times a$ stand for $\varphi(t)(a)$ for 
any objects $t$ of $\unc$ and $a$ of $\A$.
The 2-cell $\deltah$ given by \ref{uniqmor}
for the morphism from $\A$ to the strict $\unc$-module on $\A$ in $\BS$,
is in this case of the form
$\xymatrix{
\unc
\ar@/^20pt/[r]^{\vv}
\ar@{}[r]|{\Downarrow}
\ar@/_20pt/[r]_{\multh}
&
[\A,\A].
}$
Its inverse is therefore of the expected 
form to be part of a module morphism
from the strict $\unc$-module on $\A$ to the module 
$\A$ with multiplication $\multh$.\\
 
Recall that the coherence Axiom \ref{cohmdmor1b} for the pair 
$(1_{\A},\deltah)$
as a morphism from  
$\A$ with multiplication $\multh$ to the strict $\unc$-module on $\A$
is the commutativity of the diagram in $\A$
\begin{tag}\label{di1}
{\small
$$
\xymatrix{
a 
\ar[rr]^{\lml_a}
\ar@{-}[rdd]_{id}
& & 
\gen \times a
\\
\\
& \gen . a
\ar[ruu]_{\deltal_{\gen,a}}
}
$$ 
}
\end{tag} for any object $a$ where the bilinear $\deltal$ 
corresponds to $\deltah$,
whereas  Axiom \ref{cohmdmor2b} amounts to the commutation
of the diagram in $\A$
\begin{tag}\label{di2}
{\small
$$
\xymatrix{
t_1.(t_2.a)
\ar@{-}[rr]^{id}
\ar[d]_{t_1.\deltal_{t_2,a}}
& & 
(t_1.t_2).a
\ar[d]^{\deltal_{t_1.t_2,a}}
\\
t_1 . (t_2 \times a)
\ar[r]_{\deltal_{t_1, t_2 \times a}}
& 
t_1 \times (t_2 \times a)
\ar[r]_{\assml_{t_1,t_2,a}}
& 
(t_1.t_2) \times a
}
$$
}
\end{tag}
for all objects $t_1$, $t_2$ in $\unc$ 
and $a$ in $\A$.\\

Axiom \ref{cohmdmor1b} for the pair $(1_A,\deltah^{-1})$ amounts the 
commutation for any object $a$ in $\A$
of the diagram 
$$
\xymatrix{
a 
\ar[rdd]_{\lml_a}
\ar@{-}[rr]^{id}
& & 
\gen . a
\\
\\
& \gen \times a
\ar[ruu]_{\deltal^{-1}_{\gen,a}}
}
$$ 
which does commute since diagram \ref{di1} does.\\

Axiom \ref{cohmdmor2b} for the pair $(1_{\A},\delta^{-1})$ is 
the commutation for any objects $t_1,t_2$ in $\unc$ 
and $a$ in $\A$ of the diagram 
{\small
$$
\xymatrix{
t_1 \times (t_2 \times a)
\ar[rr]^{\assmh_{t_1,t_2,a}}
\ar[d]_{t_1 \times \deltal^{-1}_{t_2,a}}
& & 
(t_1.t_2) \times a
\ar[d]^{\deltal^{-1}_{t_1.t_2,a}}
\\
t_1 \times  (t_2 . a)
\ar[r]_{\deltal^{-1}_{t_1, t_2 . a}}
& 
t_1 . (t_2 . a)
\ar@{-}[r]_{id}
& 
(t_1.t_2) . a
}
$$
}
which is equivalent to the commutation of the external diagram
in the pasting
$$
\xymatrix{
& 
t_1 . (t_2 . a)
\ar@{-}[r]^{id}
\ar[ld]_{\deltal_{t_1,t_2.a}}
\ar[d]^{t_1.\deltal_{t_2,a}}
&
(t_1 . t_2) . a
\ar[dd]^{\deltal_{t_1.t_2,a}}
\\
t_1 \times (t_2.a)
\ar[rd]_{t_1 \times \deltal_{t_2,a}}
&
t_1 . (t_2 \times a)
\ar[d]^{\deltal_{t_1,t_2 \times a}}
&
\\
&
t_1 \times (t_2 \times a)
\ar[r]_-{ \assml_{t_1,t_2,a} }
&
(t_1 . t_2) \times a
}
$$
in which the left diagram commutes by naturality
of $\deltah_{t_1} : \vv(t_1) \rightarrow \multh(t_1): \unc \rightarrow \A$
and the right diagram is diagram \ref{di2}.
\epf


\end{section}

\begin{section}{Appendix}
This section contains various technical
developments.\\

\noindent{\bf Section \ref{prel}.}\\

\begin{lemma}\label{syPic}
In any symmetric Picard category
the diagram
$$
\xymatrix{
a^{\eqi} \otimes b^{\eqi}
\ar[r]^{!} 
\ar[d]_{\syc}
&
{(b \otimes a)}^{\eqi}
\ar[d]^{{\syc}^{\eqi}}
\\
b^{\eqi} \otimes a^{\eqi}
\ar[r]_{!}
&
{(a \otimes b)}^{\eqi}
}
$$ 
commutes for any objects $a$ and $b$.
\end{lemma}
\pf
By definition 
the canonical arrow 
$a^{\eqi} \otimes b^{\eqi} \rightarrow {(b \otimes a)}^{\eqi}$
is the only arrow $f$  making the diagram
$$
\xymatrix{
(a^{\eqi} \otimes b^{\eqi}) \otimes (b \otimes a)
\ar[rr]^{f \otimes 1}
\ar[d]_{\cong}
& 
& 
{(b \otimes a)}^{\eqi} \otimes (b \otimes a)
\ar[ddd]^{\jj}
\\
a^{\eqi} \otimes ((b^{\eqi} \otimes b) \otimes a)
\ar[d]_{1 \otimes (\jj \otimes 1)}
& 
\\
a^{\eqi} \otimes (\un \otimes a)
\ar[d]_{\cong}
\\
a^{\eqi} \otimes a
\ar[rr]_{\jj}
&
&
\un
}
$$ 
commute (see \cite{Lap83}-p.310). 
In the following pasting all diagrams commute,
$$\xymatrix{
(a^{\eqi} \otimes b^{\eqi})
\otimes (b \otimes a)
\ar[d]_{1 \otimes \syc}
\ar[r]^{\syc \otimes 1}
&
(b^{\eqi} \otimes a^{\eqi})
\otimes (b \otimes a)
\ar[r]^{! \otimes 1}
\ar[d]|{1 \otimes \syc}
&
{(a \otimes b)}^{\eqi}
\otimes (b \otimes a)
\ar[r]^{{\syc}^{\eqi} \otimes 1}
\ar[d]|{1 \otimes \syc}
&
{(b \otimes a)}^{\eqi} \otimes (b \otimes a)
\ar[d]^{\jj}
\\
(a^{\eqi} \otimes b^{\eqi})
\otimes (a \otimes b)
\ar[r]^{\syc \otimes 1}
&
(b^{\eqi} \otimes a^{\eqi})
\otimes (a \otimes b)
\ar[r]^{! \otimes 1}
\ar[d]|{\cong}
&
{(a \otimes b)}^{\eqi}
\otimes (a \otimes b)
\ar[r]^{\jj}
&
\un\\
& 
b^{\eqi} \otimes ((a^{\eqi} \otimes a) \otimes b)
\ar[r]_{1 \otimes (\jj \otimes 1)}
&
b^{\eqi} \otimes (\un \otimes b)
\ar[r]_{1 \otimes \rc}
&
b^{\eqi} \otimes b
\ar[u]^{\jj}
}$$
and by the coherence theorem for symmetric
monoidal categories the two left legs 
of the two diagrams above are equal.
\epf

\begin{tag}\label{pfinvmdal}
Proof of Lemma \ref{invmdal}.
\end{tag}
\pf
That the functor $\inv: \A \rightarrow \A$ is monoidal result from
points \ref{invmo1} 
that it is symmetric
amounts to Lemma \ref{syPic}.
\begin{tag}\label{invmo1}
In any symmetric Picard category the  diagram
$$\xymatrix{
a^{\eqi} \otimes (b^{\eqi} \otimes c^{\eqi})
\ar[r]^{\ac}
\ar[d]_{1 \otimes s}
&
(a^{\eqi} \otimes b^{\eqi}) \otimes c^{\eqi}
\ar[d]^{s  \otimes 1}
\\
a^{\eqi} \otimes (c^{\eqi} \otimes b^{\eqi})
\ar[d]_{1 \otimes !}
&
(b^{\eqi} \otimes a^{\eqi}) \otimes c^{\eqi}
\ar[d]^{! \otimes 1}
\\
a^{\eqi} \otimes {(b \otimes c)}^{\eqi}
\ar[d]_{s}
&
{(a \otimes b)}^{\eqi} \otimes c^{\eqi}
\ar[d]^{s}
\\
{(b \otimes c)}^{\eqi} \otimes a^{\eqi}
\ar[d]_{!}
&
c^{\eqi} \otimes {(a \otimes b)}^{\eqi}
\ar[d]^{!}
\\
{(a \otimes (b \otimes c))}^{\eqi}
\ar[r]_{{\ac}^{\eqi}}
&
{((a \otimes b) \otimes c)}^{\eqi}
}$$
commutes for any objects $a$, $b$ and $c$.
\end{tag}
\pf
By the naturality of $\syc$, the lemma is equivalent to 
the commutation of the external diagram in the pasting
$$\xymatrix{
a^{\eqi} \otimes (b^{\eqi} \otimes c^{\eqi})
\ar[r]^{\ac}
\ar[d]_{1 \otimes s}
&
(a^{\eqi} \otimes b^{\eqi}) \otimes c^{\eqi}
\ar[d]^{s  \otimes 1}
\\
a^{\eqi} \otimes (c^{\eqi} \otimes b^{\eqi})
\ar[d]_{s}
&
(b^{\eqi} \otimes a^{\eqi}) \otimes c^{\eqi}
\ar[d]^{s}
\\
(c^{\eqi} \otimes b^{\eqi}) \otimes a^{\eqi}
\ar[d]_{! \otimes 1}
\ar[r]^{\ac}
&
c^{\eqi} \otimes (b^{\eqi} \otimes a^{\eqi})
\ar[d]^{1 \otimes !}
\\
{(b \otimes c)}^{\eqi} \otimes a^{\eqi}
\ar[d]_{!}
&
c^{\eqi} \otimes {(a \otimes b)}^{\eqi}
\ar[d]^{!}
\\
{(a \otimes (b \otimes c))}^{\eqi}
\ar[r]_{{\ac}^{\eqi}}
&
{((a \otimes b) \otimes c)}^{\eqi}
}$$
where the top diagram commutes according
to the coherence for the symmetric monoidal 
structure and the bottom theorem
commutes according to the coherence theorem
for the group structure.
\epf


\begin{tag}\label{pfmonatgp}
Proof of Lemma \ref{monatgp}.
\end{tag}
\pf
The collection $\jj$ is natural by definition
of the functor $\inv$. That it is monoidal
amounts to the commutation of the diagram
$$\xymatrix{
\un 
\ar@{-}[ddd]_{\cong}
\ar[r]^-{\jj_{a \otimes b}}
&
{(a \otimes b)}^{\eqi} \otimes (a \otimes b)
\ar[d]^{! \otimes 1}
\\
& 
(b^{\eqi} \otimes a^{\eqi})
\otimes 
(a \otimes b)
\ar[d]^{\syc \otimes 1}
\\
& 
(a^{\eqi} \otimes b^{\eqi})
\otimes 
(a \otimes b)
\ar@{-}[d]^{\cong} 
\\
\un \otimes \un
\ar[r]_-{\jj_a \otimes \jj_b}
&
(a^{\eqi} \otimes a)
\otimes 
(b^{\eqi} \otimes b)
}
$$
for any objects $a$ and $b$ of $\A$, which holds
by definition of the arrow 
$!: \un: {(a \otimes b)}^{\eqi} \rightarrow b^{\eqi} \otimes a^{\eqi}$
see \cite{Lap83} p.310.
\epf

\begin{tag}\label{defisofuntag}
Definition of the isomorphism \ref{isofuntag}.
\end{tag}
The isomorphism \ref{isofuntag} is defined pointwise in $a$
as the only arrow in $\B$ making 
the diagram 
$$
\xymatrix{
\un
\ar[d]_{F^0}
\ar[r]^-{\jj_{Fa}} 
&
{(Fa)}^{\eqi} \otimes Fa
\ar@{-}[d]^{\cong \otimes 1}
\\
F\un
\ar[rd]_{F(\jj_a)}
&
F(a^{\eqi}) \otimes Fa
\ar[d]^{F^2_{a^{\eqi},a}} 
\\
& 
F(a^{\eqi} \otimes a)
}
$$
commute.\\

\begin{lemma}\label{pcanco}
Given arrows 
$\xymatrix{
\A 
\ar[r]^F
&
\B 
\ar[r]^G
&
C
}
$ in $\BS$
the diagram in $\C$
$$
\xymatrix{
{(FG(a))}^{\eqi} 
\ar@{-}[r]^-{\cong}
\ar@{-}@/_30pt/[rr]_-{\cong}
&
F( {G(a)}^{\eqi} )
\ar@{-}[r]^-{F(\cong)}
&
FG(a^{\eqi})
}
$$
where all the $\cong$ are of type \ref{isofuntag},
is commutative. 
\end{lemma}
\pf
Consider the pasting of commutative diagrams
below where all the $\cong$ are of type \ref{isofuntag} 
$$
\xymatrix{
\un
\ar@/_40pt/[ddd]_{{(FG)}^0}
\ar[rr]^-{\jj_{FGa}}
\ar[dd]^{F^0}
& & 
{(FGa)}^{\eqi}  \otimes FGa
\ar[ld]_{\cong \otimes 1}
\ar[rd]^{\cong \otimes 1}
\\
& 
F{(Ga)}^{\eqi} \otimes FGa
\ar[rr]_{F(\cong) \otimes 1}
\ar[d]_{F^2_{{(Ga)}^\eqi,Ga}}
& & 
FG(a^{\eqi}) \otimes FGa
\ar[d]_{F^2_{G(a^{\eqi}), Ga}}
\ar@/^40pt/[dd]^{{(FG)}^2_{a^{\eqi},a}}
\\
F \un
\ar[r]_-{F\jj_{Ga}}
\ar[d]^{F(G^0)}
& 
F({(Ga)}^{\eqi}   \otimes  Ga  )
\ar[rr]_{F(\cong \otimes 1)}
& & 
F(G(a^{\eqi}) \otimes Ga)
\ar[d]_{F(G^2_{ a^{\eqi},a }) }
\\
FG\un
\ar[rrr]_{FG\jj_a}
& & & 
FG(a^{\eqi} \otimes a)
}
$$
\epf

\begin{lemma}\label{sigmap}
For any monoidal transformation 
$\sigma:F \rightarrow G: \A \rightarrow \B$
where $\A$ and $\B$ are objects of $\BS$, for any object 
$a$ of $\A$, the diagram in $\B$
$$
\xymatrix{
{(Fa)}^{\eqi}
\ar[d]_{\cong}
& 
{(Ga)}^{\eqi}
\ar[l]_-{{(\sigma_a)}^{\eqi}}
\ar[d]^{\cong}
\\
F(a^{\eqi})
\ar[r]_-{\sigma_{a^{\eqi}}}
& 
G(a^{\eqi})
}
$$  
where the $\cong$ denote isomorphisms of type \ref{isofuntag},
is commutative.
\end{lemma}
\pf
Consider the pasting of diagrams in $\B$ 
$$\xymatrix{
\un 
\ar[rr]^-{\jj_{Fa}}
\ar[dd]_{F^0}
\ar@{-}[rrrd]_(.4){id}
&
&
{ (Fa) }^{\eqi} \otimes Fa
\ar@{.>}[d]^{\cong \otimes 1}
\ar@{.>}[rrrd]^{ {[{(\sigma_a)}^\eqi]}^{-1} \otimes \sigma_a }
\\
&
& 
F( a^{\eqi} ) \otimes Fa
\ar[d]_{ F^2_{a^{\eqi},a} }
\ar@{.>}[rrrd]^(.6){ \sigma_{a^{\eqi}} \otimes \sigma_a}
& 
\un
\ar[rr]^{\jj_{Ga}}
\ar[dd]_{G^0}
&
&
{(Ga)}^{\eqi} \otimes Ga
\ar@{.>}[d]^{\cong \otimes 1}
\\
F\un
\ar[rr]^-{F\jj_a}
\ar[rrrd]_{\sigma_{\un}}
&
&
F(a^{\eqi} \otimes a)
\ar[rrrd]^(.6){\sigma_{a^{\eqi} \otimes a}}
&
& 
& 
G(a^{\eqi}) \otimes Ga
\ar[d]^{G^2_{a^{\eqi},a}}
\\
& 
&
& 
G\un
\ar[rr]_{G\jj_a}
&
&
G(a^{\eqi} \otimes a)
}$$
where the $\cong$ denote isomorphisms of type
\ref{isofuntag}. All diagrams above commute
apart from the one consisting of the four
dotted arrows. Since all arrows are invertible
this last diagram commutes. The 
result follows since tensoring with $Fa$ is an 
equivalence $\B \rightarrow \B$.
\epf

\noindent{\bf A few computational lemmas for $\BS$.}\\
We present here a couple of results not stated in \cite{Sch08}.

\begin{lemma}\label{calc6}
The 2-cells in $\BS$
$$\xymatrix{
\A 
\ar[r]^-{L'}
& 
\unc \otimes \A
\ar[r]^-{U \otimes 1}
&
\B \otimes \A
\ar@{=>}[r]^-{\tau}
&
\C
}$$
and
$$\xymatrix{
\A 
\ar@{=>}[r]^-{{\Res(\tau)}^*} 
&
[\B, \C]
\ar[r]^-{[U,1]}
&
[\unc,\C]
\ar[r]^-{\ev_{\gen}}
&
\C
}$$ 
are equal for any 2-cell $\tau: \B \otimes \A \rightarrow \C$
and $U: \unc \rightarrow \B$.
\end{lemma}
\pf
The 2-cell $\tau \circ (U \otimes 1) \circ L'$ above 
rewrites successively
\begin{tabbing}
\=1. $\xymatrix@C=3pc{
\A 
\ar[r]^-{\eta^*}
&
[\unc, \unc \otimes \A]
\ar[r]^-{\ev_{\gen}}
&
\unc \otimes \A
\ar[r]^-{U \otimes 1}
&
\B \otimes \A
\ar@{=>}[r]^-{\tau}
&
\C
}
$\\
\>2. $\xymatrix@C=3pc{
\A 
\ar[r]^-{\eta^*}
&
[\unc, \unc \otimes \A]
\ar[r]^-{[1, U \otimes 1]}
&
[\unc, \B \otimes \A]
\ar@{=>}[r]^-{[1, \tau]}
&
[\unc, \B \otimes \A]
\ar@{=>}[r]^-{\ev_{\gen}}
&
\C
}
$\\
\>3. $\xymatrix@C=3pc{
\A 
\ar[r]^-{\eta^*}
&
[\B, \B \otimes \A]
\ar[r]^-{[U, 1]}
&
[\unc, \B \otimes \A]
\ar@{=>}[r]^-{[1,\tau]}
&
[\unc,\C]
\ar[r]^{\ev_{\gen}}
&
\C
}
$\\
\>4. $\xymatrix@C=3pc{
\A 
\ar[r]^-{\eta^*}
&
[\B, \B \otimes \A]
\ar@{=>}[r]^-{[1, \tau]}
&
[\B, \C ]
\ar[r]^-{[U,1]}
&
[\unc,\C]
\ar[r]^{\ev_{\gen}}
&
\C
}
$\\
\>5. $\xymatrix@C=3pc{
\A 
\ar@{=>}[r]^-{{(\Res(\tau))}^*}
&
[\B, \C ]
\ar[r]^-{[U,1]}
&
[\unc,\C]
\ar[r]^{\ev_{\gen}}
&
\C.
}$
\end{tabbing}
\epf

\begin{lemma}\label{calc51}
The 2-cells in $\BS$
$$\xymatrix{
\A 
\ar[r]^-{R'}
& 
\A \otimes \unc
\ar[r]^-{1 \otimes U}
&
\A \otimes \B
\ar@{=>}[r]^-{\tau}
&
\C
}$$
and
$$\xymatrix{
\A 
\ar@{=>}[r]^-{\Res(\tau)}
& 
[\B, \C]
\ar[r]^-{[U,1]}
&
[\unc,\C]
\ar[r]^-{\ev_{\gen}}
&
\C
}$$ 
are equal for any 2-cell $\tau: \A \otimes \B \rightarrow \C$
and any arrow $U: \unc \rightarrow \B$.
\end{lemma}
\pf
The 2-cell $\tau \circ (1 \otimes U) \circ R'$ above 
rewrites successively
\begin{tabbing}
\=1. $\xymatrix@C=3pc{
\A 
\ar[r]^-{\eta}
&
[\unc, \A \otimes \unc]
\ar[r]^-{\ev_{\gen}}
&
\A \otimes \unc
\ar[r]^-{1 \otimes U}
&
\A \otimes \B
\ar@{=>}[r]^-{\tau}
&
\C
}
$
\\
\>2.
$\xymatrix@C=3pc{
\A 
\ar[r]^-{\eta}
&
[\unc, \A \otimes \unc]
\ar[r]^-{[1, 1 \otimes U]}
&
[\unc,\A \otimes \B]
\ar@{=>}[r]^-{[1,\tau]}
&
[\unc,\C]
\ar[r]^-{\ev_{\gen}}
&
\C
}
$
\\
\>3.
$
\xymatrix@C=3pc{
\A 
\ar[r]^-{\eta}
&
[\B, \A \otimes \B]
\ar[r]^-{[U,1]}
&
[\unc,\A \otimes \B]
\ar@{=>}[r]^-{[1,\tau]}
&
[\unc,\C]
\ar[r]^-{\ev_{\gen}}
&
\C
}
$
\\
\>4.
$
\xymatrix@C=3pc{
\A 
\ar[r]^-{\eta}
&
[\B, \A \otimes \B]
\ar@{=>}[r]^-{[1,\tau]}
&
[\A,\C]
\ar[r]^-{[U,1]}
&
[\unc,\C]
\ar[r]^-{\ev_{\gen}}
&
\C
}
$
\\
\>5.
$
\xymatrix@C=3pc{
\A 
\ar@{=>}[r]^-{\Res(\tau)}
&
[\B,\C]
\ar[r]^-{[U,1]}
&
[\unc,\C]
\ar[r]^-{\ev_{\gen}}
&
\C
}
$
\end{tabbing}
where in the above derivation arrows 1. and 2. are
equal according to Corollary \cite{Sch08}-\deuxnateva.
\epf

\begin{remark}\label{evgen2cel}
Any 2-cell $\sigma: F \rightarrow G : \unc \rightarrow \A$
with $F$ {\em strict} in $\BS$ is fully determined by its component in $\gen$
since $F$ being strict the component at $F$ of counit of the adjunction 
$\epsilon_F : \vs \circ \ev_{\gen}(F) \rightarrow F$ is an identity
and by naturality of $\epsilon$ one has the commutation of 
$$
\xymatrix{
\vs (F(\gen)) 
\ar[r]^-{\vs ( \sigma_{\gen} ) } 
\ar@{-}[d]_{id}
&
 \vs (G(\gen))
\ar[d]^{\epsilon^G}
\\
F
\ar[r]_-{\sigma}
&
G
}
$$
in $\BS$.
\end{remark}

\begin{lemma}\label{evstrict}
For any {\em strict} arrow $F: \unc \rightarrow \A$
the diagram in $\BS$ 
$$
\xymatrix{
\unc 
\ar[r]^{\vv}
\ar[d]_F
&
[\A,\A]
\ar[d]^{[F,1]}
\\
\A
&
[\unc,\A]
\ar[l]^{\ev_{\gen}}
}
$$
commutes.
\end{lemma}
\pf
The arrow 
$$\xymatrix{\unc \ar[r]^{\vv} & 
[\A,\A] 
\ar[r]^{[F,1]} 
& 
[\unc,\A]
}$$
has dual 
$$\xymatrix{
\unc 
\ar[r]^F
&
\A
\ar[r]^{\vs}
&
[\unc,\A].
}$$
Since the arrow $F$ is strict it is equal
to $\vs \circ \ev_{\gen}(F)$ and the result 
follows then from Lemma \cite{Sch08}-\evdual.
\epf

\begin{remark}\label{pastevsev}
Since the units of the adjunctions $\vs \dashv \ev_{\gen}$
are identities in $\BS$, one has a bijection 
for any $\A$ and $\B$ between sets of 2-cells
in $\BS$ of the following kind
\begin{tag}\label{2ctype1}
$$\xymatrix{
&
\B
\ar[rd]^{\vs}
\ar@{=>}[d]
\\
\A
\ar[ru]^F
\ar[rr]_G
& & 
[\unc, \B]
}$$
\end{tag}
and 
\begin{tag}\label{2ctype2}
$$
\xymatrix{
\A
\ar[rr]^F
\ar[rd]_G
&
\ar@{=>}[d]
&
\B 
\\
&
[\unc,\A]
\ar[ru]_{\ev_{\gen}}
}
$$
\end{tag}
sending any 2-cell $\Xi$ of the type \ref{2ctype1}  
to 
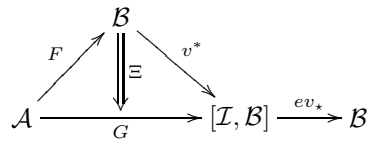
$$
\xymatrix{
&
\B
\ar[rd]^-{\vs}
\ar@{=>}[d]^-{\Xi}
\\
\A
\ar[ru]^-F
\ar[rr]_-G
& & 
[\unc, \B]
\ar[r]^-{\ev_{\gen}}
&
\B
}
$$
with inverse sending any 2-cell $\Xi'$ of type \ref{2ctype2}   
to
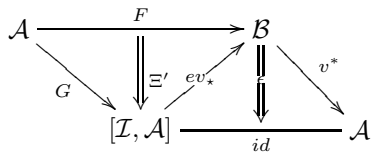
$$
\xymatrix{
\A
\ar[rr]^-F
\ar[rd]_-G
&
\ar@{=>}[d]^{\Xi'}
&
\B
\ar[rd]^{\vs}
\ar@{=>}[d]|{\epsilon}
\\
&
[\unc,\A]
\ar[ru]|{\ev_{\gen}}
\ar@{-}[rr]_-{id}
&
&
\A
}
$$
\end{remark}

We shall describe for any arrow $\varphi: \A \otimes \M \rightarrow \M$
and $U: \unc \rightarrow \A$
of $\BS$ some bijections between sets of 2-cells of the 
following kind:
\begin{tag}\label{2c1}
{\tiny
$$
\xymatrix{
\unc \M
\ar[rr]^-{U \otimes 1} 
\ar[rrdd]_{L}
&  
& 
\A \M
\ar[dd]^{\varphi}
\\
&
\ar@{=>}[ru]
&
\\
&
&
\M
}
$$ 
}
\end{tag}
\begin{tag}\label{2c2}
{\tiny
$$
\xymatrix{
\unc 
\ar[rr]^-{U} 
\ar[rrdd]_{\vv}
&  
& 
\A
\ar[dd]^{\Res(\varphi)}
\\
&
\ar@{=>}[ru]
&
\\
&
&
[\M,\M]
}
$$ 
}
\end{tag}
\begin{tag}
\label{2c3}
{\tiny
$$
\xymatrix{
\M 
\ar[rr]^-{{(\Res(\varphi))}^*} 
\ar[rrdd]_{\vs}
&  
& 
[\A,\M]
\ar[dd]^{[U,1]}
\\
&
\ar@{=>}[ru]
&
\\
&
&
[\unc,\M]
}
$$ 
}
\end{tag}
\begin{tag}
\label{2c4}
{\tiny
$$
\xymatrix{
\M
\ar@{-}[dd]_{id}
\ar[rr]^-{ {(\Res(\varphi))}^* }
& 
& 
[\A,\M]
\ar[dd]^{[U,1]}
\\
\ar@{=>}[rr] & &  
\\
\M
& 
& 
[\unc,\M]
\ar[ll]^{\ev_{\gen}}
}
$$
}
\end{tag}
and
\begin{tag}
\label{2c5}
{\tiny
$$
\xymatrix{
\M
\ar[rr]^-{L'}
\ar@{-}[dd]_{id}
& 
& 
\unc \M
\ar[dd]^{U \otimes 1}
\\
\ar@{=>}[rr] & & 
\\
\M
& & 
\A \M
\ar[ll]^-{\varphi}
}
$$ 
}
\end{tag}
as follows.\\
- Since $L$ is the image by $\Ext$ 
of  $\vv: \unc \rightarrow [\M,\M]$
and the image by $\Res$ of 
$$\xymatrix{
\unc \M 
\ar[r]^-{U \otimes 1}
&
\A \M
\ar[r]^-{\varphi}
&
\M
}$$
is 
$$
\xymatrix{
\unc 
\ar[r]^-{U}
&
\A
\ar[r]^-{\Res(\varphi)}
&
[\M,\M]
}
$$
the maps $\Res$/$\Ext$ define the bijection (and its inverse)
between sets of 2-cells $\ref{2c1}$ and $\ref{2c2}$.\\
- 2-cells \ref{2c2} and \ref{2c3} correspond by duality.\\ 
- The bijection between sets of 2-cells $\ref{2c3}$ and $\ref{2c4}$
sends any 2-cell $\Xi: \vs \rightarrow [j,1] \circ {(\Res(\varphi))}^*$
to
{\small
$$
\xymatrix{
\M
\ar[rr]^{ { (\Res(\varphi)) }^* } 
\ar[rrdd]_{\vs}
&  
& 
[\A,\M]
\ar[dd]^{[U,1]}
\\
&
\ar@{=>}[ru]^{\Xi}
&
\\
&
&
[\unc,\M]
\ar[rr]_{\ev_{\gen}}
& 
&
\M
}
$$
}
which is as expected a 2-cell
$$id \rightarrow \ev_{\gen} \circ [U,1] \circ {(\Res(\varphi))}^*$$ 
according to the adjunction \ref{unimonadj}. 
Its inverse 
sends any 2-cell 
$\Xi: 1 \rightarrow \ev_{\gen}\circ [U,1] \circ {(\Res(\varphi))}^*$ 
to the pasting
{\small
$$
\xymatrix{
\M
\ar[rr]^{id}
\ar[dd]_{{(\Res(\varphi))}^*}
& 
\ar@{=>}[dd]^{\Xi}
& 
\M
\ar[rr]^{\vs}
\ar@{=>}[rd]^{\epsilon}
& & 
[\unc,\M]
\\
& & & & 
\\
[\A,\M]
\ar[rr]_{[U,1]}
& &
[\unc,\M]
\ar[rruu]_{id}
\ar[uu]_{\ev_{\gen}}
}
$$
} 
- Eventually 2-cells $\ref{2c4}$ and $\ref{2c5}$ are the same since 
their codomains arrows respectively
$\ev_{\gen} \circ [U,1] \circ {\Res(\varphi)}^*$
and $\varphi \circ (U \otimes 1) \circ L'$
are equal according to Lemma \ref{calc6}.\\

\begin{lemma}\label{vilem}
The above bijection between sets of 2-cells \ref{2c1}/\ref{2c5} sends
any $\Xi: L \rightarrow \varphi \circ (U \otimes 1): 
\unc \M \rightarrow \M$
to the pasting
{\tiny
$$\xymatrix{
\unc \M
\ar[rr]^-{U \otimes 1} 
\ar[rrdd]_{L}
&  
& 
\A \M
\ar[dd]^{\varphi}
\\
&
\ar@{=>}[ru]^{\Xi}
&
\\
\M 
\ar[uu]^{L'}
\ar[rr]_{id}
\ar@{}[ru]|{=}
&
&
\M
}$$
}
\end{lemma}
\pf
According to Lemma \ref{calc6}, the above
2-cell 
$\xymatrix{
\A   
\ar[r]^-{L'} 
&
\unc \A 
\ar@{=>}[r]^-{\Xi} 
&
\A
}$
is
$$ 
\xymatrix{
\A
\ar@{=>}[r]^-{ {\Res(\Xi)}^* }
&
[\unc,\A]
\ar[r]^-{\ev_{\gen}}
&
\A.
}$$
Consider then the image of the above 2-cell by the bijection
\ref{2c4} $\rightarrow$ \ref{2c3},
it is
$$
\xymatrix{
&
[\A,\A]
\ar@{=>}[d]|{{(\Res(\Xi)}^*}
\ar[rd]^{[U,1]}
& & 
\A
\ar[rd]^{\vs}
\ar@{=>}[d]^{\epsilon}
&
\\
\A
\ar[rr]_-{\vs}
\ar[ru]^{\Res(\varphi)}
& & 
[\unc,\A]
\ar@{-}[rr]_-{id}
\ar[ru]^{\ev_{\gen}}
& & 
[\unc,\A]
}
$$
which is just ${(\Res(\bar{\lp}))}^*: \A \rightarrow [\unc,\A]$
since 
$\xymatrix{
\A \ar[r]^{\vs} 
&
[\unc,\A] \ar@{=>}[r]^{\epsilon} 
&
\A}$ is an identity 2-cell,
which dual has image by $\Ext$ the 2-cell
$\bar{\lp}$.
\epf

Easy computation also gives the following.
\begin{remark}\label{calc1}
For any arrows $F: \A \rightarrow [\B,\X]$
and $G: \X \rightarrow [\C,\D]$ of $\BS$
the arrow 
$$\xymatrix@C=3pc{
(\A \otimes \B) \otimes \C
\ar[r]^-{\Ext(F) \otimes 1}
&
\X \otimes \C
\ar[r]^-{\Ext(G)}
&
\D
}$$
has image by $\Res$
$$
\xymatrix{
\A \otimes \B 
\ar[r]^-{\Ext(F)}
&
\X 
\ar[r]^-G
&
[\C,\D]
}
$$
and has image by $\Res \circ \Res$ the arrow
$$ 
\xymatrix{
\A
\ar[r]^-F
&
[\B,\X]
\ar[r]^{[1,G]}
&
[\B,[\C,\D]].
}
$$
\end{remark}

\noindent{\bf Sections \ref{BSCat} and \ref{BSCattens}.}\\

\begin{tag}Proof of Equality (I) in second pasting 
of Axiom \ref{cohunassh}.
\end{tag} \label{ptequal1}
\pf To check the equality $(I)$ of arrows, consider 
the derivation of equal composite arrows in $\BS$ for
any arrows
$d: \A \rightarrow [\A,\A]$ and 
$j: \unc \rightarrow \A$.
\begin{tabbing}
1. $\xymatrix{
[\A,\B]
\ar[r]^-{[\A,-]}
&
[[\A,\A],[\A,\B]]
\ar[r]^-{[d,1]}
&
[\A,[\A,\B]]
\ar[r]^-{[j,1]}
&
[\unc,[\A,\B]]
\ar[r]^-{\ev_{\gen}}
&
[\A,\B]
}$\\
2. $\xymatrix{
[\A,\B]
\ar[r]^-{[\A,-]}
&
[[\A,\A],[\A,\B]]
\ar[r]^-{[d,1]}
&
[\A,[\A,\B]]
\ar[r]^-{\ev_{j(\gen)}}
&
[\A,\B]
}$\\
3. $\xymatrix{
[\A,\B]
\ar[r]^-{[\A,-]}
&
[[\A,\A],[\A,\B]]
\ar[r]^-{[d^*,1]}
&
[\A,[\A,\B]]
\ar[r]^-{\Dual}
&
[\A,[\A,\B]]
\ar[r]^-{\ev_{j(\gen)}}
&
[\A,\B]
}$\\
4. $\xymatrix{
[\A,\B]
\ar[r]^-{[\A,-]}
&
[[\A,\A],[\A,\B]]
\ar[r]^-{[ d^*, 1  ] }
&
[\A,[\A,\B]]
\ar[r]^-{[1,\ev_{j(\gen)}]}
&
[\A,\B]
}$\\
5. $\xymatrix{
[\A,\B]
\ar[r]^-{[\A,-]}
&
[[\A,\A],[\A,\B]]
\ar[r]^-{[1, \ev_{j(\gen)}]}
&
[\A,[\A,\B]]
\ar[r]^-{[d^*,1]}
&
[\A,\B]
}$\\
6. $\xymatrix{
[\A,\B]
\ar[r]^-{[\ev_{j(\gen)},1]}
&
[[\A,\A],\B]
\ar[r]^-{[d^*,1]}
&
[\A,\B]
}$\\
7. $\xymatrix{
[\A,\B]
\ar[r]^-{[\ev_{\gen},1]}
&
[[\unc,\A],\B]
\ar[r]^-{[ [j,1],1  ] }
&
[[\A,\A],\B]
\ar[r]^-{[d^*,1]}
&
[\A,\B]
}$
\end{tabbing}
where the equalities between arrows stand for the 
following reasons:\\
- 1. and 2. : by the naturality in $\A$ of the collection of 
arrows $q_{\A}: \A \rightarrow [[\A,\C],\C]$\\
- 2. and 3. : Lemma \cite{Sch08}-\tbctrois;\\ 
- 3. and 4. : Lemma \cite{Sch08}-\evdual;\\
- 5. and 6. : Lemma \cite{Sch08}-\homevtrois;\\
- 6. and 7. : by the naturality of the collections of arrows $q$;\\
\epf

\begin{tag}\label{defTCcomp1fyz}
Definition of the 2-cell $\TCcompu_{f:x' \rightarrow x,y,z}$.
\end{tag}
The domain of this 2-cell rewrites
{\tiny
\begin{tabbing}
\=1. \=$\xymatrix{
\A_{y,z}
\ar[r]^-{\A(x',-)}
&
[\A_{x',y}, \A_{x',z}]
\ar[r]^-{[\A_{x',x},-]}
&
[[\A_{x',x},\A_{x',y}],[\A_{x',x},\A_{x',z}]]
\ar[r]^-{[\A(x',-),1]}
&
[\A_{x,y},[\A_{x',x}, \A_{x',z}]]
\ar[r]^-{[1,[f,1]]}
&
[\A_{x,y}, [\unc, \A_{x',z}]]
\;...
}$\\
\> \>$\xymatrix{... 
\ar[r]^-{[1,\ev_{\gen}]}
&
[\A_{x,y},\A_{x',z}]
}$\\
\>2. $\xymatrix{
\A_{y,z}
\ar[r]^-{\A(x',-)}
&
[\A_{x',y}, \A_{x',z}]
\ar[r]^-{[\A_{x',x},-]}
&
[[\A_{x',x},\A_{x',y}],[\A_{x',x},\A_{x',z}]]
\ar[r]^-{[1,[f,1]]}
&
[[\A_{x',x},\A_{x',y}],[\unc,\A_{x',z}]]
\ar[r]^-{[1,\ev_{\gen}]}
&
...}$\\
\> \>
$\xymatrix{...\;\;
[[\A_{x',x},\A_{x',y}],\A_{x',z}]
\ar[r]^-{[\A(x',-),1]}
&
[\A_{x,y},\A_{x',z}]
}$\\
\>3. $\xymatrix{
\A_{y,z}
\ar[r]^-{\A(x',-)}
&
[\A_{x',y}, \A_{x',z}]
\ar[r]^-{[\unc,-]}
&
[[\unc,\A_{x',y}],[\unc,\A_{x',z}]]
\ar[r]^-{[[f,1],1]}
&
[[\A_{x',x},\A_{x',y}],[\unc,\A_{x',z}]]
\ar[r]^-{[1,\ev_{\gen}]}
& ...}$\\
\> \>$\xymatrix{...\;\;
[[\A_{x',x},\A_{x',y}],\A_{x',z}]
\ar[r]^-{[\A(x',-),1]}
&
[\A_{x,y},\A_{x',z}]
}$\\
\>4. $\xymatrix{
\A_{y,z}
\ar[r]^-{\A(x',-)}
&
[\A_{x',y}, \A_{x',z}]
\ar[r]^-{[\unc,-]}
&
[[\unc,\A_{x',y}],[\unc,\A_{x',z}]]
\ar[r]^-{[[1,\ev_{\gen}]}
&
[[\unc,\A_{x',y}],\A_{x',z}]
\ar[r]^-{[[f,1],1]}
&
[[\A_{x',x},\A_{x',y}],\A_{x',z}]
\;...}$\\
\>\>$\xymatrix{...
\ar[r]^-{[\A(x',-),1]}
&
[\A_{x,y},\A_{x',z}]
}$\\
\>5. $\xymatrix{
\A_{y,z}
\ar[r]^-{\A(x',-)}
&
[\A_{x',y}, \A_{x',z}]
\ar[r]^-{[\ev_{\gen},1]}
&
[[\unc,\A_{x',y}],\A_{x',z}]
\ar[r]^-{[[f,1],1]}
&
[[\A_{x',x},\A_{x',y}],\A_{x',z}]
\ar[r]^-{[\A(x',-),1]}
&
[\A_{x,y},\A_{x',z}]
}$\\
\>6. $\xymatrix{
\A_{y,z}
\ar[r]^-{\A(x',-)}
&
[\A_{x',y}, \A_{x',z}]
\ar[r]^-{[\A(f,1),1]}
&
[\A_{x,y},\A_{x',z}]
.}$\\
\end{tabbing}
}
In the previous derivation arrows 2. and 3. are equal
according to Lemma \cite{Sch08}-\impronatPosttrois and arrows 4. 
and 5. are equal according to Lemma \cite{Sch08}-\homevtrois. \\

The codomain of $\TCcompu_{f,y,z}$ rewrites
\begin{tabbing}
\=$\xymatrix@C=3pc{
\A_{y,z}
\ar[r]^-{\A(x,-)}
&
[\A_{x,y}, \A_{x,z}]
\ar[r]^-{[1,\A(x',-)]}
&
[\A_{x,y}, [\A_{x',x}, \A_{x',z}]]
\ar[r]^-{[1,[f,1]]}
&
[\A_{x,y},[\unc,\A_{x',z}]]
\ar[r]^-{[1,\ev_{\gen}]}
&
[\A_{x,y}, \A_{x',z}]
}$\\
\>$\xymatrix{
\A_{x,y}
\ar[r]^-{\A(x,-)}
&
[\A_{x,y}, \A_{x,z}]
\ar[r]^-{[1,\A(f,1)]}
&
[\A_{x,y}, \A_{x',z}]
}
$
\end{tabbing}

\begin{remark}\label{remtoto2}
For any arrows $f: \A \rightarrow \B$ and 
$\tilde{f}: \unc \rightarrow [\A,\B]$ such that 
$\ev_{\gen}(\tilde{f}) = f$ and any object $D$ in $\BS$,
all the diagrams
in the pasting below commute
$$
\xymatrix{
[\A,\B]
\ar[r]^-{[\D,-]}
\ar[d]_{[f,\B]}
&
[[\D,\A],[\D,\B]]
\ar[d]^{[\tilde{f},1]}
\ar[ld]|{\ev_f}
\\
[\D,\B]
&
[\unc,[\D,\B]]
\ar[l]^-{\ev_{\gen}}
.}
$$
The top left one commutes according to Corollary
\cite{Sch08}-\evhomun
and the bottom right one does according to the 
2-naturality of the collection of arrows $q$.
\end{remark}

\begin{tag}\label{defTCcomp2xgz}
Definition of the 2-cell $\TCcompd_{x,g:y' \rightarrow y,z}$.
\end{tag}
Observe that the image by 
$\ev_{\gen}: \BS(\unc, [\A_{x,y'},\A_{x,y}])
\rightarrow \BS(\A_{x,y'},\A_{x,y})$ 
of 
$$\xymatrix{
\unc 
\ar[r]^-{\tilde{g}} 
&
\A_{y',y}
\ar[r]^-{\A(x,-)}
&
[\A_{x,y'}, \A_{x,y}]
}$$
is equal to
$\A(1,g): \A_{x,y'} \rightarrow \A_{x,y}$.
Therefore according to Remark \ref{remtoto2} the domain
of this 2-cell which is 
{\small
\begin{tabbing}
\=$\xymatrix{
\A_{y,z}
\ar[r]^-{\A(x,-)}
&
[\A_{x,y},\A_{x,z}]
\ar[r]^-{[\A_{x,y'},-]}
&
[[\A_{x,y'},\A_{x,y}],[\A_{x,y'}, \A_{x,z}]]
\ar[r]^-{[\A(x,-),1]}
&
[\A_{y,y},[\A_{x,y'},\A_{x,z}]]
\ar[r]^-{[g,1]}
&
[\unc,[\A_{x,y'},\A_{x,z}]]
\;...}$\\
\>$\xymatrix{ ...\;
\ar[r]^-{\ev_{\gen}}
&
[\A_{x,y'},\A_{x,z}]
}$
\end{tabbing}
}
is equal to 
$$
\xymatrix@C=3pc{
\A_{y,z}
\ar[r]^-{\A(x,-)}
&
[\A_{x,y},\A_{x,z}]
\ar[r]^-{[\A(1,g),1]}
&
[\A_{x,y'},\A_{x,z}]
}
$$
The codomain of the 2-cell 
$\TCcompd_{x,g:y' \rightarrow y,z}$
rewrites successively
\begin{tabbing}
\=1.$\xymatrix{
\A_{y,z}
\ar[r]^-{\A(y'-)}
&
[\A_{y',y},\A_{y',z}]
\ar[r]^-{[1,\A(x,-)]}
&
[\A_{y',y},[\A_{x,y'},\A_{x,z}]]
\ar[r]^-{[g,1]}
&
[\unc, [\A_{x,y'}, \A_{x,z}]]
\ar[r]^-{\ev_{\gen}}
&
[\A_{x,y'}, \A_{x,z}]
}$\\
\>2. $\xymatrix@C=3pc{
\A_{y,z}
\ar[r]^-{\A(y'-)}
&
[\A_{y',y},\A_{y',z}]
\ar[r]^-{[g,1]}
&
[\unc, \A_{y',z}]
\ar[r]^-{[1,\A(x,-)]}
&
[\unc,[\A_{x,y'},\A_{x,z}]]
\ar[r]^-{\ev_{\gen}}
&
[\A_{x,y'}, \A_{x,z}]
}$\\
\>3. $\xymatrix{
\A_{y,z}
\ar[r]^-{\A(y'-)}
&
[\A_{y',y},\A_{y',z}]
\ar[r]^-{[g,1]}
&
[\unc, \A_{y',z}]
\ar[r]^-{\ev_{\gen}}
&
\A_{y',z}
\ar[r]^-{\A(x,-)}
&
[\A_{x,y'}, \A_{x,z}]
}$\\
\>4. $\xymatrix{
\A_{y,z}
\ar[r]^-{\A(g,1)}
&
\A_{y',z}
\ar[r]^-{\A(x,-)}
&
[\A_{x,y'}, \A_{x,z}]
}$\\
\end{tabbing}
where in the above derivation arrows 2. and 3. are equal 
due to Lemma \cite{Sch08}-\deuxnateva.\\

\begin{tag}\label{defTCcomp3xyh}
Definition of the 2-cell $\TCcompt_{x, y, h:z \rightarrow z'}$.
\end{tag}
The 2-cell
$$
\xymatrix{
\unc 
\ar[r]^-h
&
\A_{z,z'}
\ar[r]^-{\A(x,-)}
&
[\A_{x,z},\A_{x,z'}]
\ar[r]^-{[\A_{x,y},-]}
&
[[\A_{x,y},\A_{x,z}],
[\A_{x,y}, \A_{x,z'}]]
\ar[r]^-{[\A(x,-),1]}
&
[\A_{y,z},[\A_{x,y},\A_{x,z'}]]
}$$
has image by $\ev_{\gen}$
$$
\xymatrix{
\A_{y,z}
\ar[r]^-{\A(x,-)}
&
[\A_{x,y},\A_{x,z}]
\ar[r]^-{[1,\A(1,h)]}
&
[\A_{x,y},\A_{x,z'}]
}
$$
which is the domain of $\TCcompt_{x, y, h:z \rightarrow z'}$.\\
The 2-cell
$$
\xymatrix@C=3pc{
\unc 
\ar[r]^h
&
\A_{z,z'}
\ar[r]^-{\A(y,-)}
&
[\A_{y,z},\A_{y,z'}]
\ar[r]^-{[1,\A(x,-)]}
&
[\A_{y,z},[\A_{x,y}, \A_{x,z'}]]
}
$$
has image by $\ev_{\gen}$
$$
\xymatrix{
\A_{y,z}
\ar[r]^-{\A(1,h)}
&
\A_{y,z'}
\ar[r]^-{\A(x,-)}
&
[\A_{x,y},\A_{x,z'}]
}
$$
which is the codomain of $\TCcompt_{x, y, h:z \rightarrow z'}$.\\

\begin{tag}\label{pteqax1}
Proof of the equivalence of Axioms
\ref{cohassp}/\ref{cohassh}.
\end{tag}
The first of the 2-cell of Axiom \ref{cohassp}
decomposes as the product $\Xi_2 \circ \Xi_1$ 
where $\Xi_1$ is
{\tiny
$$
\xymatrix{
((\A_{t,u}  \A_{z,t})
\A_{y,z})  \A_{x,y}
\ar[r]^-{A'}
&
(\A_{t,u}  \A_{z,t})
 (\A_{y,z} \A_{x,y})
\ar[r]^-{1 \otimes \comp_{x,y,z}}
&
(\A_{t,u}  \A_{z,t}) \A_{x,z}
\ar@{=>}[r]^-{\assp_{x,z,t,u}}
&
\A_{x,u}
}
$$
}
and $\Xi_2$ is 
{\tiny
$$\xymatrix{
((\A_{t,u}  \A_{z,t})
 \A_{y,z})  \A_{x,y}
\ar[r]^-{(\comp_{z,t,u} \otimes 1) \otimes 1}
&
(\A_{z,u}  \A_{y,z})  \A_{x,y}
\ar@{=>}[r]^-{\assp_{x,y,z,u}}
&
\A_{x,u}
}
$$
}

The 2-cell $\Xi_1$ has a strict domain which image 
by $\Res$ is strict since the arrows $A'$, $\Res(A')$
and $\comp$ are strict. 
According to Lemma \cite{Sch08}-\calcdeux, the 2-cell $\Xi_1$ 
has image by $\Res \circ \Res$
{\tiny
$$\xymatrix{
\A_{t,u} \otimes \A_{z,t}
\ar@{=>}[r]^{\Res(\assp_{x,z,t,u})}
&
[\A_{x,z},\A_{x,u}]
\ar[r]^-{[\A_{x,y},-]}
&
[[\A_{x,y},\A_{x,z}],[\A_{x,y},\A_{x,u}]]
\ar[r]^-{[\A(x,-),1]}
&
[\A_{y,z},[\A_{x,y},\A_{x,u}]]
.}$$
}
This 2-cell has again strict domain and its 
image by $\Res$ is the 2-cell
{\tiny
$$\xymatrix{
\A_{t,u} 
\ar@{=>}[r]^-{\assh}
&
[\A_{z,t},[\A_{x,z},\A_{x,u}]]
\ar[r]^-{[1,[\A_{x,y},-]]}
&
[\A_{z,t},[[\A_{x,y},\A_{x,z}],[\A_{x,y},\A_{x,u}]]]
\ar[r]^-{[1,[\A(x,-),1]]}
&
[\A_{z,t},[\A_{y,z},[\A_{x,y},\A_{x,u}]]].
}
$$
}
The image by $\Res$ of $\Xi_2$ 
is 
{\tiny
$$\xymatrix{
(\A_{t,u} \otimes \A_{z,t}) \otimes \A_{y,z}
\ar[r]^-{\comp_{z,t,u} \otimes 1}
&
\A_{z,u} \otimes \A_{y,z}
\ar@{=>}[r]^-{\Res(\assp_{x,y,z,u})}
&
[\A_{x,y},\A_{x,u}]
}
$$
}
which has image by $\Res$ the 2-cell
{\tiny
$$\xymatrix{
\A_{t,u} \otimes \A_{z,t}
\ar[r]^-{\comp_{z,t,u}}
&
\A_{z,u}
\ar@{=>}[r]^-{\assh_{x,y,z,u}}
&
[\A_{y,z},[\A_{x,y},\A_{x,u}]]
}
$$
}
which has image by $\Res$ the 2-cell
{\tiny
$$
\xymatrix{
\A_{t,u}
\ar[r]^-{\A(z,-)}
&
[\A_{z,t}, \A_{z,u}]
\ar@{=>}[r]^-{[1,\assh_{x,y,z,u}]}
&
[\A_{z,t},[\A_{y,z},[\A_{x,y},\A_{x,u}]]].
}
$$
}
The second 2-cell from Axiom \ref{cohassp} decomposes
as $\Xi_5 \circ \Xi_4 \circ \Xi_3$ where 
$\Xi_3$ is
{\tiny
$$
\xymatrix{
((\A_{t,u} \A_{z,t}) \A_{y,z})
\A_{x,y}
\ar[r]^-{A' \otimes 1}
&
(\A_{t,u} (\A_{z,t} \A_{y,z}) )\A_{x,y}
\ar[r]^-{A'}
&
\A_{t,u} ((\A_{z,t} \A_{y,z}) )\A_{x,y})
\ar@{=>}[r]^-{1 \otimes \assp_{x,y,z,t}}
&
\A_{t,u} \A_{x,t}
\ar[r]^-{\comp}
&
\A_{x,u}
}
$$
}
$\Xi_4$ is
{\tiny
$$
\xymatrix{
((\A_{t,u} \A_{z,t}) \A_{y,z})
\A_{x,y}
\ar[r]^-{A' \otimes 1}
&
(\A_{t,u} (\A_{z,t} \A_{y,z}) )\A_{x,y}
\ar[r]^-{(1 \otimes \comp_{y,z,t}) \otimes 1}
&
(\A_{t,u} \A_{y,t}) \A_{x,y}
\ar@{=>}[r]^-{\assp_{x,y,t,u}}
&
\A_{x,u}
}
$$
}
and $\Xi_5$ is
{\tiny
$$\xymatrix{
((\A_{t,u}
\A_{z,t})
\A_{y,z})
\A_{x,y}
\ar@{=>}[r]^-{\assp_{y,z,t,u} \otimes 1}
&
\A_{y,u} \A_{x,y}
\ar[r]^-{\comp_{x,y,u}}
&
\A_{x,u}
}
$$
}
The image by $\Res \circ \Res \circ \Res$ of
the 2-cell $\Xi_3$ is 
{\tiny
$$
\xymatrix{
\A_{t,u} 
\ar[r]^-{\eta}
&
[\A_{z,t} \A_{y,z}, \A_{t,u} (\A_{z,t} \A_{y,z})]
\ar[r]^-{\Res}
&
[\A_{z,t},[\A_{y,z},\A_{t,u} ( \A_{z,t} \A_{y,z})]
\ar[r]^-{[1,[1,\Res(A')]]}
&
...}$$
$$\xymatrix{...\;
[\A_{z,t},[\A_{y,z},[\A_{x,y}, \A_{t,u} ((\A_{z,t} \A_{y,z}) \A_{x,y})]]]
\ar@{=>}[r]^-{[1,[1,[1,1 \otimes \assp_{x,y,z,t}]]]} 
&
[\A_{z,t},[\A_{y,z},[\A_{x,y},\A_{t,u} \A_{x,t}]]]
\ar[r]^-{[1,[1,[1,\comp_{x,t,u}]]]}
&
[\A_{z,t},[\A_{y,z},[\A_{x,y},\A_{x,u}]]]
}
$$
}
which rewrites
{\tiny
\begin{tabbing} 
\=1. \=$\xymatrix{
\A_{t,u} 
\ar[r]^-{\eta}
&
[\A_{z,t}\A_{y,z},\A_{t,u}(\A_{z,t}\A_{y,z})]
\ar[r]^-{[1,\Res(A')]}
&
[\A_{z,t}\A_{y,z},[\A_{x,y},\A_{t,u}((\A_{z,t}\A_{y,z})\A_{x,y})]]
\ar@{=>}[r]^-{[1,[1,1 \otimes \assp_{x,y,z,t}]]}
&
...}$\\
\> \>$\xymatrix{...\; 
[\A_{z,t}\A_{y,z},[\A_{x,y},\A_{t,u}\A_{x,t}]]
\ar[r]^-{[1,[1,\comp_{x,t,u}]]}
&
[\A_{z,t}\A_{y,z},[\A_{x,y},\A_{x,u}]]
\ar[r]^-{\Res}
&
[\A_{z,t},[\A_{y,z},[\A_{x,y},\A_{x,u}]]]
}$\\
\>2. \>$\xymatrix{
\A_{t,u}
\ar[r]^-{\eta}
&
[(\A_{z,t}\A_{y,z})\A_{x,y},\A_{t,u}((\A_{z,t}\A_{y,z}) \A_{x,y}) ]
\ar[r]^-{\Res}
&
[\A_{z,t}\A_{y,z},[\A_{x,y},\A_{t,u}((\A_{z,t}\A_{y,z}) \A_{x,y})] ]
\ar@{=>}[r]^-{  [1,[1, 1 \otimes \assp_{x,y,z,t}]]   }
&...}$\\
\> \>$\xymatrix{...\;
[\A_{z,t}\A_{y,z},[\A_{x,y},\A_{t,u}\A_{x,t}]]
\ar[r]^-{[1,[1,\comp_{x,t,u}]]}
&
[\A_{z,t}\A_{y,z},[\A_{x,y},\A_{x,u}]]
\ar[r]^-{\Res}
&
[\A_{z,t},[\A_{y,z},[\A_{x,y},\A_{x,u}]]
}
$\\
\>3. \>$\xymatrix{ 
\A_{t,u}
\ar[r]^-{\eta}
&
[(\A_{z,t}\A_{y,z}) \A_{x,y},  \A_{t,u} ((\A_{z,t}\A_{y,z}) \A_{x,y})]
\ar@{=>}[r]^-{[1, 1 \otimes \assp_{x,y,z,t}]}
&
[(\A_{z,t}\A_{y,z})\A_{x,y}, \A_{t,u}\A_{x,t}]
\ar[r]^-{[1,\comp_{x,t,u}]}
&
...}$\\
\> \>$\xymatrix{...\;
[(\A_{z,t}\A_{y,z})\A_{x,y},\A_{x,u}]
\ar[r]^-{\Res}
&
[\A_{z,t}\A_{y,z}, [\A_{x,y},\A_{x,u}]]
\ar[r]^-{\Res}
&
[\A_{z,t},[\A_{y,z},[\A_{x,y},\A_{x,u}]]]
}$\\\>4.
\>$\xymatrix{ 
\A_{t,u}
\ar[r]^-{\eta}
&
[\A_{x,t},  \A_{t,u}\A_{x,t}]
\ar@{=>}[r]^-{ [\assp_{x,y,z,t},1] }
&
[(\A_{z,t}\A_{y,z})\A_{x,y}, \A_{t,u}\A_{x,t}]
\ar[r]^-{[1,\comp_{x,t,u}]}
&
[(\A_{z,t}\A_{y,z})\A_{x,y},\A_{x,u}]
\ar[r]^-{\Res}
&
...}$\\
\>\>$\xymatrix{...\;
[\A_{z,t}\A_{y,z}, [\A_{x,y},\A_{x,u}]]
\ar[r]^-{\Res}
&
[\A_{z,t},[\A_{y,z},[\A_{x,y},\A_{x,u}]]]
}$\\
\>5. \>$\xymatrix{ 
\A_{t,u}
\ar[r]^-{\A(x,-)}
&
[\A_{x,t}, \A_{x,u}]
\ar[r]^-{[\assp_{x,y,z,t},1]}
&
[(\A_{z,t} \A_{y,z}) \A_{x,y},\A_{x,u}]
\ar@{=>}[r]^-{[\A_{x,y},-] } 
&
[[\A_{x,y},(\A_{z,t}\A_{y,z})\A_{x,y}], [\A_{x,y},\A_{x,u}]]
\ar[r]^-{[\eta,1]}
&
...}$\\
\>\>$\xymatrix{...\;
[\A_{z,t}\A_{y,z}, [\A_{x,y},\A_{x,u}]]
\ar[r]^-{\Res}
&
[\A_{z,t},[\A_{y,z},[\A_{x,y},\A_{x,u}]]]
}$\\
\>6. \>$\xymatrix{ 
\A_{t,u}
\ar[r]^-{\A(x,-)}
&
[\A_{x,t}, \A_{x,u}]
\ar[r]^-{[\A_{x,y},-]}
&
[[\A_{x,y},\A_{x,t}],[\A_{x,y},\A_{x,u}]]
\ar@{=>}[r]^-{[[1,\assp_{x,y,z,t}],1] } 
&
[[\A_{x,y},(\A_{z,t}\A_{y,z})\A_{x,y}], [\A_{x,y},\A_{x,u}]]
\ar[r]^-{[\eta,1]}
&
...}$\\
\>\>$\xymatrix{...\;
[\A_{z,t}\A_{y,z}, [\A_{x,y},\A_{x,u}]]
\ar[r]^-{\Res}
&
[\A_{z,t},[\A_{y,z},[\A_{x,y},\A_{x,u}]]]
}$\\
\>7. \>$\xymatrix{ 
\A_{t,u}
\ar[r]^-{\A(x,-)}
&
[\A_{x,t}, \A_{x,u}]
\ar[r]^-{[\A_{x,y},-]}
&
[[\A_{x,y},\A_{x,t}],[\A_{x,y},\A_{x,u}]]
\ar@{=>}[r]^-{[\Res(\assp_{x,y,z,t}),1] } 
&
[\A_{z,t}\A_{y,z}, [\A_{x,y},\A_{x,u}]]
\ar[r]^-{[\A_{y,z},-]}
&
...}$\\
\>\>$\xymatrix{...\;
[[\A_{y,z}, \A_{z,t} \A_{y,z}],
[\A_{y,z},[\A_{x,y},\A_{x,u}]]]
\ar[r]^-{[\eta,1]}
&
[\A_{z,t},[\A_{y,z},[\A_{x,y},\A_{x,u}]]]
}$\\
\>8. \>$\xymatrix{ 
\A_{t,u}
\ar[r]^-{\A(x,-)}
&
[\A_{x,t}, \A_{x,u}]
\ar[r]^-{[\A_{x,y},-]}
&
[[\A_{x,y},\A_{x,t}],[\A_{x,y},\A_{x,u}]]
\ar[r]^-{[\A_{y,z},-]}
&
[[\A_{y,z},[\A_{x,y},\A_{x,t}]],[\A_{y,z},[\A_{x,y},\A_{x,u}]]]
\;...}$\\
\>\>$\xymatrix{...
\ar@{=>}[r]^-{[\assh_{x,y,z,t},1]}
&
[\A_{z,t},[\A_{y,z},[\A_{x,y},\A_{x,u}]]]
}$
\end{tabbing}
}

The image by $\Res$ of the 2-cell 
$\Xi_4$ is 
$$
\xymatrix@C=3pc{
(\A_{t,u} \A_{z,t})\A_{y,z}
\ar[r]^-{A'}
&
\A_{t,u} (\A_{z,t} \A_{y,z})
\ar[r]^-{1 \otimes \comp_{y,z,t}}
&
\A_{t,u} \A_{y,t}
\ar@{=>}[r]^-{\Res(\assp_{x,y,t,u})}
&
[\A_{x,y} ,\A_{x,u}  ]
}
$$
which according to Lemma \cite{Sch08}-\calcdeux 
as image by $\Res \circ \Res$ the 2-cell
$$
\xymatrix{
\A_{t,u}
\ar@{=>}[r]^-{\assh_{x,y,t,u}}
&
[\A_{y,t},[\A_{x,y},\A_{x,u}]]
\ar[r]^-{[\A_{y,z},-]}
&
[[\A_{y,z},\A_{y,t}],[\A_{y,z},[\A_{x,y}, A_{x,u}]]]
\ar[r]^-{[\A(y,-),1]}
&
[\A_{z,t},[\A_{y,z},[\A_{x,y}, \A_{x,u}]]].
}
$$

The image by $\Res \circ \Res \circ \Res$ of 
the 2-cell $\Xi_5$ is 
$$
\xymatrix@C=3pc{
\A_{t,u}
\ar@{=>}[r]^-{\assh_{y,z,t,u}}
&
[\A_{z,t},[\A_{y,z},\A_{y,u}]]
\ar[r]^-{[1,[1,\A(x,-)]]}
&
[\A_{z,t},[\A_{y,z},[\A_{x,y},\A_{x,u}]]]
}
$$
\epf

\begin{tag}\label{pteqax2}
Proof of the equivalence of Axioms \ref{cohunassp} and \ref{cohunassh}. 
\end{tag}
\pf
It is easy to check that the image by $\Res$ 
of the 2-cell $\Xi_1$ of Axiom \ref{cohunassp} 
is ${A(x,-)}_{y,z} * \rh_{y,z}$
the first 2-cell of the Axiom \ref{cohunassh}.\\

The image by $\Res$ of the 2-cell 
$\Xi_2$ is 
$$
\xymatrix{
\A_{y,z} 
\ar[r]^-{\A(x,-)}
&
[\A_{x,y},\A_{x,z}]
\ar[d]_{[\ev_{\gen},1]}
\ar@{-}[rr]^-{id}
& 
\ar@{=>}[d]|{[\lh,1]}
&
[\A_{x,y},\A_{x,z}]
\\
&
[[\unc,\A_{x,y}],\A_{x,z}]
\ar[rr]_-{[[\unit_y,1],1]}
& & 
[[\A_{y,y},\A_{x,y}],\A_{x,z}]
\ar[u]_{[\A(-,y),1]}
}
$$
The 2-cell $\Xi_3$, namely 
$$
\xymatrix{
\A_{y,z} \otimes \A_{x,y}
\ar[r]^-{R' \otimes 1}
&
(\A_{y,z} \otimes \unc  ) \otimes \A_{x,y}
\ar[r]^-{(1 \otimes \unit_y) \otimes 1}
&
\A_{y,z} \otimes \A_{y,y}
\ar@{=>}[r]^-{\assp_{x,y,y,z}}
&
\A_{x,z}
}
$$
has image by $\Res$ the 2-cell
$$
\xymatrix{
\A_{y,z}
\ar[r]^-{R'}
&
\A_{y,z} \otimes \unc
\ar[r]^-{1 \otimes \unit_y}
&
\A_{y,z} \otimes \A_{y,y}
\ar[r]^-{\Res(\assp_{x,y,y,z})}
&
[\A_{x,y},\A_{x,z}]
}
$$
which is according to Lemma \ref{calc51}
$$
\xymatrix{
\A_{y,z}
\ar@{=>}[r]^-{\assh_{x,y,y,z}}
&
[\A_{y,y},[\A_{x,y},\A_{x,z}]]
\ar[r]^-{[\unit_y,1]}
&
[\unc,[\A_{x,y},\A_{x,z}]]
\ar[r]^-{[\ev_{\gen},1]}
&
[\A_{x,y},\A_{x,z}]
}
$$
\epf

\begin{tag}\label{eqcofun1hp}
Equivalence of Axioms \ref{cohfun1p} and \ref{cohfun1h}.
\end{tag}
First let us remark that the two 2-cells of Axiom \ref{cohfun1h} 
have the same domain. This results from the following sequence of 
equal arrows.
{\tiny
\begin{tabbing}
\=1. \=
$\xymatrix{
[\B_{Fx,Fz}, \B_{Fx,Ft}]
\ar[r]^-{[\B_{Fx,Fy},-]}
& 
[[\B_{Fx,Fy},\B_{Fx,Fz}], [\B_{Fx,Fy},\B_{Fx,Ft}]]
\ar[r]^-{[\B(Fx,-),1]}
&
[\B_{Fy,Fz},[\B_{Fx,Fy},\B_{Fx,Ft}]]
\ar[r]^-{[F_{y,z},1]} & 
...}$\\
\> \>$\xymatrix{ ...\;
[\A_{y,z},[\B_{Fx,Fy},\B_{Fx,Ft}]]
\ar[r]^-{[1,[F_{x,y},1]]}
&
[\A_{y,z},[\A_{x,y},\B_{Fx,Ft}]]}$\\
\>2. \>
$\xymatrix{
[\B_{Fx,Fz}, \B_{Fx,Ft}]
\ar[r]^-{[\B_{Fx,Fy},-]}
& 
[[\B_{Fx,Fy},\B_{Fx,Fz}], [\B_{Fx,Fy},\B_{Fx,Ft}]]
\ar[r]^-{[\B(Fx,-),1]}
&
[\B_{Fy,Fz},[\B_{Fx,Fy},\B_{Fx,Ft}]]
\ar[r]^-{[1,[F_{x,y},1]]}
&
...}$\\
\> \>$\xymatrix{...\;
[\B_{Fy,Fz},[\A_{x,y},\B_{Fx,Ft}]]
\ar[r]^-{[F_{y,z},1]}
&
[\A_{y,z},[\A_{x,y},\B_{Fx,Ft}]]}$\\
\>3. \>
$\xymatrix{
[\B_{Fx,Fz}, \B_{Fx,Ft}]
\ar[r]^-{[\B_{Fx,Fy},-]}
& 
[[\B_{Fx,Fy},\B_{Fx,Fz}], [\B_{Fx,Fy},\B_{Fx,Ft}]]
\ar[r]^-{[1,[F_{x,y},1]]}
&
[[\B_{Fx,Fy}, \B_{Fx,Fz}],[\A_{x,y},\B_{Fx,Ft}]]
\ar[r]^-{[\B(Fx,-),1]} & ... }$\\
\> \>$\xymatrix{... \;
[\B_{Fy,Fz},[\A_{x,y},\B_{Fx,Ft}]]
\ar[r]^-{[F_{y,z},1]}
&
[\A_{y,z},[\A_{x,y},\B_{Fx,Ft}]]}$\\
\>4. \>
$\xymatrix{
[\B_{Fx,Fz}, \B_{Fx,Ft}]
\ar[r]^-{[\A_{x,y},-]}
& 
[[\A_{x,y},\B_{Fx,Fz}], [\A_{x,y},\B_{Fx,Ft}]]
\ar[r]^-{[[F_{x,y},1],1]}
&
[[\B_{Fx,Fy}, \B_{Fx,Fz}],[\A_{x,y},\B_{Fx,Ft}]]
\ar[r]^-{[\B(Fx,-),1]} & ...  }$\\
\> \> $\xymatrix{... \;
[\B_{Fy,Fz},[\A_{x,y},\B_{Fx,Ft}]]
\ar[r]^-{[F_{y,z},1]}
&
[\A_{y,z},[\A_{x,y},\B_{Fx,Ft}]].}$
\end{tabbing} 
}
In the above derivation arrows 3. and 4. are equal
according to Lemma \cite{Sch08}-\impronatPosttrois.\\

The 2-cell 
{\small
$$
\xymatrix@C=3pc{
(\A_{z,t} \otimes \A_{y,z}) \otimes \A_{x,y}
\ar[r]^-{(F_{z,t} \otimes F_{y,z}) \otimes 1}
&
(\B_{Fz,Ft} \otimes \B_{Fy,Fz}) \otimes \A_{x,y}
\ar[r]^-{1 \otimes F_{x,y}}
& 
(\B_{Fz,Ft} \otimes \B_{Fy,Fz}) \otimes \B_{Fx,Fy}
\ar@{=>}[r]^-{\assp_{Fx,Fy,Fz,Ft}}
&
\B_{Fx,Ft}
}
$$
}
has a strict domain and image by $\Res$
{\small
$$
\xymatrix@C=3pc{
\A_{z,t} \otimes \A_{y,z}
\ar[r]^-{F_{z,t} \otimes F_{y,z}}
&
\B_{Fz,Ft} \otimes \B_{Fy,Fz}
\ar[r]^-{\Res(\assh_{Fx,Fy,Fz,Ft})}
& 
[\B_{Fx,Fy}, \B_{Fx,Ft}]
\ar@{=>}[r]^{[F_{x,y},1]}
&
[\A_{x,y},\B_{Fx,Ft}]
}
$$
}
which has a strict domain and image by $\Res$
{\small
$$
\xymatrix{
\A_{z,t}
\ar[r]^-{F_{z,t}}
& 
\B_{Fz,Ft}
\ar@{=>}[r]^-{\assh}
&
[\B_{Fy,Fz},[\B_{Fx,Fy},\B_{Fx,Ft}]]
\ar[r]^-{[F_{y,z},1]}
&
[\A_{y,z},[\B_{Fx,Fy},\B_{Fx,Ft}]]
\ar[r]^-{[1,[F_{x,y},1]]}
&
[\A_{y,z},[\A_{x,y},\B_{Fx,Ft}]].
}
$$
}
The 2-cell 
{\small
$$
\xymatrix{
(\A_{z,t} \otimes \A_{y,z}) \otimes \A_{x,y}
\ar@{=>}[r]^-{F^2_{y,z,t} \otimes 1}
&
\B_{Fy,Ft} \otimes \A_{x,y}
\ar[r]^-{1 \otimes F_{x,y}}
& 
\B_{Fy,Ft} \otimes \B_{Fx,Fy}
\ar[r]^-{\comp}
&
\B_{Fx,Ft}
}
$$
}
has image by $\Res$ that rewrites
{\small
\begin{tabbing}
\=$\xymatrix@C=3pc{
\A_{z,t} \otimes \A_{y,z}
\ar@{=>}[r]^-{F^2}
&
\B_{Fy,Ft}
\ar[r]^-{\Res(\comp \otimes (1 \otimes F_{x,y}))}
&
[\A_{x,y},\B_{Fx,Ft}]
}
$\\
\>$
\xymatrix{
\A_{z,t} \otimes \A_{y,z}
\ar@{=>}[r]^-{F^2}
&
\B_{Fy,Ft}
\ar[r]^-{\B(Fx,-)}
&
[\B_{Fx,Fy},\B_{Fx,Ft}]
\ar[r]^-{[F_{x,y},1]}
&
[\B_{Fx,Fy},\B_{Fx,Ft}]
}
$
\end{tabbing}
}
The image by $\Res$ of this last arrow is
{\small
$$
\xymatrix@C=3pc{
\A_{z,t} 
\ar@{=>}[r]^-{{F'}^2}
&
[\A_{y,z},\B_{Fy,Ft}]
\ar[r]^-{[1,\B(Fx,-)]}
&
[\A_{y,z},[\B_{Fx,Fy},\B_{Fx,Ft}]]
\ar[r]^-{[1,[F_{x,y},1]]}
&
[\A_{y,z},[\B_{Fx,Fy},\B_{Fx,Ft}]].
}
$$
}
The 2-cell
{\small
$$
\xymatrix{
(\A_{z,t} \otimes \A_{y,z}) \otimes \A_{x,y}
\ar[r]^-{\comp \otimes 1}
& 
\A_{y,t} \otimes \A_{x,y}
\ar@{=>}[r]^-{F^2_{x,y,t}}
&
\B_{Fx,Ft}
}
$$
}
has image by $\Res$ the 2-cell
{\small
$$
\xymatrix{
\A_{z,t} \otimes \A_{y,z}
\ar[r]^-{\comp}
&
\A_{y,t}
\ar@{=>}[r]^-{{F'}^2}
& 
[\A_{x,y},\B_{Fx,Ft}]
}
$$
}
which has image by $\Res$
{\small
$$
\xymatrix{
\A_{z,t}
\ar[r]^-{\A(y,-)}
& 
[\A_{y,z}, \A_{y,t}]
\ar@{=>}[r]^-{[1, {F'}^2_{x,y,t}]}
&
[\A_{y,z},[\A_{x,y}, \B_{Fx,Ft}]].
}
$$
}

The 2-cell
{\small
$$
\xymatrix{
\A_{z,t}
\otimes 
\B_{Fx,Fz}
\ar[r]^-{F_{z,t} \otimes 1}
&
\B_{Fz,Ft} \otimes \B_{Fx,Fz}
\ar[r]^-{\comp_{Fx,Fz,Ft}}
&
\B_{Fx,Ft}
}
$$
}
has image by $\Res$ the 2-cell
{\small
$$
\xymatrix{
\A_{z,t}
\ar[r]^-{F_{z,t}}
&
\B_{Fz,Ft}
\ar[r]^-{\B(Fx,-)}
&
[\B_{Fx,Fz}, \B_{Fx, Ft}]
}
$$
}
therefore according to Lemma \ref{calc1} the 2-cell
{\tiny
$$
\xymatrix{
(\A_{z,t} \otimes \A_{y,z}) \otimes \A_{x,y}
\ar[r]^{A'}
& 
\A_{z,t} \otimes (\A_{y,z} \otimes \A_{x,y} )
\ar@{=>}[r]^{1 \otimes F^2_{x,y,z}}
& 
\A_{z,t} \otimes \B_{Fx,Fz}
\ar[r]^{F_{z,t} \otimes 1}
&
\B_{Fz,Ft} \otimes \B_{Fx,Fz}
\ar[r]^{\comp_{Fx,Fz,Ft}}
&
\B_{Fx,Ft} 
}
$$
}
has image by $\Res \circ \Res$ the 2-cell
{\tiny
$$
\xymatrix{
\A_{z,t}
\ar[r]^-{F_{z,t}}
&
\B_{Fz,Ft}
\ar[r]^-{\B(Fx,-)}
&
[\B_{Fx,Fz},\B_{Fx,Ft}]
\ar[r]^-{[\A_{x,y},-]}
&
[[\A_{x,y},\B_{Fx,Fz}],[\A_{x,y},\B_{Fx,Ft}]]
\ar@{=>}[r]^-{[{F'}^2_{x,y,z},1]}
&
[\A_{y,z},[\A_{x,y},\B_{Fx,Ft}]]
.}
$$
}
According to Lemma \ref{calc1} the 2-cell
{\small
$$
\xymatrix{
(\A_{z,t} \otimes \A_{y,z}) \otimes \A_{x,y}
\ar[r]^-{A'}
&
\A_{z,t} \otimes (\A_{y,z} \otimes \A_{x,y})
\ar[r]^-{1 \otimes \comp_{x,y,z}}
&
\A_{z,t} \otimes \A_{x,z}
\ar@{=>}[r]^-{F^2_{x,z,t}}
&
\B_{Fx,Ft} 
}
$$
}
has image by $\Res \circ \Res$ the 2-cell
{\small
$$
\xymatrix@C=3pc{
\A_{z,t}
\ar@{=>}[r]^-{{F'}^2_{x,z,t}}
&
[\A_{x,z}, \B_{Fx,Ft}]
\ar[r]^-{[\A_{x,y},-]}
&
[[\A_{x,y},\A_{x,z}],[\A_{x,y} \B_{Fx,Ft}]]
\ar[r]^-{[\A(x,-),1]}
& 
[\A_{y,z},[\A_{x,y},\B_{Fx,Ft}]]
.}
$$
}
\epf

\begin{tag}\label{eqcofun2hp}
Equivalence of Axioms \ref{cohfun2p} and \ref{cohfun2h}.
\end{tag}
\pf
According to Lemma \ref{calc51}, the arrow
$$
\xymatrix{
\A_{x,y}
\ar[r]^-{R'}
& 
\A_{x,y} \otimes \unc
\ar[r]^-{1 \otimes \unit_x}
&
\A_{x,y} \otimes \A_{x,x}
\ar@{=>}[r]^-{F^2_{x,x,y}}
&
\A_{x,y}
}
$$
is equal to 
$$
\xymatrix{
\A_{x,y}
\ar@{=>}[r]^-{{F'}^2_{x,x,y}}
&
[\A_{x,x}, \B_{Fx,Fy}]
\ar[r]^-{[\unit_x,1]}
&
[\unc,\B_{Fx,Fy}]
\ar[r]^-{\ev_{\gen}}
&
\B_{Fx,Fy}.
}
$$
Since the image by $\Res$ of 
$$
\xymatrix{
\A_{x,y} \otimes \B_{Fx,Fx}
\ar[r]^-{F_{x,y} \otimes 1}
&
\B_{Fx,Fy} \otimes \B_{Fx,Fx}
\ar[r]^-{\comp}
&
\B_{Fx,Fy}
}
$$
is
$$
\xymatrix{
\A_{x,y}
\ar[r]^-{F_{x,y}}
& 
\B_{Fx,Fy}
\ar[r]^-{\B(Fx,-)}
& 
[\B_{Fx,Fx}, \B_{Fx,Fy}],
}
$$
the arrow 
$$
\xymatrix{
\A_{x,y}
\ar[r]^-{R'}
&
\A_{x,y} \otimes \unc
\ar@{=>}[r]^-{1 \otimes F^0_x}
&
\A_{x,y} \otimes \B_{Fx,Fx}
\ar[r]^-{F_{x,y} \otimes 1}
&
\B_{Fx,Fy} \otimes \B_{Fx,Fx}
\ar[r]^-{\comp}
&
\B_{Fx,Fy}
}
$$
is equal to 
$$
\xymatrix{
\A_{x,y}
\ar[r]^-{F_{x,y}}
&
\B_{Fx,Fy}
\ar[r]^-{\B(Fx,-)}
&
[\B_{Fx,Fx},\B_{Fx,Fy}]
\ar@{=>}[r]^-{[F^0_x,1]}
&
[\unc,\B_{Fx,Fy}]
\ar[r]^-{\ev_{\gen}}
&
\B_{Fx,Fy}
}
$$
according to Lemma \ref{calc51}.
\epf

\begin{tag}\label{eqcofun3hp}
Equivalence of Axioms \ref{cohfun3p} and \ref{cohfun3h}.
\end{tag}
\pf
The arrow 
$$\xymatrix{
\B_{Fy,Fy} \otimes \A_{x,y}
\ar[r]^-{1 \otimes F_{x,y}}
& 
\B_{Fy,Fy} \otimes \B_{Fx,Fy}
\ar[r]^-{\comp}
&
\B_{Fx,Fy}
}
$$
has image by $\Res$
$$
\xymatrix{
\B_{Fy,Fy}
\ar[r]^-{\B(Fx,-)}
&
[\B_{Fx,Fy}, \B_{Fx,Fy}]
\ar[r]^-{[F_{x,y},1]}
&
[\A_{x,y},\B_{Fx,Fy}]
}
$$
which has dual
$$
\xymatrix{
\A_{x,y}
\ar[r]^-{F_{x,y}}
& 
\B_{Fx,Fy}
\ar[r]^-{\B(-,Fy)}
&
[\B_{Fy,Fy},\B_{Fx,Fy}].
}
$$
Therefore according to Lemma \ref{calc6},
the 2-cell
$$
\xymatrix{
\A_{x,y}
\ar[r]^-{L'}
& 
\unc \otimes \A_{x,y}
\ar@{=>}[r]^-{F^0_y \otimes 1}
& 
\B_{Fy,Fy} \otimes \A_{x,y}
\ar[r]^-{1 \otimes F_{x,y}}
& 
\B_{Fy,Fy} \otimes \B_{Fx,Fy}
\ar[r]^-{\comp}
&
\B_{Fx,Fy}
}
$$
is equal to
$$
\xymatrix{
\A_{x,y}
\ar[r]^-{F_{x,y}}
&
\B_{Fx,Fy}
\ar[r]^-{\B(-,Fy)}
&
[\B_{Fy,Fy},\B_{Fx,Fy}]
\ar@{=>}[r]^{[F^0_y,1]}
& 
[\unc,\B_{Fx,Fy}]
\ar[r]^{\ev_{\gen}}
&
\B_{Fx,Fy}. 
}
$$
According to Lemma \ref{calc6} the 2-cell
$$
\xymatrix{
\A_{x,y}
\ar[r]^-{L'}
&
\unc \otimes \A_{x,y}
\ar[r]^-{\unit_y \otimes 1}
&
\A_{y,y} \otimes \A_{x,y}
\ar[r]^-{F^2_{x,y,y}}
&
\A_{x,y}
}
$$
is equal to 
$$
\xymatrix{
\A_{x,y}
\ar[r]^-{{({F'}^2_{x,y,y})}^*}
&
[\A_{y,y}, \A_{x,y}]
\ar[r]^-{[\unit_y,1]}
&
[\unc,\A_{x,y}]
\ar[r]^-{\ev_{\gen}}
& 
\A_{x,y}
}
$$
\epf

\noindent{\bf Section \ref{Exples}.}\\

We will need the following characterization of bilinear 
natural transformations.
\begin{remark}\label{monatABC}
For any symmetric monoidal functors
$F,G: \A \rightarrow [\B,\C]$ with respective 
underlying functors $F',G': \A \times \B \rightarrow \C$.
Any 2-cell $\sigma: F \rightarrow G: \A \rightarrow [\B,\C]$
in $\BS$ corresponds to a collection of arrows 
$\sigma_{a,b}: F'(a,b) \rightarrow G'(a,b)$ in $\C$, natural
in $a$ and $b$ and that satisfies the following two 
conditions:\\
- For all objects $a$, $a'$, $b$ and $b'$ of $\A$
the following diagram commutes
$$\xymatrix{
F'(a,b)+F'(a,b')
\ar[d]|{\sigma_{a,b} + \sigma_{a,b'}} 
\ar[r]^{{ {F^*}^2_{b,b'} }_a}
&
F'(a,b+b')
\ar[d]|{\sigma_{a,b+b'}}
\\
F'(a',b)+F'(a',b') 
\ar[r]^{{ {F^*}^2_{b,b'} }_a'}
&
F'(a',b+b');
}
$$ 
- For all objects $a$, $a'$, $b$ and $b'$ of $\A$
the following diagram commutes
$$
\xymatrix{
F'(a,b)+F'(a',b)
\ar[d]_{\sigma_{a,b} + \sigma_{a',b}} 
\ar[r]^{{ F^2_{a,a'} }_b}
&
F'(a+a',b)
\ar[d]_{\sigma_{a+a',b}}
\\
F'(a,b')+F'(a',b') 
\ar[r]_{{ F^2_{a,a'} }_{b'}}
&
F'(a+a',b')}
$$
\end{remark}

\begin{tag}\label{pfeqdef2rgobj}
Proof that 2-rings in the sense of \ref{defobjJP} are 
exactly one point $\BS$-categories.
\end{tag}
\pf
Let us start with a 2-ring $\A$ as defined in \ref{defobjJP}.
Its corresponds to a one-point $\BS$-category 
as follows. Since monoidal functors $\unc \rightarrow \A$ 
are in one-to-one correspondence with objects of $\A$, the 
object $\unm$ correspond to a strict arrow 
$\unit: \unc \rightarrow \A$ in $\BS$.
That the multiplication ``$.$'', which is already a functor 
$\A \times \A \rightarrow \A$, defines an arrow
$\multh: \A \rightarrow [\A,\A]$ in $\BS$ corresponds to
the existence of the natural arrows $\underline{a}_{b,b'}$
and $\overline{b}_{a,a'}$
and the commutation of Diagrams \ref{amofun1},
\ref{bmofun1}, \ref{asymofun}, \ref{bsymofun} and \ref{lastaxiom}.
Precisely for any object $a$ of $\A$, 
the multiplication on the left by $a$,
defines a functor $a.-: \A \rightarrow \A$.
That this one is monoidal corresponds to the 
existence of the arrows $\underline{a}_{b,b'}$ natural in
$b$ and $b'$ and such that Diagrams \ref{amofun1} 
commute for all $b$, $b'$ and $b''$. 
That the monoidal $a.-$ is symmetric corresponds to  
the commutation of \ref{asymofun} for all $b$ and $b'$.
That for any arrow $f: a \rightarrow a'$, the natural 
transformation $f.-: a.- \rightarrow a'.-: \A \rightarrow \A $ 
is monoidal for the structures described previously
on $a.-$ and $a'.-$ corresponds to the naturality in
the argument $a$ of the maps $\underline{a}_{b,b'}$. 
The assignments $a \mapsto a.-$ and 
$(f:a \rightarrow a') \mapsto (f.-: a.- \rightarrow a'.-)$
define therefore a functor say $\multh: \A \rightarrow \BS(\A,\A)$.
The existence of a natural transformation
$(a.-) + (a'.-) \rightarrow (a+a').-$
corresponds to the existence of arrows $\overline{b}_{a,a'}$
natural in $b$. This transformation is monoidal
since Diagrams \ref{lastaxiom} commute.
Eventually that the collection of these arrows
for all $a$ and $a'$ defines a symmetric monoidal
structure on the above functor $\multh$ results from the 
the naturality of the collection in the arguments
$a$ and $a'$ and the commutation Diagrams
\ref{bmofun1} and \ref{bsymofun}.\\

According to Remark \ref{monatABC}, a
2-cell $\assh$ in $\BS$ as in \ref{assh} corresponds
to a natural collection of arrows $\assmu_{a,b,c}: 
a.(b.c) \rightarrow (a.b).c$  in $\A$ such that 
Diagrams \ref{bJP1}, \ref{bJP2} and \ref{bJP3} commute.
A 2-cell $\rh: 1 \rightarrow \ev_{\gen} \circ [\unit,1] \circ \multh: 
\A \rightarrow \A$ in $\BS$
corresponds to a natural collection of arrows 
$\rmu_a : a.\unm \rightarrow a$ such that Diagrams
\ref{bJP4} commute. Similarly a 2-cell 
$\lh: 1 \rightarrow \ev_{\gen} \circ [\unit,1] \circ 
\multh^*: \A \rightarrow \A$ amounts to a natural collection of 
arrows $\lmu_a : \unm. a \rightarrow a$ such that Diagrams
\ref{bJP5} commute.\\

Then the coherence conditions \ref{cohassh} and \ref{cohunassh} 
for $\assh$, $\rh$ and $\lh$ above
are equivalent to the coherence conditions for the associativity 
and unit laws of the monoidal category $(\A,.,\un,\assmu, \rmu,\lmu)$. 
\epf

\begin{tag}\label{pfeqdef2rgmor}
Proof that 2-ring morphisms in the sense 
of \ref{defmorJP} are exactly $\BS$-functors.
\end{tag}
\pf
According to the two observations below the result becomes 
clear after inspection 
of Axioms \ref{cohfun1h}, \ref{cohfun2h} and \ref{cohfun3h}.\\

Observe first that according to Remark \ref{monatABC} 
a 2-cell 
$$\xymatrix{
\A 
\ar[dd]_H
\ar[r]^-\multh
& 
[\A,\A]
\ar[rd]^{[1,H]}
\\
& 
& 
[\A,\B]
\\
\B
\ar[r]_-{\multh}
\ar@{=>}[ruu]
&
[\B,\B]
\ar[ru]_{[\multh,1]}
}$$ in $\BS$
is the same thing than a collection of arrows 
$\Xi_{a,b} : Ha. Hb \rightarrow H(a+b)$
natural in $a$ and $b$ and satisfying 
the conditions that for any objects $a$, $b$, $c$ of $\A$,
the two diagrams in $\B$ below commute
\begin{tag}\label{diag1}
$$
\xymatrix{
H(a).H(b)
+
H(a).H(c)
\ar[r]^-{\underline{H(a)}_{H(b),H(c)} }
\ar[d]_{ \Xi_{a,b} + \Xi_{a,c} }
&
H(a).(H(b)+H(c))
\ar[r]^-{H(a).H^2_{b,c}}
&
H(a).H(b+c) 
\ar[d]^{\Xi_{a,b+c}}
\\
H(a. b) + H(a.c)
\ar[r]_-{H^2_{a.b,a.c}}
&
H(a.b + a.c)
\ar[r]_-{H(\underline{a}_{b,c})}
& 
H(a.(b + c))
}
$$
\end{tag}
\begin{tag}\label{diag2}
$$
\xymatrix{
H(a).H(c)
+
H(b).H(c)
\ar[r]^-{\overline{H(c)}_{H(a),H(b)} }
\ar[d]_{ \Xi_{a,c} + \Xi_{b,c} }
&
(H(a)+H(b)).H(c)
\ar[r]^{H^2_{a,b}.H(c)}
&
H(a+b).H(c) 
\ar[d]^{\Xi_{a + b , c}}
\\
H(a.c) + H(b.c)
\ar[r]_{H^{2}_{a.c,b.c}}
&
H(a.c + b.c)
\ar[r]_{H(\overline{c}_{a,b})}
& 
H((a + b).c)
}
$$  
\end{tag}
where we write the $\multh$ as products. 
Also according to Remark \ref{evgen2cel} in Appendix 
any 2-cell
$$
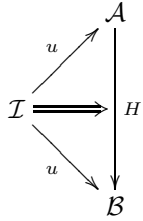
\xymatrix{
& \A \ar[dd]^H \\
\unc \ar[ru]^{\unit} 
\ar[rd]_{\unit}
\ar@{=>}[r]
& \\
& \B
}$$ in $\BS$ 
is fully determined by its component
at the generator $\gen$ which is an arrow 
$1_{\B} \rightarrow H(1_{\A})$ if $1_{\A}$ and
$1_\B$ denote the images $\unit(\gen)$, 
and conversely
every such arrow corresponds to a 2-cell as above in this
way.
\epf



\begin{tag}\label{pfPicsEn}
Proof of Proposition \ref{PicsEn}.
\end{tag} 
\pf
That there are identity 2-cells \ref{assh} amounts to the fact
that the diagram in $\BS$
$$
\xymatrix{
& 
[\C,\D]
\ar[ld]_-{[\B,-]}
\ar[rd]^-{[\A,-]}
&
\\
[[\B,\C],[\B,\D]]
\ar[d]_{[1,[\A,-]]}
&
& 
[[\A,\C],[\A,\D]]
\ar[d]^{[[\A,\B],-]}
\\
[[\B,\C], [[\A,\B],[\A,\D]]
&
& 
[[[\A,\B],[\A,\C]], [[\A,\B],[\A,\D]]
\ar[ll]^-{[[\A,-],1]}
}
$$ 
commutes
for any objects $\A$, $\B$, $\C$ and $\D$.
Such a diagram involves only strict
arrows in $\BS$ and its underlying diagram in $\Cat$ is
$$
\xymatrix{
& 
\BS(\C,\D)
\ar[ld]_-{\Post}
\ar[rd]^-{\Post}
&
\\
\BS([\B,\C],[\B,\D])
\ar[d]_{\BS(1,[\A,-])}
&
& 
\BS([\A,\C],[\A,\D])
\ar[d]^{\Post}
\\
\BS([\B,\C], [[\A,\B],[\A,\D]])
& 
& 
\BS( [[\A,\B],[\A,\C]], [[\A,\B],[\A,\D]] )
\ar[ll]^-{\BS([\A,-],1)}
}
$$ 
which commutes according to Lemma 
\cite{Sch08}-\impronatPost.\\

One has identity 2-cells $\rh$ 
since for any objects $\A$ and $\B$ in $\BS$ the composite
$$\xymatrix{
[\A,\B]
\ar[r]^-{[\A,-]}
&
[[\A,\A], [\A,\B]]
\ar[r]^-{[\vv,1]}
&
[\unc,[\A,\B]]
\ar[r]^-{\ev_{\gen}}
&
[\A,\B]
}$$
which is strict and has underlying functor
that is an identity, hence is an identity in $\BS$.\\ 

One has identity 2-cells $\lh$ since the composite
$$\xymatrix{
[\A,\B]
\ar[r]^-{[-,\B]}
&
[[\B,\B],[\A,\B]]
\ar[r]^-{[\vv,1]}
&
[\unc,[\A,\B]]
\ar[r]^-{\ev_{\gen}}
&
[\A,\B]
}$$
is the identity at $[\A,\B]$ as shown below.
The arrow 
$\xymatrix{[\A,\B]
\ar[r]^-{[-,\B]}
&
[[\B,\B],[\A,\B]]
\ar[r]^-{[\vv,1]}
&
[\unc,[\A,\B]]}$
is $\vs$ since it has dual
$$\xymatrix{
\unc
\ar[r]^-{\vv}
&
[\B,\B]
\ar[r]^-{[\A,-]}
&
[[\A,\B],[\A,\B]]
}$$
that is $\vv$ according to
Lemma \cite{Sch08}-\lemvvun. One concludes 
since $\ev_{\gen} \circ \vs = 1$.
\epf

\begin{tag} Proof of Proposition \ref{unit2rg} \end{tag}\label{pfunit2rg}
\pf
That one has an identity 2-cell \ref{assh}
amounts to the commutation of the diagram
$$
\xymatrix{
& 
\unc
\ar[rd]^{\vv}
\ar[ld]_{\vv}
&
\\
[\unc,\unc]
\ar[d]_{[\unc,-]}
& & 
[\unc,\unc]
\ar[d]^{[1,\vv]}
\\
[[\unc,\unc],[\unc,\unc]]
\ar[rr]_{[\vv,1]}
& & 
[\unc,[\unc,\unc]]
.}
$$
This diagram commutes since 
the arrow
$$\xymatrix{
\unc
\ar[r]^-{\vv}
&
[\unc, \unc]
\ar[r]^-{[\unc,-]}
&
[[\unc,\unc], [\unc,\unc]]
\ar[r]^-{[\vv,1]}
&
[\unc,[\unc,\unc]]
}$$
is
$$\xymatrix{
\unc
\ar[r]^-{\vv}
&
[[\unc,\unc],[\unc,\unc]]
\ar[r]^-{[\vv,1]}
&
[\unc,[\unc,\unc]]
}$$
according to Lemma \cite{Sch08}-\lemvvun,
which is equal to 
$$\xymatrix{
\unc
\ar[r]^-{\vv}
&
[\unc,\unc]
\ar[r]^-{[1,\vv]}
&
[\unc,[\unc,\unc]].
}$$
according Lemma \cite{Sch08}-\lemvvdeux.\\

One has an identity 2-cell $\rh$ since 
the composite 
$\xymatrix{
\unc 
\ar[r]^{\vv}
& 
[\unc,\unc]
\ar[r]^{\ev_{\gen}}
&
\unc
}$ 
is the identity.\\

One has also an identity 2-cell $\lh$
since the composite 
$$
\xymatrix{
\unc 
\ar[r]^{\vs}
&
[\unc,\unc]
\ar[r]^{\ev}
&
\unc
}$$ 
is an identity according to Lemma \cite{Sch08}-\vvgen.
\epf

\begin{tag}\label{pflemtoto}
Proof of Lemma \ref{lemtoto}. \end{tag}
\pf
According to Lemma \cite{Sch08}-\homevtroisdual, the diagram in $\BS$ 
$$
\xymatrix{
[\unc,[\A,\B]]
\ar[r]^-{\ev_{\gen}}
\ar[d]_{[-,[\A,\C]]}
&
[\A,\B]
\ar[d]^{q}
\\
[[[\A,\B],[\A,\C]],[\unc,[\A,\C]]
\ar[r]_{[1,\ev_{\gen}]}
&
[[[\A,\B],[\A,\C]], [\A,\C]]
}$$
commutes, therefore also does the diagram 
$$
\xymatrix{
[[\A,\B],[\A,\C]]
\ar[r]^{[\tilde{F},1]}
\ar[rd]_{\ev_{\gen}}
&
[\unc'[\A,\C]]
\ar[d]^{\ev_{\gen}}
\\
&
[\A,\C]
.}
$$
From this and according to Corollary \cite{Sch08}-\evhomun
all diagrams in the pasting below
$$
\xymatrix{
[\B,\C]
\ar[r]^{[F,\C]}
\ar[d]_{[\A,-]}
&
[\A,\C]
\\
[[\A,\B],[\A,\C]]
\ar[r]_-{[\tilde{F},1]}
\ar[ru]|{\ev_F}
&
[\unc,[\A,\C]]
\ar[u]_{\ev_{\gen}}
}
$$
commute and according to Lemma \ref{calc51} the commutation of the 
external diagram above is equivalent to the commutation of the first 
diagram of the Lemma.\\
In the pasting 
$$\xymatrix{
[\C,\A]
\ar[r]^{[\C,F]}
\ar[d]_{[-,\B]}
&
[\C,\B]
\\
[[\A,\B],[\C,\B]]
\ar[r]_-{[\tilde{F},1]}
\ar[ru]|{\ev_{F}}
&
[\unc,[\C,\B]]
\ar[u]_{\ev_{\gen}}
}$$
the top-left diagram commutes according to Corollary
\cite{Sch08}-\evhomdeux and  we have already seen that the bottom-left diagram
commutes. The commutation of the external diagram above is
equivalent to the commutation of the second diagram of the Lemma
according to Lemma \ref{calc6}.
\epf

\noindent{\bf Section \ref{Mod}.}\\ 

\begin{tag}\label{def2ctheta}
Definition of the 2-cell $\theta$ from Axiom \ref{cohmdobj3p}.
\end{tag}
\pf
The arrows
$$\xymatrix{
\A \otimes \M 
\ar[r]^-{L' \otimes 1}
&
(\unc \otimes \A) \otimes \M
\ar[r]^-{A'}
&
\unc \otimes (\A \otimes \M)}$$
is strict and, as shown below, it has the same image by $\Res$
as the arrow $L': \A \otimes \M \rightarrow \unc \otimes (\A \otimes \M)$.
The 2-cell $\theta$ corresponds then via the adjunction
\ref{punitens} to the identity 2-cell.\\
The image by $\Res$ of $L'_{\A}$ is 
$$
\xymatrix{
\A 
\ar[r]^-{\eta}
&
[\M, \A  \M]
\ar[r]^-{[1,\eta^*]}
&
[\M,[\unc,\unc  (\A \M)]]
\ar[r]^{[1,\ev_{\gen}]}
&
[\M, \unc (\A \M)]
.}
$$
The arrow $A' \circ L'\otimes 1$ above as 
image by $\Res$ that rewrites
{\tiny
\begin{tabbing}
\small
\=1.  \=$\xymatrix{
\A \ar[r]^-{L'}
& \unc \A
\ar[r]^-{\Res(A')}
&
[\M,\unc(\A \M)]
}$\\
\>2.  \>$\xymatrix{
\A \ar[r]^-{\eta^*}
&
[\unc, \unc \A]
\ar[r]^-{\ev_{\gen}}
& 
\unc \A
\ar[r]^-{\Res(A')}
&
[\M,\unc(\A \M)]
}$\\
\>3. \>$\xymatrix@C=3pc{
\A 
\ar[r]^-{\eta^*}
&
[\unc, \unc \A]
\ar[r]^-{[1,\Res(A')]}
& 
[\unc, [\M, \unc(\A \M)]]
\ar[r]^-{\ev_{\gen}}
&
[\M,\unc(\A \M)]
}$\\
\>4.  \>$\xymatrix@C=4pc{
\A 
\ar[r]^-{{\Res(\Res(A'))}^*}
& 
[\unc, [\M, \unc(\A \M)]]
\ar[r]^-{\ev_{\gen}}
&
[\M,\unc(\A \M)]
}$\\
\>5.  \>$\xymatrix@C=3pc{
\A \ar[r]^-{\eta}
&
[\M, \A \M]
\ar[r]^-{[-,\unc(\A\M)]}
& 
[[\A\M,\unc(\A\M)], [\M, \unc(\A \M)]]
\ar[r]^-{[\eta,1]}
&
[\unc,[\M,\unc(\A \M)]]
\ar[r]^-{\ev_{\gen}}
&
[\M,\unc (\A \M)]
}$\\
\>6.  \>$\xymatrix@C=3pc{
\A \ar[r]^-{\eta}
&
[\M, \A \M]
\ar[r]^-{[1,\eta^*]}
&
[\M,[\unc,\unc(\A\M)]]
\ar[r]^{\Dual}
&
[\unc,[\M,\unc(\A \M)]]
\ar[r]^-{\ev_{\gen}}
&
[\M,\unc (\A \M)]
}$\\
\>7.  \>$\xymatrix@C=3pc{
\A \ar[r]^-{\eta}
&
[\M, \A \M]
\ar[r]^-{[1,\eta^*]}
&
[\M,[\unc,\unc(\A\M)]]
\ar[r]^{[1,\ev_{\gen}]}
&
[\M,\unc (\A \M)]
}$
\end{tabbing}
}
where in the above derivation the equalities
between arrows hold for the following reasons:\\
- 4. and 5. by definition of $A'$ since its image 
by $\Res \circ \Res$ is in this case\\
{\tiny $\xymatrix{ 
\unc 
\ar[r]^-{\eta}
&
[\A\M,\unc(\A\M)]
\ar[r]^-{[\M,-]} 
&
[[\M,\A\M],[\M,\unc(\A\M)]]
\ar[r]^-{[\eta,1]}
&
[\A,[\M,\unc(\A\M)]]
}$ 
}
that has dual\\
{\tiny
$\xymatrix{
\A 
\ar[r]^-{\eta}
&
[\M,\A\M]
\ar[r]^-{[-,\unc(\A\M)]}
&
[[\A \M,\unc(\A\M)], [\M,\unc(\A\M)]]
\ar[r]^-{[\eta,1]}
&
[\unc, [\M,\unc(\A \M)]]
}$;
}\\
- 5. and 6. by Lemma \cite{Sch08}-\tbcun;\\
- 6. and 7. by Lemma \cite{Sch08}-\evdual.
\epf

\begin{tag}\label{pfptm1}
Proof of \ref{ptm1}.
\end{tag}
\pf
The first of the 2-cells of Axiom \ref{cohmdobj1p}
can be decomposed as the composite $\Xi_2 \circ \Xi_1$
where $\Xi_1$ is 
$$
\xymatrix{
((\A\A)\A)\M
\ar[r]^{A'}
&
(\A\A)(\A\M)
\ar[r]^{1 \otimes \mult}
&
(\A\A)\M
\ar@{=>}[r]^-{\assmd}
&
\M
}
$$
and
$\Xi_2$ is
$$
\xymatrix{
((\A \A)\A)\M
\ar[r]^-{(\comp \otimes 1) \otimes 1}
&
(\A \A) \M
\ar@{=>}[r]^-{\assmd}
&
\M.
}
$$
The 2-cell of \ref{ptm1} decomposes as
$\Xi_4 \circ \Xi_3$
where 
$\Xi_3$ is
{\small
$$\xymatrix@C=3pc{
(\A \otimes \A) \otimes \A 
\ar@{=>}[r]^-{\assm \otimes 1} 
&  
[\M,\M] \otimes \A 
\ar[r]^-{1 \otimes \multh} 
&
[\M,\M] \otimes [\M,\M]
\ar[r]^-{\comp}
&
[\M,\M]
}$$
}
and $\Xi_4$ is
$$
\xymatrix{
(\A \A) \A
\ar[r]^{\comp \otimes 1}
&
\A \A
\ar@{=>}[r]^-{\assm} 
&
[\M,\M].
}
$$
The 2-cell $\Xi_3$ has a strict domain 
and an easy computation gives that its 
image by $\Res$ is
$$
\xymatrix@C=3pc{
\A \otimes \A
\ar@{=>}[r]^-{\assm}
&
[\M,\M]
\ar[r]^-{[\M,-]}
&
[[\M,\M],[\M,\M]]]
\ar[r]^-{[\multh,1]}
&
[\A,[\A,[\M,\M]]]
.}
$$
According to Lemma \ref{calc6} the 2-cell $\Xi_1$
has an image by $\Res$ with a strict domain 
and has image by $\Res \circ \Res$ the 2-cell
$$
\xymatrix{
\A \A
\ar@{=>}[r]^-{\Res(\assmd)}
&
[\M,\M]
\ar[r]^-{[\M,-]}
&
[[\M,\M],[\M,\M]]
\ar[r]^-{[\multh,1]}
&
[\A,[\M,\M]].
}
$$
Therefore the 2-cell $\Xi_1$ has image by $\Res$
the 2-cell $\Xi_3$.

An easy computation shows that the 2-cell $\Xi_2$ 
has image by $\Res$
$$\xymatrix{
(\A\A)\A
\ar[r]^-{\comp \otimes 1}
&
\A \A
\ar@{=>}[r]^-{\Res(\assmd)}
&
[\M,\M]
}$$
which is $\Xi_4$.
\epf

\begin{tag}\label{pfptm2}
Proof of \ref{ptm2}.
\end{tag}
\pf
The second 2-cells of Axiom \ref{cohmdobj1p}
can be decomposed as the composite $\Xi_3 \circ \Xi_2 \circ \Xi_1$
where 
$\Xi_1$ is 
$$
\xymatrix{
((\A\A)\A)\M
\ar[r]^-{A' \otimes 1}
&
(\A(\A\A))\M
\ar[r]^-{A'}
&
\A((\A\A)\M)
\ar@{=>}[r]^-{1 \otimes \assmd}
&
\A\M
\ar[r]^-{\mult}
&
\M
}
$$
$\Xi_2$ is
$$
\xymatrix{
((\A\A)\A)\M
\ar[r]^-{A' \otimes 1}
&
(\A(\A\A))\M
\ar[r]^-{(1 \otimes \comp) \otimes 1}
&
(\A\A)\M
\ar@{=>}[r]^-{\assmd}
&
\M
}
$$
and 
$\Xi_3$ is 
$$
\xymatrix{
((\A\A)\A)\M
\ar@{=>}[r]^-{\assp \otimes 1}
&
\A \M
\ar[r]^-{\mult}
&
\M.
}
$$

The 2-cell of \ref{ptm2} decomposes as
$\Xi_6 \circ \Xi_5 \circ \Xi_4$
where 
$\Xi_4$ is the 2-cell
$$
\xymatrix{
(\A\A)\A 
\ar[r]^-{A'}
&
\A (\A\A)
\ar@{=>}[r]^-{1 \otimes \assm}
&
\A \otimes [\M,\M]
\ar[r]^{\multh \otimes 1}
&
[\M,\M] \otimes [\M,\M]
\ar[r]^-{\comp}
&
[\M,\M],
}
$$ 
$\Xi_5$ is 
$$
\xymatrix{
(\A\A)\A
\ar[r]^-{A'}
&
\A (\A\A)
\ar[r]^-{1 \otimes \comp}
&
\A\A
\ar@{=>}[r]^-{\assm}
&
[\M,\M]
}
$$
and $\Xi_6$ is 
$$
\xymatrix{
(\A\A)\A \ar@{=>}[r]^-{\assp} 
&
\A 
\ar[r]^-{\multh}
& 
[\M,\M]. 
}
$$

The 2-cell $\Xi_1$ has image by $\Res$ 
the 2-cell 
$$\xymatrix{
(\A\A)\A \ar[r]^-{A'}
&
\A(\A\A)
\ar[r]^-{\Res(A')}
&
[\M,\A((\A\A)\M)]
\ar@{=>}[r]^-{[1,1 \otimes \assmd]}
&
[\M,\A \M]
\ar[r]^-{[1,\mult]}
&
[\M,\M]
}$$

The 2-cell 
$$\Xi_7 =
\xymatrix{
\A(\A\A)
\ar[r]^-{\Res(A')}
&
[\M, \A((\A\A)\M)]
\ar@{=>}[r]^-{[1,1 \otimes \assmd]}
&
[\M,\A \M]
\ar[r]^-{[1,\mult]}
&
[\M,\M]
}$$
has a strict domain and
according to Lemma \cite{Sch08}-\calcdeux
its image by $\Res$ is
$$\Xi_8 =
\xymatrix{
\A
\ar[r]^-{\multh}
&
[\M, \M]
\ar[r]^-{[ \M,-]}
&
[[\M,\M],[\M,\M]]
\ar[r]^-{[\assm,1]}
&
[\A \A ,[ \M,\M ]]
.}
$$ 
The 2-cell 
$$\Xi_9 =
\xymatrix{
\A (\A\A)
\ar@{=>}[r]^-{1 \otimes \assm}
&
\A \otimes [\M,\M]
\ar[r]^-{\multh \otimes 1}
&
[\M,\M] \otimes [\M,\M]
\ar[r]^-{\comp}
&
[\M,\M]
}
$$
has also a strict domain and 
its image by $\Res$ is $\Xi_8$,
therefore $\Xi_7 = \Xi_9$ and 
$\Res(\Xi_1) = \Xi_7 * A' = \Xi_9 * A' = \Xi_4$.\\

The image by $\Res$ of $\Xi_2$ is 
$$\xymatrix{
(\A\A)\A 
\ar[r]^{A'}
&
\A(\A\A)
\ar[r]^-{1 \otimes \comp}
&
\A \A
\ar[r]^-{\Res(\assmd)}
&
[\M,\M]
}$$
which is $\Xi_5$.\\

Eventually the image by $\Res$ of $\Xi_3$ is trivially $\Xi_6$.
\epf

To prove \ref{ptm4}, we shall use the following lemma.
\begin{lemma}\label{dualevn}
For any objects $\A$, $\B$, $\C$ and $\D$, any
2-cell $\xymatrix{ \C \ar@{=>}[r]^-{\tau} & [\D,\A]}$
of $\BS$
and any $c \in \C$, the 2-cells in $\BS$
$$\xymatrix{
[\A,\B]
\ar[r]^-{[\D,-]}
&
[[\D,\A],[\D,\B]]
\ar@{=>}[r]^-{[\tau,1]}
&
[\C,[\D,\B]]
\ar[r]^-{\ev_c}
&
[\D,\B]
}$$
and 
$$
\xymatrix{
[\A,\B]
\ar[r]^-{[\ev_c,1]}
&
[[\C,\A],\B]
\ar@{=>}[r]^-{[\tau^*,1]}
&
[\D,\B]
}
$$
are equal.
\end{lemma}
\pf
According to Lemma \cite{Sch08}-\evdual,
The first of the 2-cell rewrites
\begin{tabbing}
\=1. \=$\xymatrix{
[\A,\B]
\ar[r]^-{[\D,-]}
&
[[\D,\A],[\D,\B]]
\ar@{=>}[r]^-{[\tau,1]}
&
[\C,[\D,\B]]
\ar[r]^-{\Dual}
&
[\D,[\C,\B]]
\ar[r]^-{[1,\ev_c]}
&
[\D,\B]
}$\\
\> according to Lemma \cite{Sch08}-\tbctrois,\\
\>2. \>$\xymatrix{
[\A,\B]
\ar[r]^-{[\C,-]}
&
[[\C,\A],[\C,\B]]
\ar[r]^-{[\tau^*,1]}
&
[\D,[\C,\B]]
}
$
\end{tabbing}
and this last arrow has dual
$$
\xymatrix{
\D \ar@{=>}[r]^-{\tau^*} 
& 
[\C,\A]
\ar[r]^-{[-,\B]}
&
[[\A,\B],[\C,\B]]
\ar[r]^-{[1,\ev_c]}
&
[\A,\B],\B]
.}
$$
One the other hand the dual of the second 2-cell of the lemma
rewrites
\begin{tabbing}
\=1. \=$\xymatrix{
\D
\ar@{=>}[r]^-{\tau^*}
&
[\C,\A]
\ar[r]^-{{[\ev_c,1]}^*}
&
[[\A,\B],\B]
}
$\\
\>2. \>$\xymatrix{
\D
\ar@{=>}[r]^-{\tau^*}
&
[\C,\A]
\ar[r]^-{\ev_c}
&
\A
\ar[r]^-{q}
&
[[\A,\B],\B]
}
$\\
\>3. \>$\xymatrix{
\D
\ar@{=>}[r]^-{\tau^*}
&
[\C,\A]
\ar[r]^-{[-,\B]}
&
[[\A,\B],[\C,\B]]
\ar[r]^-{[1,\ev_c]}
&
[[\A,\B],\B]
}
$
\end{tabbing}
where in the above derivation 
arrows 2. and 3. are equal by Corollary \cite{Sch08}-\homevtroisdual.
\epf

\begin{tag}\label{pfptm4} Proof of \ref{ptm4}. \end{tag} 
\pf
The 2-cell from \ref{ptm4}
decomposes as $\zeta_2 \circ \zeta_1$
where 
$\zeta_1$ $=$
$$\xymatrix{
\A 
\ar[r]^-{R'}
&
\A \otimes \unc
\ar@{=>}[r]^-{1 \otimes \lmz}
&
\A \otimes [\M,\M]
\ar[r]^-{\multh \otimes 1}
&
[\M,\M] \otimes [\M,\M]
\ar[r]^-{\comp}
&
[\M,\M]
}$$
and
$\zeta_2$ $=$ 
$$
\xymatrix{
\A
\ar[r]^-{R'}
&
\A \otimes \unc
\ar[r]^-{1 \otimes \unit}
&
\A \otimes \A
\ar@{=>}[r]^-{\assm}
&
[\M,\M].
}
$$
We show below that the image by $\Res$
of the 2-cell $\Xi_2$ is $\zeta_1$ whereas
the image by $\Res$ of $\Xi_3$ is $\zeta_2$.\\

The image by $\Res$ of 
the 2-cell $\Xi_2$ rewrites
{\small
\begin{tabbing}
\=1. \=$\xymatrix{
\A \ar[r]^-{\eta} 
&
[\A,\A \M]
\ar@{=>}[r]^-{[\lmd,1]}
&
[\M,\A \otimes \M]
\ar[r]^-{[1,\mult]}
&
[\M,\M]
}$\\
\>2. \>$\xymatrix{
\A \ar[r]^-{\eta} 
&
[\M,\A \M]
\ar[r]^-{[1,\mult]}
&
[\M,\M]
\ar@{=>}[r]^-{[\lmd,1]}
&
[\M,\M]
}$\\
\>3. \>$\xymatrix{
\A \ar[r]^-{\multh} 
&
[\M,\M]
\ar@{=>}[r]^-{[\lmd,1]}
&
[\M,\M]
}$\\
\end{tabbing}
}
According to Lemma \ref{dualevn} and since
the 2-cell $\lmd$ is $\ev_{\gen} * {(\lmz)}$, 
the 2-cell 3. above is 
$$
\xymatrix{
\A
\ar[r]^-{\multh}
&
[\M,\M]
\ar[r]^-{[\M,-]}
&
[[\M,\M],[\M,\M]]
\ar@{=>}[r]^-{[\lmz,1]}
&
[\unc,[\M,\M]]
\ar[r]^-{\ev_{\gen}}
&
[\M,\M].
}
$$
This last 2-cell is actually $\zeta_1$. To check this
use Lemma \ref{calc1}
and the fact that the image by $\Res$ of the arrow
$$\xymatrix{
\A \otimes [\M,\M]
\ar[r]^-{\multh \otimes 1}
&
[\M,\M] \otimes [\M,\M]
\ar[r]^-{\comp}
&
[\M,\M]
}$$
is 
$$\xymatrix{
\A
\ar[r]^-{\multh}
&
[\M,\M]
\ar[r]^-{[\M,-]}
&
[[\M,\M],[\M,\M]].
}$$

The image by $\Res$ of the 2-cell $\Xi_3$ is
$\xymatrix{
\A
\ar[r]^-{R'}
&
\A \otimes \unc
\ar[r]^-{1 \otimes \unit}
&
\A \otimes \A
\ar@{=>}[r]^-{\Res(\assmd)}
&
[\M,\M]
}$
which is $\zeta_2$.\\
\epf

\begin{tag}\label{pfptm6}
Proof of \ref{ptm6}.
\end{tag}
\pf
The 2-cell $\Xi_2$ rewrites 
$$
\xymatrix{
\A \otimes \M
\ar[r]^-{\mult}
&
\M
\ar@{-}[rr]^{id}
\ar[d]_{\multh^*}
&
\ar@{=>}[d]^{\lmh}
&
\M
\\
& 
[\A,\M]
\ar[rr]_-{[\unit,1]}
& & 
[\unc,\M]
\ar[u]_{\ev_{\gen}}
}
$$
which image by $\Res$ is the 2-cell
$$
\xymatrix{
\A
\ar[r]^{\multh}
&
[\M,\M]
\ar@{-}[rr]^-{id}
\ar[d]_{[1,\multh^*]}
& 
\ar@{=>}[d]^{[1,\lmh]}
& 
[\M,\M]\\
& 
[\M,[\A,\M]]
\ar[rr]_-{[1,[\unit,1]]}
&
&
[\M,[\unc,\M]]
\ar[u]_{[1,\ev_{\gen}]}.
}
$$
On the other hand
the 2-cell
$$\xymatrix{
\A 
\ar[r]^-{L'}
& 
\unc \otimes \A
\ar[r]^-{1 \otimes \multh}
&
\unc \otimes [\M,\M]
\ar@{=>}[r]^-{\lmz \otimes 1}
&
[\M,\M] \otimes [\M,\M]
\ar[r]^-{\comp}
&
[\M,\M]
}$$
rewrites 
\begin{tabbing}
\=1. \=$\xymatrix{
\A 
\ar[r]^-{\multh}
& 
[\M,\M]
\ar[r]^-{L'}
&
\unc \otimes [\M,\M]
\ar@{=>}[r]^{\lmz \otimes 1}
&
[\M,\M] \otimes [\M,\M]
\ar[r]^-{\comp}
&
[\M,\M]
}$\\
\>2. \>$
\xymatrix{
\A
\ar[r]^-{\multh}
&
[\M,\M]
\ar[r]^-{[-,\M]}
&
[[\M,\M],[\M,\M]]
\ar@{=>}[r]^-{[\lmz,1]}
&
[\unc,[\M,\M]]
\ar[r]^-{\ev_{\gen}}
&
[\M,\M]
}
$\\
\>3. \>$
\xymatrix{
\A
\ar[r]^-{\multh}
&
[\M,\M]
\ar@{=>}[r]^-{[1,{\lmz}^*]}
&
[\M,[\unc,\M]]
\ar[r]^-{\Dual}
&
[\unc,[\M,\M]]
\ar[r]^-{\ev_{\gen}}
&
[\M,\M]
}
$\\
\>4. \>$
\xymatrix{
\A
\ar[r]^{\multh}
&
[\M,\M]
\ar@{=>}[r]^-{[1,{\lmz}^*]}
&
[\M,[\unc,\M]]
\ar[r]^-{[1,\ev_{\gen}]}
&
[\M,\M]
}
$\\
\>5. \>$
\xymatrix{
\A
\ar[r]^-{\multh}
&
[\M,\M]
\ar@{=>}[r]^-{[1,\lmh]}
&
[\M,\M].
}
$
\end{tabbing}
In the above derivation the equality between
arrows hold for the following reasons:\\
- 1. and 2. by Lemma \ref{calc6},\\
- 2. and 3. by Lemma \cite{Sch08}-\tbcun (which 
can be improved to take 2-cells into account),\\
- 3. and 4. by Lemma \cite{Sch08}-\evdual.\\

By definition of the 2-cell $\theta$, the 2-cell $\Xi_3$ has image 
by $\Res$ an identity 2-cell.
Eventually it is rather straightforward that 2-cell $\Xi_4$ has image by $\Res$
$$
\xymatrix{
\A 
\ar[r]^-{L'}
&
\unc \otimes \A
\ar[r]^-{\unit \otimes 1}
&
\A \otimes \A
\ar@{=>}[r]^-{\assm}
&
[\M,\M]
}
$$
since the 2-cell $\assm$ is $\Res(\assmd)$.
\epf

\begin{tag}\label{pfptm7}
PROOF of \ref{ptm7}
\end{tag}
\pf
The 2-cell
$$\xymatrix{
(\A \A) \M
\ar[r]^-{1 \otimes \mor}
&
\A \A \N
\ar@{=>}[r]^-{\assmd}
& 
\N 
}$$
has a strict domain
$$\xymatrix{ 
\ar[r]^-{1 \otimes \mor}
& 
\ar[r]^-{A'}
&
\ar[r]^-{1 \otimes \psi}
&
\ar[r]^-{\psi}
&
}$$
and has image by $\Res$
$$
\xymatrix{
\A \A 
\ar@{=>}[r]^-{\assm}
&
[\N,\N]
\ar[r]^-{[\mor,1]}
&
[\M,\N]
}
$$
which is image by $\Res$ is $\Xi_1$.
The 2-cell
$$\xymatrix{
(\A \A) \M
\ar[r]^-{\comp \otimes 1}
&
\A \M
\ar@{=>}[r]^-{\deltad}
&
\N
}$$
has image by $\Res$ the 2-cell
$$
\xymatrix{
\A \A 
\ar[r]^-{\comp}
&
\A
\ar@{=>}[r]^-{\deltah}
&
[\M,\N]
}$$
which image by $\Res$ is $\Xi_2$.
\epf

\begin{tag}\label{pfptm8}
Proof of \ref{ptm8}.
\end{tag}
\pf
According to Lemma \ref{calc1}, the 2-cell 
$$
\xymatrix{
(\A\A)\M
\ar[r]^-{A'}
&
\A(\A\M)
\ar@{=>}[r]^-{1 \otimes \deltad}
&
\A \N
\ar[r]^-{\psi}
&
\N
}
$$
has image by $\Res \circ \Res$ the 2-cell $\Xi_4$.
According to Lemma \ref{calc1}, the 2-cell
$$
\xymatrix{
(\A \A) \M 
\ar[r]^-{A'}
&
\A (\A \M) 
\ar[r]^-{1 \otimes \mult}
&
\A \M 
\ar@{=>}[r]^-{\deltad}
&
\N
}
$$
has image by $\Res \circ \Res$ the 2-cell $\Xi_6$.
Eventually the 2-cell 
$$
\xymatrix{
(\A \A) \M
\ar@{=>}[r]^-{\beta}
&
\M
\ar[r]^-{\sigma}
&
\N
}
$$
has image by $\Res \circ \Res$ the 2-cell $\Xi_7$.
\epf

\begin{tag}\label{pfptm9}
Proof of \ref{ptm9}.
\end{tag}
\pf 
The 2-cell 
$$\lhpp_{\mor}: 
\xymatrix{\mor \rightarrow [1,\mor] \circ \vv:  \unc \rightarrow [\M,\N]}$$
has image by $\ev_{\gen}$ an identity 2-cell and 
the 2-cell
$$ 
\xymatrix{
\unc 
\ar@{=>}[r]^-{\lmz}
&
[\M,\M]
\ar[r]^-{[1,\mor]}
& 
[\M,\N]
}$$
has image by $\ev_{\gen}$ the 2-cell
$$
\xymatrix{
\M
\ar@{=>}[r]^-{\lmd}
&
\M
\ar[r]^-{\mor}
& 
\N
}
$$
\epf

\begin{tag}\label{pfptm10} Proof of \ref{ptm10}. \end{tag}
\pf
The arrow $\mor: \unc \rightarrow [\M,\N]$ is strict
has image by the functor $\ev_{\gen}$ 
the arrow $\mor: \M \rightarrow \N$. 
The 2-cell 
$$\xymatrix{
\unc
\ar@{=>}[r]^-{\lmz}
&
[\N,\N]
\ar[r]^-{[\mor,1]}
&
[\M,\N]
}$$
has image by $\ev_{\gen}$ 
the 2-cell
$$
\xymatrix{
\M
\ar[r]^-{\mor}
&
\N
\ar@{=>}[r]^-{\lmd}
&
\N
.}
$$
The 2-cell
$$
\xymatrix{
\unc 
\ar[r]^-{\unit}
&
\A
\ar@{=>}[r]^-{\deltah}
&
[\M,\N]
}
$$
has image by $\ev_{\gen} = SMC(1,\ev_{\gen}) \circ \Dual$
the 2-cell
$$
\xymatrix{
\M
\ar@{=>}[r]^-{{\Res(\delta)}^*}
&
[\A,\N]
\ar[r]^-{[\unit,1]}
&
[\unc,\N]
\ar[r]^-{\ev_{\gen}}
&
\N
}
$$
which is according to Lemma \ref{calc6}
$$
\xymatrix{
\M
\ar[r]^-{L'}
&
\unc \otimes \M
\ar[r]^-{\unit \otimes 1}
& 
\A \otimes \M
\ar@{=>}[r]^-{\delta}
&
\M
}
$$
\epf

\end{section}

\end{document}